%% file: Paper_3.tex
\documentclass[final,leqno]{siamltex}
\usepackage{fullpage}
\usepackage{graphicx}
\usepackage[T1]{fontenc}
\usepackage{psfrag}
\usepackage{amsmath}
\usepackage{amssymb}
\usepackage{color}
\usepackage{times}
\usepackage{afterpage}
\usepackage{mathtools}
\usepackage{bbm}
\usepackage{mathdots}
\usepackage{diagbox}

\newtheorem{remark}{Remark}[section]

\DeclareMathOperator{\sech}{sech}
\DeclareMathOperator{\cosech}{cosech}

\usepackage[letterpaper=true,colorlinks=true,
linkcolor=red,filecolor=green,citecolor=blue,
pdfpagemode=None]{hyperref}

\usepackage{xcolor}
\colorlet{linkequation}{cyan}
\newcommand*{\SavedEqref}{}
\let\SavedEqref\eqref
\renewcommand*{\eqref}[1]{%
  \begingroup
    \hypersetup{
      linkcolor=linkequation,
      linkbordercolor=linkequation,
    }%
    \SavedEqref{#1}%
  \endgroup
}

\makeatletter
\newcommand{\clonelabel}[2]{\@bsphack
  \expandafter\ifx\csname r@#2\endcsname\relax
  \else\protected@write\@auxout{}{\string\newlabel{#1}%
    {\csname r@#2\endcsname}}%
  \fi
  \expandafter\ifx\csname r@#2@cref\endcsname\relax
  \else\protected@write\@auxout{}{\string\newlabel{#1@cref}%
    {\csname r@#2@cref\endcsname}}%
  \fi
  \@esphack}
\makeatother

\begin{document}

\title{Optimal Control of a Subdiffusion model using Neumann-Neumann and Dirichlet-Neumann Waveform Relaxation
  Algorithms}

\author{Soura Sana\footnotemark[1]\and Bankim C. Mandal\footnotemark[2]}

\maketitle

\begin{abstract}
This article investigates the  convergence patterns exhibited by Neumann-Neumann as well as Dirichlet-Neumann waveform relaxation algorithms when applied to distributed optimal control problems subject to sub-diffusive and diffusive partial differential equation (PDE) constraints. Focusing on regular 1D domains with multiple subdomains, we examine the influence of varying values of the diffusion coefficient on algorithms' convergence. By adopting semi-discrete approach, we discretize in time while maintaining continuity in space. Our simple but elegant approach of a both-sided graded mesh in time simplifies analysis and enhances convergence compared to traditional mesh grading used in time-fractional PDEs.
\end{abstract}



\renewcommand{\thefootnote}{\fnsymbol{footnote}}
\footnotetext[1]{School of Basic Sciences, IIT Bhubaneswar, India ({\tt ss87@iitbbs.ac.in}).}
\footnotetext[2]{School of Basic Sciences, IIT Bhubaneswar, India ({\tt bmandal@iitbbs.ac.in}).}

\input{Introduction}

\input{Model_Problem_and_Algorithms}

\input{Lemmas}

\input{Mesh}

\input{Proof_DNWR}

\input{Proof_NNWR}

\input{Numerical}

\input{Conclusion}


\bibliographystyle{siam}
\bibliography{paper}

\end{document}

%% file: Introduction.tex
\section{Introduction}
This article is inspired by Gabriele Ciaramella, Laurence Halpern, and Luca Mechelli's work titled "Convergence Analysis and Optimization of a Robin Schwarz Waveform Relaxation Method for Time-Periodic Parabolic Optimal Control Problems" \cite{ciaramella2024convergence}. In their study, they addressed time-periodic reaction-diffusion partial differential equations as constraints for distributed optimal control problems and employed the Optimized Schwarz waveform relaxation algorithm for solution. In our work, we extend their approach with some modification by employing both  Neumann-Neumann and Dirichlet-Neumann waveform relaxation (NNWR, DNWR) algorithms to solve distributed optimal control problems subject to sub-diffusion and diffusion partial differential equation constraints. Furthermore, we relax the time-periodic assumption.

Fractional calculus has become pivotal in modeling complex systems in contemporary research \cite{rossikhin2010application,povstenko2020fractional,qin2017multi,quiroga2015adjoint}. These applications are particularly significant where the current state of evolution depends not only on the immediate previous state but on the entire history, leading to computationally expensive and memory-intensive numerical procedures \cite{qian2017certified}. To tackle these challenges, distributing computations across multiple processors becomes essential. Domain Decomposition (DD) methods \cite{schwarz1870ueber,lions1990schwarz,dryja1990some,cai1999restricted} are highly efficient in this regard. DD methods can be applied in discrete as well as continuous settings. When time evolution is involved, waveform relaxation techniques coupled with DD is introduced \cite{gander1998space,giladi2002space}. DNWR, NNWR, OSWR, and SWR are examples of substructuring waveform relaxation techniques\cite{gander2021dirichlet,etna_vol45_pp424-456,mandal2017neumann,gander2007optimized}. The effectiveness of DNWR and NNWR algorithms for sub-diffusion and diffusion-wave equations has been extensively investigated in \cite{sana2023dirichlet,sana2023dirichlet_a}.

The history of optimal control problems spans over three centuries. Modern optimal control theory, as we know it, was developed by L. S. Pontryagin \cite{boltyanskiy1962mathematical}. J. L. Lions played a substantial role in advancing the theory of optimal control of PDEs \cite{lions1971optimal}. For a comprehensive historical survey of optimal control, refer to the works of Sussmann et al. \cite{sussmann1997300} and Fernandez \cite{fernandez2003control}. Solutions to optimal control problems generally fall into two categories: direct and indirect methods. In the direct method, also known as "discretize then optimize," the continuous objective functional and constraints are first converted into discrete algebraic functions and equations, which are then solved using nonlinear programming techniques. In contrast, the indirect method, also known as "optimize then discretize," involves deriving optimality conditions, leading to a coupled system of PDEs, which are subsequently discretized using numerical techniques. If the control variable acts on the boundary, it is termed as boundary control, whereas if it affects the PDE as a forcing condition, the control problem is referred to as distributed control. In this article, we employ the indirect method to solve distributed optimal control problem. 
 
The study of optimal control in continuous cases subject to sub-diffusive PDEs was pioneered by Mophou et.al. \cite{mophou2011optimal_a,mophou2011optimal}. The main themes of their studies include the existence and uniqueness of solutions to control problems and the derivation of coupled PDEs. We apply DNWR and NNWR algorithms to solve these coupled systems.  DD methods have proven effective in handling PDE-constrained optimization problems \cite{bensoussan1973methode,benamou1996domain,benamou1997domain}.

This article is structured in the following manner: In Section \ref{model}, we first consider the distributed optimal control subject to sub-diffusion and diffusion constraints. Then, we use first-order optimality conditions to derive the state and adjoint equations. Sections \ref{sec_DNWR} through \ref{sec_NNWR} present the derivation of the DNWR and NNWR algorithms for coupled PDE systems. In Section \ref{Section4}, we establish some lemmas required for the  primary analysis of convergence. Section \ref{Section_mesh} introduces a both-sided graded mesh. Sections \ref{sectionDNWR} and \ref{NNWR_con} provide the main convergence results for the DNWR and NNWR, respectively.  Finally, in Section \ref{Numerical_ex}, we numerically verify all analytical results.

%% file: Model_Problem_and_Algorithms.tex
\section{Optimal Control of a subdiffusion model} \label{model}
Let, $\Omega$ be the open bounded subset in $\mathbb{R}^d$ with boundary $\partial\Omega$, and time window $[0,T]$ with $T > 0$. Let the quadratic cost functional be 
\begin{equation}\label{model_problem_1}
	J(y,u) := \frac{1}{2}\|y - y_Q\|^2_{L^2((0,T)\times \Omega)} + \frac{\sigma}{2}\|u\|^2_{L^2((0,T)\times \Omega)},
\end{equation}
where $u(\boldsymbol{x},t)$ is the control variable and $y_Q(\boldsymbol{x},t)$ be the target state. $\sigma$ be the regularization parameter, which mainly use to penalize the control variable. Our objective is to achieve the target state $y_Q(\boldsymbol{x},t)$ using the control $u(\boldsymbol{x},t)$, so that the cost functional $J(y,u)$ will be minimum. The state variable $y(\boldsymbol{x},t)$ is governed by the fractional diffusion equation of order $\alpha > 0$ given as follows:
\begin{equation}\label{model_problem_2}
	\begin{cases}
		D^{\alpha}_{RL}y(t,\boldsymbol{x}) = \nabla\cdot\left(\kappa(t,\boldsymbol{x})\nabla y(t,\boldsymbol{x})\right)+f(t,\boldsymbol{x}) + u(t,\boldsymbol{x}), &  \textrm{in}\; (0,T)\times\Omega,\\
		y(t,\boldsymbol{x}) = g(t,\boldsymbol{x}), & \textrm{on}\; (0,T)\times\partial\Omega,\\
		I^{1-\alpha}_{RL}y(0^{+},\boldsymbol{x}) = y_{0}(\boldsymbol{x}), & \textrm{in}\; \Omega.
	\end{cases}
\end{equation}
Here, $\kappa(t,\boldsymbol{x})>0$ is the dimensionless diffusion coefficient and $D^{\beta}_{RL}$ be the left sided  Riemann-Liouville fractional derivative  \cite{podlubny1998fractional} defined for any fractional order $\beta$, $n-1<\beta<n$, $n \in \mathbb{N}$  as follows:
\begin{equation*}
 D^{\beta}_{RL}f(t) := \frac{1}{\Gamma(n-\beta)}\frac{d}{dt^n}\int_0^t(t-\tau)^{n-\beta-1}f(\tau)d\tau.
\end{equation*}
And $I^{\beta}_{RL}$ represents the Riemann-Liouville fractional integral of order $\beta > 0$ on a continuous function $f(t)$ defined as follows:
\begin{equation*}
	I^{\beta}_{RL}f(t) := \frac{1}{\Gamma(\beta)}\int_0^t(t-\tau)^{\beta-1}f(\tau)d\tau.
\end{equation*}
For $0 < \alpha <1$ and $f,u \in L^2((0,T)\times \Omega)$ the problem \eqref{model_problem_2} has a unique solution $y \in L^2((0,T);H^2(\Omega) \cap H^1_0(\Omega))$ for zero Dirichlet boundary condition, see ~\cite{mophou2011optimal}. And the objective functional $J$ is convex, so unique infimum exist. If $y \in L^2((0,T)\times \Omega)$ be the solution of  \eqref{model_problem_1} then using the Euler-Lagrange optimality conditions we have adjoint variable $p \in L^2((0,T)\times \Omega)$ such that $(u,y,p)$ satisfy the optimality system ~\cite{mophou2011optimal_a}:
\begin{equation}\label{model_problem_3}
	\begin{cases}
		D^{\alpha}_{RL}y(t,\boldsymbol{x}) = \nabla\cdot\left(\kappa(t,\boldsymbol{x})\nabla y(t,\boldsymbol{x})\right)+f(\boldsymbol{x},t) + u(t,\boldsymbol{x}), & \textrm{in}\; (0,T)\times\Omega,\\
		y(t,\boldsymbol{x}) = g(t,\boldsymbol{x}), & \textrm{on}\; (0,T)\times\partial\Omega,\\
		I^{1-\alpha}_{RL}y(0^{+},\boldsymbol{x}) = y_{0}(\boldsymbol{x}), & \textrm{in}\; \Omega.
	\end{cases}
\end{equation}
\begin{equation}\label{model_problem_4}
	\begin{cases}
		_{T}D^{\alpha}_{t}p(t,\boldsymbol{x}) = \nabla\cdot\left(\kappa(t,\boldsymbol{x})\nabla p(t,\boldsymbol{x})\right)+ y(t,\boldsymbol{x}) - y_Q(t,\boldsymbol{x}), & \textrm{in}\; (0,T)\times\Omega,\\
		p(t,\boldsymbol{x}) = 0, & \textrm{on}\; (0,T)\times\partial\Omega,\\
		p(T,\boldsymbol{x}) = 0, & \textrm{in}\; \Omega.
	\end{cases}
\end{equation}
with the condition $\sigma u = - p$. $_{T}D^{\beta}_{t}$ be the right sided Caputo fractional derivative \cite{podlubny1998fractional} of order $\beta$, $n-1<\alpha<n$, $n \in \mathbb{N}$  defined as follows:
\begin{equation*}
	_{T}D^{\beta}_{t}f(t) := -\frac{1}{\Gamma(n-\beta)}\int_t^T(\tau - t)^{n-\beta-1}f^{(n)}(\tau)d\tau.
\end{equation*}

\section{The Dirichlet-Neumann Waveform Relaxation algorithm} \label{sec_DNWR}
The DNWR technique \cite{gander2021dirichlet} is an iterative algorithm that operates semi-parallelly, except in the scenario of two sub-domains. It merges the domain decomposition method for spatial resolution with the waveform relaxation method for temporal resolution. Suppose we have two non-overlapping sub-domains, $\Omega_{1}$ and $\Omega_{2}$, within $\Omega$, with their interface denoted as $\Gamma := \partial\Omega_{1}\cap\partial\Omega_{2}$.
The solution pair $(y,p)$ of state and adjoint equations and diffusion parameter $\kappa(t,\boldsymbol{x})$ in \eqref{model_problem_3}-\eqref{model_problem_4} are restricted to $(y_j,p_j)$ and constant $\kappa_j$ on $\Omega_j$, where $j=1,2$. Let $\boldsymbol{n}_{j}$ denote the unit outward normal vector on $\Gamma$ corresponding to $\Omega{j}$.
A starting estimate $\xi^{(0)}(t,\boldsymbol{x})$ for the state equation and $\zeta^{(0)}(t,\boldsymbol{x})$ for the adjoint equation are chosen across the interface $(0,T)\times\Gamma$. The Dirichlet and Neumann steps for state problems yield iterative solutions $y_j^{(k)}, \text{ where } j = 1,2, \text{ and } k = 1,2,\ldots$, given by

\begin{equation}\label{DNWR_state}
	\begin{array}{rcll}
	\begin{cases}
	D^{\alpha}_{RL}y_{1}^{(k)} = \nabla\cdot\left(\kappa_1\nabla y_{1}^{(k)}\right)+f - \frac{p_1^{(k)}}{\sigma}, & \textrm{in}\; \Omega_{1},\\
	I^{1-\alpha}_{RL}y_{1}^{(k)}(0,\boldsymbol{x}) = y_{0}(\boldsymbol{x}), & \textrm{in}\; \Omega_{1},\\
	y_{1}^{(k)} = \xi^{(k-1)}, & \textrm{on}\; \Gamma,\\
	y_{1}^{(k)}=g, & \textrm{on}\; \partial\Omega_{1}\setminus\Gamma,
	\end{cases}
	\end{array}\
	\begin{array}{rcll}
	\begin{cases}
	D^{\alpha}_{RL}y_{2}^{(k)} = \nabla\cdot\left(\kappa_2\nabla y_{2}^{(k)}\right)+f - \frac{p_2^{(k)}}{\sigma}, & \textrm{in}\; \Omega_{2},\\
	I^{1-\alpha}_{RL}y_{2}^{(k)}(0,\boldsymbol{x}) = y_{0}(\boldsymbol{x}), & \textrm{in}\; \Omega_{2},\\
	\kappa_2\partial_{\boldsymbol{n}_{2}} y_{2}^{(k)}  =  -\kappa_1\partial_{\boldsymbol{n}_{1}} y_{1}^{(k)}, & \textrm{on}\; \Gamma,\\
	y_{2}^{(k)}=g, & \textrm{on}\;\partial\Omega_{2}\setminus\Gamma,
	\end{cases}
	\end{array}
\end{equation}
and for the adjoint problems:
\begin{equation}\label{DNWR_adjoint}
		\arraycolsep0.008em
		\begin{array}{rcll}
			\begin{cases}
				_{T}D_{t}^{\alpha}p_{1}^{(k)} = \nabla\cdot\left(\kappa_1\nabla p_{1}^{(k)}\right) + y_1^{(k)} - y_Q|_1 , & \textrm{in}\; \Omega_{1},\\
				p_{1}^{(k)}(T,\boldsymbol{x}) = 0, & \textrm{in}\; \Omega_{1},\\
				p_{1}^{(k)} = \zeta^{(k-1)}, & \textrm{on}\; \Gamma,\\
				p_{1}^{(k)}=0, & \textrm{on}\; \partial\Omega_{1}\setminus\Gamma,
			\end{cases}
		\end{array}\
		\begin{array}{rcll}
			\begin{cases}
				_{T}D_{t}^{\alpha}p_{2}^{(k)} = \nabla\cdot\left(\kappa_2\nabla p_{2}^{(k)}\right)+ y_2^{(k)} - y_Q|_2, & \textrm{in}\; \Omega_{2},\\
				p_{2}^{(k)}(T,\boldsymbol{x}) = 0, & \textrm{in}\; \Omega_{2},\\
				\kappa_2\partial_{\boldsymbol{n}_{2}} p_{2}^{(k)}  =  -\kappa_1\partial_{\boldsymbol{n}_{1}} p_{1}^{(k)}, & \textrm{on}\; \Gamma,\\
				p_{2}^{(k)}=0, & \textrm{on}\;\partial\Omega_{2}\setminus\Gamma.
			\end{cases}
		\end{array}
\end{equation}
Subsequently, modify the interface after incorporating the relaxation parameters $\theta, \phi \in (0,1)$, by
\begin{align}\label{DNWR2}
  \xi^{(k)}(t,\boldsymbol{x})&=\theta y_{2}^{(k)}\left|_{(0,T)\times\Gamma}\right.+(1-\theta)\xi^{(k-1)}(t,\boldsymbol{x}),\\
  \zeta^{(k)}(t,\boldsymbol{x})&=\phi p_{2}^{(k)}\left|_{(0,T)\times\Gamma}\right.+(1-\phi)\zeta^{(k-1)}(t,\boldsymbol{x}).
\end{align}
The primary goal of our study is to examine the convergence to zero of the update part errors $\omega^{(k-1)}(t,\boldsymbol{x}) := y|_{(0,T)\times\Gamma} - \xi^{(k-1)}(t,\boldsymbol{x})$ and $\upsilon^{(k-1)}(t,\boldsymbol{x}) := p|_{(0,T)\times\Gamma} - \zeta^{(k-1)}(\boldsymbol{x},t)$, where $y|_{(0,T)\times\Gamma}$ and $p|_{(0,T)\times\Gamma}$ represent the exact interface solution. Linearity allows us to consider $f(t,\boldsymbol{x})=0$, $g(t,\boldsymbol{x})=0$, and $y_{0}(\boldsymbol{x})=0$ for the convergence of the error equations in \eqref{DNWR_state} and \eqref{DNWR_adjoint}, respectively. This analysis will be conducted in section \ref{sectionDNWR}.

\section{The Neumann-Neumann Waveform Relaxation algorithm}\label{sec_NNWR}
We are exploring the NNWR algorithm for a scenario involving multiple subdomains. In this context, the model problem \eqref{model_problem_3}- \eqref{model_problem_4} is defined over the spatial domain $\Omega$, which is partitioned into multiple disjoint subdomains $\Omega_{i}$, where $1\leq i\leq N$, with no cross point. We designate $\Gamma_{i,j}$ as the boundary where $\partial\Omega_{i}$ intersects with $\partial\Omega_{j}$ for $j=i-1,i+1$. Within each subdomain $\Omega_i$, we define unit outward normals denoted as $\boldsymbol{n}_{i,i-1}$ and $\boldsymbol{n}_{i,i+1}$ on the boundaries $\Gamma_{i,i-1}$ and $\Gamma_{i,i+1}$, respectively. Given initial guesses $\xi_{i,j}^{(0)}(\boldsymbol{x},t)$ for the state equation and $\zeta_{i,j}^{(0)}(\boldsymbol{x},t)$ for the adjoint equation along the interfaces $(0,T)\times \Gamma_{i,j}$, the NNWR algorithm for optimal control problem can be articulated in the following manner:

\begin{equation}\label{NNWR_state}
  \begin{array}{rcll}
  \begin{cases}
   D^{\alpha}_{RL}y_{i}^{(k)} = \nabla\cdot\left(\kappa_i\nabla y_{i}^{(k)}\right) + f - \frac{p_i^{(k)}}{\sigma}, & \mbox{in $\Omega_{i}$},\\
    I^{1-\alpha}_{RL}y_{i}^{(k)}(0,\boldsymbol{x}) = y_{0}(\boldsymbol{x}), & \mbox{in $\Omega_{i}$},\\
    y_{i}^{(k)} = g, & \mbox{on $\partial\Omega_{i}\cap\partial\Omega$},\\
    y_{i}^{(k)} = \xi_{i,j}^{(k-1)}, & \mbox{on $\Gamma_{i,j}$},
   \end{cases}
  \end{array}
\begin{array}{rcll}
\begin{cases}
	D^{\alpha}_{RL}\psi_{i}^{(k)} = \nabla\cdot\left(\kappa_i\nabla \psi_{i}^{(k)}\right) - \frac{\chi_i^{(k)}}{\sigma} , & \mbox{in $\Omega_{i}$},\\
   I^{1-\alpha}_{RL}\psi_{i}^{(k)}(0,\boldsymbol{x}) = 0, & \mbox{in $\Omega_{i}$},\\
   \psi_{i}^{(k)} = 0, & \mbox{on $\partial\Omega_{i}\cap\partial\Omega$},\\
   \kappa_i\partial_{\boldsymbol{n}_{i,j}}\psi_{i}^{(k)} = \kappa_i\partial_{\boldsymbol{n}_{i,j}}y_{i}^{(k)}+\kappa_j\partial_{\boldsymbol{n}_{j,i}}y_{j}^{(k)}, & \mbox{on $\Gamma_{i,j}$}.
\end{cases}
\end{array}
\end{equation}
\begin{equation}\label{NNWR_adjoint}
	\begin{array}{rcll}
		\begin{cases}
			_{T}D^{\alpha}_{t}p_{i}^{(k)} = \nabla\cdot\left(\kappa_i\nabla p_{i}^{(k)}\right) + y_i^{(k)} - y_Q|_i, & \mbox{in $\Omega_{i}$},\\
			p_{i}^{(k)}(T,\boldsymbol{x}) = 0, & \mbox{in $\Omega_{i}$},\\
			p_{i}^{(k)} = 0, & \mbox{on $\partial\Omega_{i}\cap\partial\Omega$},\\
			p_{i}^{(k)} = \zeta_{i,j}^{(k-1)}, & \mbox{on $\Gamma_{i,j}$},
		\end{cases}
	\end{array}
	\begin{array}{rcll}
		\begin{cases}
			_{T}D^{\alpha}_{t}\chi_{i}^{(k)} = \nabla\cdot\left(\kappa_i\nabla \chi_{i}^{(k)}\right) + \psi_i^{(k)}, & \mbox{in $\Omega_{i}$},\\
			\chi_{i}^{(k)}(T,\boldsymbol{x}) = 0, & \mbox{in $\Omega_{i}$},\\
			\chi_{i}^{(k)} = 0, & \mbox{on $\partial\Omega_{i}\cap\partial\Omega$},\\
			\kappa_i\partial_{\boldsymbol{n}_{i,j}}\chi_{i}^{(k)} = \kappa_i\partial_{\boldsymbol{n}_{i,j}}p_{i}^{(k)}+\kappa_j\partial_{\boldsymbol{n}_{j,i}}p_{j}^{(k)}, & \mbox{on $\Gamma_{i,j}$}.
		\end{cases}
	\end{array}
\end{equation}
Subsequently, the update formula is employed to refresh the interface values to get
\begin{align}\label{NNWR2}
  \xi_{i,j}^{(k)}(t,\boldsymbol{x})&=\xi_{i,j}^{(k-1)}(t,\boldsymbol{x})-\theta_{i,j}
  \left( \psi_{i}^{(k)}\left|_{(0,T)\times\Gamma_{i,j}}\right.+\psi_{j}^{(k)}\left|_{(0,T)\times\Gamma_{i,j}}\right.\right),\\
  \zeta_{i,j}^{(k)}(t,\boldsymbol{x})&=\zeta_{i,j}^{(k-1)}(t,\boldsymbol{x})-\phi_{i,j}
  \left( \chi_{i}^{(k)}\left|_{(0,T)\times\Gamma_{i,j}}\right.+\chi_{j}^{(k)}\left|_{(0,T)\times\Gamma_{i,j}}\right.\right).
\end{align}
Our objective is now to investigate the convergence of the error equations derived from \eqref{NNWR_state}-\eqref{NNWR_adjoint}. This analysis involves considering zero forcing terms and initial conditions, namely $f(t,\boldsymbol{x}) = 0$ and $y_{0}(\boldsymbol{x}) = 0$, respectively, as well as imposing zero natural boundary conditions. However, we introduce error components on the artificial interfaces defined as $\omega_{i,j}^{(k-1)}(t,\boldsymbol{x}) := y_i|_{(0,T)\times\Gamma_{i,j}} - \xi_{i,j}^{(k-1)}(t,\boldsymbol{x})$ and  $\upsilon_{i,j}^{(k-1)}(t,\boldsymbol{x}) := p_i|_{(0,T)\times\Gamma_{i,j}} - \zeta_{i,j}^{(k-1)}(t,\boldsymbol{x})$ as $k \to \infty$.


%% file: Lemmas.tex
\section{Auxiliary Results}\label{Section4}
Now we need some auxiliary result which are necessary for the convergent of DNWR and NNWR algorithm.

\begin{lemma}[Theorem 2.10 \cite{mophou2011optimal}]
	\label{rl_eqi_c}
	The Riemann-Liouville time derivative $D^{\alpha}_{RL}y(t)$ with homogeneous initial condition $I^{1-\alpha}_{RL}y(0^{+}) = 0$ is equivalent to the Caputo fractional time derivative $_0D_t^{\alpha}y(t)$ with zero initial condition $y(0) = 0$.
\end{lemma}
\begin{proof}
	We know from the relation between Riemann-Liouville and Caputo fractional derivative that
	\begin{equation*}
		D^{\alpha}_{C}y(t) = D^{\alpha}_{RL}y(t) + \frac{t^{-\alpha}}{\Gamma(1-\alpha)}y(0).
	\end{equation*}
	After integrating both side with order $\alpha$ i.e. $I^{\alpha}$ we have $y(t) - y(0) = y(t) - \frac{t^{-\alpha}}{\Gamma(1-\alpha)}I^{(1-\alpha)}y(0)$. Thus $I^{(1-\alpha)}y(0) = 0$ implies $y(0) = 0$. Inversely from the same relation $y(0) = 0$ implies $D^{\alpha}_{C}y(t) = D^{\alpha}_{RL}y(t)$. Hence equivalent.
\hfill\end{proof}
\begin{lemma} \label{stable}
	Stability of an elliptic pde depends on the least absolute eigenvalue of the discretize spatial matrix. 
\end{lemma}
\begin{proof}
	The proof is straight forward, hence omitted.
\hfill\end{proof}
\begin{lemma} \label{matrix}
	The matrix $M = \begin{bmatrix}
		A & \frac{1}{\sigma}B\\
		-B & A
	\end{bmatrix}$ is similar to the matrix  $\bar{M} = \begin{bmatrix}
	A +\frac{i}{\sqrt{\sigma}}B & 0 \\
	0 &  A -\frac{i}{\sqrt{\sigma}}B 
	\end{bmatrix}$.
\end{lemma}
\begin{proof}
	Let $X = \begin{bmatrix}
		X_{11} & X_{12}\\
		X_{21} & X_{22}
	\end{bmatrix}$ be the block matrix of order $2n$. Now if there exist an invertible matrix $X$ such that $MX = X\bar{M}$ then we say that $M$ and $\bar{M}$ are similar. Therefore,
	\begin{align*}
		 \begin{bmatrix}
			A & \frac{1}{\sigma}B\\
			-B & A
		\end{bmatrix} \begin{bmatrix}
		X_{11} & X_{12}\\
		X_{21} & X_{22}
	\end{bmatrix} &= \begin{bmatrix}
	X_{11} & X_{12}\\
	X_{21} & X_{22}
	\end{bmatrix}  \begin{bmatrix}
	A +\frac{i}{\sqrt{\sigma}}B & 0 \\
	0 &  A -\frac{i}{\sqrt{\sigma}}B 
	\end{bmatrix} \\
	\begin{bmatrix}
		AX_{11} + \frac{1}{\sigma}BX_{21} & AX_{12} + \frac{1}{\sigma}BX_{22}\\
		AX_{21} - BX_{11} & AX_{22} - BX_{12}
	\end{bmatrix} &= \begin{bmatrix}
	X_{11}(A +\frac{i}{\sqrt{\sigma}}B) & X_{12}(A -\frac{i}{\sqrt{\sigma}}B)\\
	X_{21}(A +\frac{i}{\sqrt{\sigma}}B) & X_{22}(A -\frac{i}{\sqrt{\sigma}}B)
	\end{bmatrix}
	\end{align*}
	Now if we choose $X_{11} = I,  X_{12} = I, X_{21} = i\sqrt{\sigma}I, X_{22} = -i\sqrt{\sigma}I$, where $I$ is the identity matrix of order $n$ then left and right side will be equal. Hence $X = \begin{bmatrix}
	I & I\\
	i\sqrt{\sigma}I & -i\sqrt{\sigma}I 
	\end{bmatrix}$ which is invertiable so the matrix $M$ is similar to matrix $\bar{M}$ by the matrix by $X$.
\hfill\end{proof}

The above result conclude that the diagonalization of matrix $M$ is separately depends on the diagonalization of matrix block  $A +\frac{i}{\sqrt{\sigma}}B$ and $A -\frac{i}{\sqrt{\sigma}}B$. The next result gives the lower bound on the real part of eigenvalues.
\begin{lemma} \label{positive_eigenvalues}
	Let $M = A + iB$ , where $A,B$ are square matrices of order $n$ in $\mathbb{R}$ and $B = B^T$. If the matrix $(A + A^T)$ is positive definite then the real part of the eigenvalues of $M$ are strictly positive. 
\end{lemma}
\begin{proof}
	For this result we apply the general version of Lyapunov's theorem on stability \cite{ostrowski1962some}, which state that: Let $X$ be a complex matrix of order $n$ and $Y$ be a positive definite Hermitian matrix of order $n$. Then there exist a negative definite Hermitian matrix $Z$ for which $XZ + ZX^* = Y$ ($'*'$ imply the conjugate transpose) holds, if and only if real part of all the eigenvalues of $X$ are negative.
	
	Now if we choose $Z = -I$ i.e. negative identity matrix and $X = -M$ then  we have $Y = (A + iB)+(A^T - iB^T) = A + A^T$ as $B = B^T$. From assumption $(A + A^T)$ is positive definite so the real part of the eigenvalues of $-M$ are negative. Hence the real part of the eigenvalues of $M$ are strictly positive.
	\hfill\end{proof}
\begin{corollary} \label{first_derivative}
	Let $A$ is the first order backward difference matrix with step size $\Delta t$ for $n$ nodes and $B$ be the exchange matrix of same order, i.e.
	\begin{equation*}
		A = \begin{bmatrix}
			1 &  & & & \\
			-1 & 1 &  & & \\
			& \ddots & \ddots & & \\
			& & -1 & 1 &  \\
			& & & -1 & 1
		\end{bmatrix}, 
		\hspace{10mm}
		\textit{ and }
		\hspace{10mm}
		B = \begin{bmatrix}
			& & & & 1 \\
			& & & 1 & \\
			& & \iddots & & \\
			& 1 & & & \\
			1 & & & &
		\end{bmatrix}.
	\end{equation*} 
	Clearly $B$ be symmetric and $-(A + A^T)$ be the central difference matrix with the eigenvalues $\lambda_j = -\frac{4}{\Delta t^2} \sin^2(\frac{\pi}{2}\frac{j}{n+1}),$ $j = 1,2,\cdots n$. Therefore the eigenvalues of $A + A^T$ are positive, hence the matrix is positive definite. Therefore the real part of the eigenvalues of $M = A + iB$ are strictly positive. 
\end{corollary}
\begin{corollary} \label{frac_derivative}
	Let $D^{\alpha}_{apx}u^m$ be the L1 approximation of Caputo derivative of $D^{\alpha}_tu(t_m)$ on the temporal graded mesh $0 = t_0 <t_1 < \cdots < t_n = T$, \cite{stynes2021survey}. Then
	\begin{align*}
		D^{\alpha}_{apx}u^m &= \frac{1}{\Gamma(1-\alpha)} \sum_{j = 0}^{m-1} \int_{t_j}^{t_{j+1}}\frac{du}{ds}(t_m - s)^{-\alpha} ds \hspace{10mm}\textit{ for } m = 1,2,\cdots, M \\
		&= \frac{1}{\Gamma(2-\alpha)} \sum_{j = 0}^{m-1} \frac{u^{j+1} - u^{j}}{t_{j+1} - t_{j}} \left[(t_m - t_j)^{1-\alpha} - (t_m - t_{j+1})^{1-\alpha}\right].
	\end{align*}
	Choosing the initial condition $u^0 = 0$ and define 
	\begin{equation} \label{diagonal_element}
			d_{m,j} := \frac{(t_m - t_{m-j})^{1-\alpha} - (t_m - t_{m-j+1})^{1-\alpha}}{t_{m-j+1} - t_{m-j}} = \frac{1-\alpha}{t_{m-j+1} - t_{m-j}} \int_{t_{m-j}}^{t_{m-j+1}}(t_m - s)^{-\alpha} ds,
	\end{equation} we have
	\begin{align*}
		D^{\alpha}_{apx}u^m = \frac{1}{\Gamma(2-\alpha)} \left[d_{m,1} u^m - d_{m,m} u^0 + \sum_{j = 1}^{m-1} \left(d_{m,j+1} - d_{m,j}\right)u^{m-j}\right].
	\end{align*}
	Now $(t_m - s)^{-\alpha}$ is monotonic increasing function for $s < t_m$ and $0<\alpha<1$. Hence by using integral mean value theorem one get $0 < d_{m,j+1} < d_{m,j}$ for all $m$ and $0 \leq j \leq m$. Therefore the difference matrix corresponding to the L1 scheme takes the form
	\begin{equation}\label{cap_mat_1}
		A = \begin{bmatrix}
			d_{1,1} &  & & & \\
			d_{2,2} - d_{2,1} & d_{2,1} &  & & \\
			\vdots& \ddots & \ddots & & \\
			d_{n-1,n-1} - d_{n-1,n-2} & \cdots &  & d_{n-1,1} & \\
			d_{n,n} - d_{n,n-1} & \cdots &  & d_{n,2} - d_{n,1} & d_{n,1}
		\end{bmatrix}.
	\end{equation}
	Clearly the $m^{th}$ row sum of absolute elements of matrix $A$ except the diagonal element be $R_m := |(d_{m,2} - d_{m,1})| + \cdots + |(d_{m,m} - d_{m,m-1})| = d_{m,1} - d_{m,m} < d_{m,1}$. Hence by using Gershgorin disk theorem we can say that $A+A^T$ be the strictly positive definite matrix. Therefore using lemma \ref{positive_eigenvalues} we can conclude that the real part of the eigenvalues of matrix $M = A + iB$, where $B$ takes the form as in corollary \ref{first_derivative}, are strictly positive.
\end{corollary}
\begin{lemma} \label{eigenvalues}
	Let $M = A + B$ , where $M,A,B$ are square matrices of order $n$ in $\mathbb{C}$. Let $B$ is diagonalizable, i.e. $B = XDX^{-1}$. Then each eigenvalue of matrix $M$ is contained within at least one of the Gershgorin discs of matrix $X^{-1}MX$.
\end{lemma}
\begin{proof}
	We have $M = A + B$ and $B$ is diagonalizable by $X$. Therefore
	\begin{align*}
		M &= A + XDX^{-1}\\
		X^{-1}MX &= X^{-1}AX + D \\
				 &=: R + D.
	\end{align*}
 As $M$ and $X^{-1}MX$ are similar matrices. Hence using Gershgorin circle theorem we conclude that each eigenvalue of matrix $M$ is contained within at least one of the Gershgorin discs $\mathcal{D}(r_{ii} + d_{ii}, R_i)$, where  $D = (d_{ij})_{n \times n}$, $R = (r_{ij})_{n \times n}$ and $R_i = \sum_{j = 1, j \neq i}^{n} |r_{ij}|$.
\hfill\end{proof}

The next result gives the estimate on Gershgorin radius $R_i$ for specific situation.
\begin{lemma}\label{radius}
	Let $M = S + J$, where $S$ and $J$ are respectively the lower triangular matrix and exchange matrix of order $n$. Then Greshgorin radius can be easily found out from the relation \eqref{r1} for $n$ even and from the relation \eqref{r2} for $n$ odd.	
\end{lemma}
\begin{proof}
	Let $n$ be even  and define $m := \frac{n}{2}$. Lower triangular matrix $S = \begin{bmatrix}
		A_{m,m} & O_{m,m}\\
		B_{m,m} & C_{m,m}
	\end{bmatrix}$. $O_{m,m}$ be the zero matrix. $J$ be the exchange matrix which is one form of permutation matrix therefore it can be written in the cycle notation as $(1 \,n)(2 \,(n-1))...(\lceil \frac{n}{2} \rceil \, \lceil \frac{n+1}{2} \rceil)$. When $n$ is even there are $\frac{n}{2}$number of 2-cycle. Let $(\alpha \, \beta)$ is one of such 2-cycle. Then this 2-cycle contributes to the spectrum the eigenvectors $\tilde{x}_1$ and $\tilde{x}_2$ corresponding to the eigenvalues $1$ and $-1$ respectively, where both the eigenvectors $\tilde{x}_1$ and $\tilde{x}_2$ have only $2$ nonzero components at $\alpha$ and $\beta$ positions with value $\frac{1}{\sqrt{2}}, \, \frac{1}{\sqrt{2}}$ and $\frac{1}{\sqrt{2}}, \, -\frac{1}{\sqrt{2}}$, i.e. $\tilde{x}_1 = \frac{1}{\sqrt{2}}[0 \, \cdots 0 \, 1_{\alpha^{th}} \, 0 \cdots \, 0 \,  1_{\beta^{th}} \, 0 \cdots \, 0]$ and $\tilde{x}_2 = \frac{1}{\sqrt{2}}[0 \, \cdots 0 \, 1_{\alpha^{th}} \, 0 \cdots \, 0 \,  -1_{\beta^{th}} \, 0 \cdots \, 0]$.
	Therefore the eigenvector matrix of $J$ can be written as $X = \frac{1}{\sqrt{2}}\begin{bmatrix}
		I_{m,m} & J_{m,m}\\
		J_{m,m} & -I_{m,m}
	\end{bmatrix}$ corresponding to the eigenvalue matrix  $\Lambda_J = \begin{bmatrix}
	I_{m,m} & O_{m,m}\\
	O_{m,m} & -I_{m,m}
	\end{bmatrix}$.  As $J$ is also symmetric matrix so it is diagonalizable by orthogonal matrix. Hence $J = X\Lambda_J X^T$. Therefore using the lemma \ref{eigenvalues} we can find the Gershgorin radius from the matrix
	\begin{align} \label{r1}
		X^TSX &= \frac{1}{2}\begin{bmatrix}
			I_{m,m} & J_{m,m}\\
			J_{m,m} & -I_{m,m}
		\end{bmatrix}\begin{bmatrix}
		A_{m,m} & O_{m,m}\\
		B_{m,m} & C_{m,m}
	\end{bmatrix}\begin{bmatrix}
		I_{m,m} & J_{m,m}\\
		J_{m,m} & -I_{m,m}
	\end{bmatrix} \\
	&= \frac{1}{2}\begin{bmatrix}
		A_{m,m} + J_{m,m}B_{m,m} + J_{m,m}C_{m,m}J_{m,m} & A_{m,m}J_{m,m} - J_{m,m}C_{m,m} + J_{m,m}B_{m,m}J_{m,m}\\
		J_{m,m}A_{m,m} - C_{m,m}J_{m,m} - B_{m,m} & C_{m,m} - B_{m,m}J_{m,m} + J_{m,m}A_{m,m}J_{m,m}
	\end{bmatrix}. \nonumber
	\end{align}
	 Now when $n$ is odd there are $m := \frac{n-1}{2}$number of 2-cycle and only one 1-cycle. Let $(\gamma)$ be such an 1-cycle, then eigenvalue corresponding to this cycle is $1$ and the eigenvector is $[0 \, \cdots 0 \, 1_{\gamma^{th}} \, 0 \cdots \, 0 ]$. Therefore the eigenvector matrix of $J$ can be written as $X = \frac{1}{\sqrt{2}}\begin{bmatrix}
	 	I_{m,m} & O_{m,1} & J_{m,m}\\
	 	O_{1,m} & \sqrt{2}I_{1,1} & O_{1,m}\\
	 	J_{m,m} & O_{m,1}  & -I_{m,m}
	 \end{bmatrix}$ corresponding to the eigenvalue matrix  $\Lambda_J = \begin{bmatrix}
	 	I_{m,m} & O_{m,1} & O_{m,m}\\
	 	O_{1,m} & I_{1,1} & O_{1,m}\\
	 	O_{m,m} & O_{m,1}  & -I_{m,m}
	 \end{bmatrix}$. And let $S = \begin{bmatrix}
	 A_{m,m} & O_{m,1} & O_{m,m}\\
	 B_{1,m} & cI_{1,1} & O_{1,m}\\
	 C_{m,m} & D_{m,1}  & E_{m,m}
 	\end{bmatrix}$. Therefore 
	\begin{align} \label{r2}
		&X^TSX \\ &= \frac{1}{2}\begin{bmatrix}
			I_{m,m} & O_{m,1} & J_{m,m}\\
			O_{1,m} & \sqrt{2}I_{1,1} & O_{1,m}\\
			J_{m,m} & O_{m,1}  & -I_{m,m}
		\end{bmatrix}  \begin{bmatrix}
		A_{m,m} & O_{m,1} & O_{m,m}\\
		B_{1,m} & cI_{1,1} & O_{1,m}\\
		C_{m,m} & D_{m,1}  & E_{m,m}
		\end{bmatrix}\begin{bmatrix}
		I_{m,m} & O_{m,1} & J_{m,m}\\
		O_{1,m} & \sqrt{2}I_{1,1} & O_{1,m}\\
		J_{m,m} & O_{m,1}  & -I_{m,m}
		\end{bmatrix} \nonumber\\ 
		&= \frac{1}{2}\begin{bmatrix}
			A_{m,m} + J_{m,m}C_{m,m} + J_{m,m}E_{m,m}J_{m,m} & \sqrt{2}J_{m,m}D_{m,1} & A_{m,m}J_{m,m}-J_{m,m}E_{m,m} + J_{m,m}C_{m,m}J_{m,m}\\
			\sqrt{2}B_{1,m} & 2cI_{1,1} & \sqrt{2}B_{1,m}J_{m,m}\\
			J_{m,m}A_{m,m} - E_{m,m}J_{m,m} - C_{m,m} & \sqrt{2}D_{m,1}  & J_{m,m}A_{m,m}J_{m,m} - C_{m,m}J_{m,m} + E_{m,m}
		\end{bmatrix}. \nonumber
	\end{align}
\hfill\end{proof}
\begin{corollary} \label{col_4}
	Let $S$ be the Caputo finite difference matrix as in \eqref{cap_mat_1} and $n$ is odd and $m:= \frac{n-1}{2}$. From relation \eqref{r2} we have:
	\begin{align*}
		&X^TSX \\ 
		&= \frac{1}{2}\begin{bmatrix}
			A_{m,m} + J_{m,m}C_{m,m} + J_{m,m}E_{m,m}J_{m,m} & \sqrt{2}J_{m,m}D_{m,1} & A_{m,m}J_{m,m}-J_{m,m}E_{m,m} + J_{m,m}C_{m,m}J_{m,m}\\
			\sqrt{2}B_{1,m} & 2cI_{1,1} & \sqrt{2}B_{1,m}J_{m,m}\\
			J_{m,m}A_{m,m} - E_{m,m}J_{m,m} - C_{m,m} & \sqrt{2}D_{m,1}  & J_{m,m}A_{m,m}J_{m,m} - C_{m,m}J_{m,m} + E_{m,m}
		\end{bmatrix}, \nonumber
	\end{align*}
	where
	\begin{align*}
		A_{m,m} &= \begin{bmatrix}
			d_{1,1} &  & & & \\
			d_{2,2} - d_{2,1} & d_{2,1} &  & & \\
			\vdots&  & \ddots & & \\
			d_{m-1,m-1} - d_{m-1,m-2} & \cdots &  & d_{m-1,1} & \\
			d_{m,m} - d_{m,m-1} & \cdots &  & d_{m,2} - d_{m,1} & d_{m,1}
		\end{bmatrix},\\
		E_{m,m} &= \begin{bmatrix}
			d_{m+2,1} &  & & & \\
			d_{m+3,2} - d_{m+3,1} & d_{m+3,1} &  & & \\
			\vdots&  & \ddots & & \\
			d_{n-1,m-1} - d_{n-1,m-2} & \cdots &  & d_{n-1,1} & \\
			d_{n,m} - d_{n,m-1} & \cdots &  & d_{n,2} - d_{n,1} & d_{n,1}
		\end{bmatrix},\\
		C_{m,m} &= \begin{bmatrix}
			d_{m+2,m+2}-d_{m+2,m+1} &  & \cdots & & d_{m+2,3}-d_{m+2,2}\\
			\vdots&  & \vdots & & \vdots \\
			d_{n,n} - d_{n,n-1} &  & \cdots &  & d_{n,m+2}-d_{n,m+1}
		\end{bmatrix}.
	\end{align*}
	and $B_{1,m} = [(d_{m+1,m+1} - d_{m+1,m}) \, \cdots \, (d_{m+1,2} - d_{m+1,1})]$, $c = d_{m+1,1}$, $D_{m,1} = [(d_{m+2,2}-d_{m+2,1}) \, \cdots \, (d_{n,m+1}-d_{n,m})]^T$.
   	Hence from the upper left block of $X^TSX$ i.e.
	\begin{align*}
		&A_{m,m} + J_{m,m}C_{m,m} + J_{m,m}E_{m,m}J_{m,m} \\
		&= \begin{bmatrix}
			d_{1,1} + d_{n,1} & d_{n,2} - d_{n,1} &  &\cdots & d_{n,m} - d_{n,m-1} \\
			d_{2,2} - d_{2,1} & d_{2,1} +d_{n-1,1} &  &\cdots & d_{n-1,m-1} - d_{n-1,m-2}\\
			\vdots&  & \ddots & & \vdots\\
			d_{m-1,m-1} - d_{m-1,m-2} &  & \cdots & d_{m-1,1} + d_{m+3,1} & d_{m+3,2} - d_{m+3,1} \\
			d_{m,m} - d_{m,m-1} &  & \cdots & d_{m,2} - d_{m,1} & d_{m,1} + d_{m+2,1}
		\end{bmatrix} \\
		&+ \begin{bmatrix}
			d_{n,n} - d_{n,n-1} &  & \cdots &  & d_{n,m+2}-d_{n,m+1}\\
			\vdots&  & \vdots & & \vdots \\
			d_{m+2,m+2}-d_{m+2,m+1} &  & \cdots & & d_{m+2,3}-d_{m+2,2}
		\end{bmatrix},
	\end{align*}
	the Gershgorin radius corresponding to the center $\bar{C}_1 = d_{1,1} + d_{n,1}, \bar{C}_2 = d_{2,1} +d_{n-1,1}, \cdots, \bar{C}_m = d_{m,1} + d_{m+2,1}$ are given by:
	\begin{align*}
		\bar{R}_1 &= (|d_{n,2} - d_{n,1}| + \cdots +|d_{n,m} - d_{n,m-1}|) + (|d_{n,m+2}-d_{n,m+1}| + \cdots + |d_{n,n} - d_{n,n-1}|)\\
		&= 	(d_{n,1}-d_{n,m})+(d_{n,m+1}-d_{n,n}),\\
		\vdots\\
		\bar{R}_m &= 	(|d_{m,2} - d_{m,1}| + \cdots + |d_{m,m} - d_{m,m-1}|) + (|d_{m+2,3}-d_{m+2,2}| + \cdots + |d_{m+2,m+2}-d_{m+2,m+1}|)\\
		&= (d_{m,1} - d_{m,m}) + (d_{m+2,2} -d_{m+2,m+2}).
	\end{align*} 
And from the lower-right block
	\begin{align*}
		&J_{m,m}A_{m,m}J_{m,m} - C_{m,m}J_{m,m} + E_{m,m}\\
		&= \begin{bmatrix}
			d_{m+2,1} +d_{m,1} & d_{m,2}-d_{m,1}  & &\cdots &d_{m,m} - d_{m,m-1} \\
			d_{m+3,2} - d_{m+3,1} & d_{m+3,1} + d_{m-1,1} &  & \cdots&d_{m-1,m-1}-d_{m-1,m-2} \\
			\vdots&  & \ddots & &\vdots \\
			d_{n-1,m-1} - d_{n-1,m-2} & \cdots &  & d_{n-1,1} + d_{2,1} & d_{2,2} - d_{2,1} \\
			d_{n,m} - d_{n,m-1} & \cdots &  & d_{n,2} - d_{n,1} & d_{n,1} + d_{1,1}
		\end{bmatrix}\\
		&- \begin{bmatrix}
			d_{m+2,3}-d_{m+2,2}	 &  & \cdots & & d_{m+2,m+2}-d_{m+2,m+1} \\
			\vdots&  & \vdots & & \vdots \\
			d_{n,m+2}-d_{n,m+1}	 &  & \cdots &  & d_{n,n} - d_{n,n-1},
		\end{bmatrix}
	\end{align*}
		we take the Gershgorin radius corresponding to the center $\underline{C}_{m+2} = d_{m+2,1} + d_{m,1}, \underline{C}_{m+3} = d_{m+3,1} + d_{m-1,1}, \cdots, \underline{C}_n = d_{n,1} +d_{1,1}$ as
	\begin{align*}
		\underline{R}_{m+2} &= 	(|d_{m,2} - d_{m,1}| + \cdots + |d_{m,m} - d_{m,m-1}|) + (|d_{m+2,3}-d_{m+2,2}| + \cdots + |d_{m+2,m+2}-d_{m+2,m+1}|)\\
		&= (d_{m,1} - d_{m,m}) + (d_{m+2,2} -d_{m+2,m+2}), \\ \vdots \\
		\underline{R}_{n} &= (|d_{n,2} - d_{n,1}| + \cdots +|d_{n,m} - d_{n,m-1}|) + (|d_{n,m+2}-d_{n,m+1}| + \cdots + |d_{n,n} - d_{n,n-1}|)\\
		&= 	(d_{n,1} - d_{n,m}) + (d_{n,m+1} - d_{n,n}).
	\end{align*} 
From the upper-right block
	\begin{align*}
		&A_{m,m}J_{m,m} - J_{m,m}E_{m,m} + J_{m,m}C_{m,m}J_{m,m}\\
		&= \begin{bmatrix}
			-(d_{n,m}-d_{n,m-1}) &\cdots  & & -(d_{n,2}-d_{n,1})& d_{1,1}- d_{n,1}\\
			-(d_{n-1,m-1}-d_{n-1,m-2})&\cdots  &  &d_{2,1} - d_{n-1,1} & d_{2,2} - d_{2,1}\\
			&  & \iddots & &\vdots \\
			-(d_{m+3,2}-d_{m+3,1}) & d_{m-1,1}-d_{m+3,1} &  & \cdots & d_{m-1,m-1} - d_{m-1,m-2}\\
			d_{m,1} - d_{m+2,1} & d_{m,2} - d_{m,1} &  & \cdots &  d_{m,m} - d_{m,m-1}
		\end{bmatrix}\\
		&+ \begin{bmatrix}
			d_{n,m+2}-d_{n,m+1} &  & \cdots &  &  d_{n,n} - d_{n,n-1}\\
			\vdots&  & \vdots & & \vdots \\
			d_{m+2,3}-d_{m+2,2} &  & \cdots & &  d_{m+2,m+2}-d_{m+2,m+1}
		\end{bmatrix},
	\end{align*}
		the contribution to the Gershgorin radius are:
	\begin{align*}
		\bar{R}_1 &= \left(|d_{n,2} - d_{n,1}| + \cdots +|d_{n,m} - d_{n,m-1}|\right) + \left(|d_{n,m+2}-d_{n,m+1}| + \cdots + |d_{n,n} - d_{n,n-1}|\right) + |d_{1,1} - d_{n,1}|\\
		&= (d_{n,1} - d_{n,m}) +(d_{n,m+1} - d_{n,n}) + |d_{1,1} - d_{n,1}|,\\
		\vdots \\
		\bar{R}_m &= |d_{m,1} - d_{m+2,1}| + \left(|d_{m,2} - d_{m,1}| + \cdots + |d_{m,m} - d_{m,m-1}| \right) + \left(|d_{m+2,3}-d_{m+2,2}| + \cdots + |d_{m+2,m+2}-d_{m+2,m+1}|\right) \\
		&= (d_{m+2,2} - d_{m+2,m+2}) + |d_{m,1} - d_{m+2,1}| + (d_{m,1} - d_{m,m}).
	\end{align*}
And from the lower-left block
	\begin{align*}
	&J_{m,m}A_{m,m} - E_{m,m}J_{m,m} - C_{m,m}\\
	&= \begin{bmatrix}
		d_{m,m} - d_{m,m-1} & \cdots &  &  d_{m,2} - d_{m,1} &   d_{m,1} - d_{m+2,1} \\
		d_{m-1,m-1} - d_{m-1,m-2} & \cdots &  &  d_{m-1,1}-d_{m+3,1} &  -(d_{m+3,2}-d_{m+3,1})\\
		&  & \iddots & &\vdots \\
		d_{2,2} - d_{2,1}	&d_{2,1} - d_{n-1,1}  &  & \cdots&  -(d_{n-1,m-1}-d_{n-1,m-2})\\
		d_{1,1}- d_{n,1} & -(d_{n,2}-d_{n,1}) & &\cdots &   -(d_{n,m}-d_{n,m-1})		
	\end{bmatrix}\\
	&-\begin{bmatrix}
		d_{m+2,m+2}-d_{m+2,m+1} &  & \cdots & & d_{m+2,3}-d_{m+2,2}\\
		\vdots&  & \vdots & & \vdots \\
		d_{n,n} - d_{n,n-1} &  & \cdots &  & d_{n,m+2}-d_{n,m+1}
	\end{bmatrix},
	\end{align*}
	the contribution to the Gershgorin radius are:
	\begin{align*}
		\bar{R}_{m+2} &= (d_{m,1} - d_{m,m}) + |d_{m,1} - d_{m+2,1}| + (d_{m+2,2} - d_{m+2,m+2}),\\
		\vdots \\
		\bar{R}_{n} &= |d_{1,1} - d_{n,1}| + (d_{n,1} - d_{n,m}) + (d_{n,m+1}- d_{n,n}).
	\end{align*}
	Also we have
	$J_{m,m}D_{m,1} = [(d_{n,m+1}-d_{n,m}) \, \cdots \, (d_{m+2,2}-d_{m+2,1})]^T$, $B_{1,m}J_{m,m} = [(d_{m+1,2} - d_{m+1,1}) \, \cdots \, (d_{m+1,m+1} - d_{m+1,m})]$. 
		Hence the Gershgorin radius of the matrix $X^TSX$ corresponding to the center $C_1 = \frac{d_{1,1} + d_{n,1}}{2}, \cdots, C_m = \frac{d_{m,1} + d_{m+2,1}}{2}, C_{m+1} = d_{m+1,1}, C_{m+2} = \frac{d_{m+2,1} + d_{m,1}}{2},  \cdots, C_n = \frac{d_{n,1} +d_{1,1}}{2}$ are given by:
	\begin{align*}
		R_1 &= (d_{n,1}-d_{n,m}) + (d_{n,m+1}- d_{n,n}) + \frac{1}{\sqrt{2}} (d_{n,m+1}-d_{n,m}) + \frac{1}{2}|d_{1,1} - d_{n,1}| \leq \frac{3}{2}(d_{1,1} + d_{n,1}),\\
		\vdots \\
		R_m &= (d_{m+2,2} - d_{m+2,m+2}) + \frac{1}{\sqrt{2}}(d_{m+2,2} - d_{m+2,1}) + \frac{1}{2}|d_{m,1} - d_{m+2,1}| + (d_{m,1} - d_{m,m})\leq \frac{3}{2}(d_{m,1} + d_{m+2,1}), \\
		R_{m+1} &= \frac{1}{\sqrt{2}}(|d_{m+1,m+1} - d_{m+1,m}| + \cdots + |d_{m+1,2} - d_{m+1,1}|) = \frac{1}{\sqrt{2}}(d_{m+1,1} - d_{m+1,m+1}) \leq \frac{1}{\sqrt{2}}d_{m+1,1},\\
		R_{m+2} &= (d_{m,1} - d_{m,m}) + \frac{1}{2}|d_{m,1} - d_{m+1,1}| + (d_{m+2,2} - d_{m+2,m+2}) + \frac{1}{\sqrt{2}} (d_{m+2,2} - d_{m+2,1}) \leq \frac{3}{2}(d_{m,1} + d_{m+2,1}), \\ \vdots \\
		R_{n} &= \frac{1}{2}|d_{1,1} - d_{n,1}| + (d_{n,1} - d_{n,m}) + (d_{n,m+1} - d_{n,n}) + \frac{1}{\sqrt{2}}(d_{n,m+1} - d_{n,m})\leq \frac{3}{2}(d_{1,1} + d_{n,1}).
	\end{align*}
\end{corollary}
\begin{lemma} \label{eigenvector}
	Let $S = \begin{bmatrix}
		S_1 & O\\
		S_2 & s
	\end{bmatrix}$ be a real lower triangular matrix of order $n+1$ and $S_1$ block is of order $n$. $B= \begin{bmatrix}
		B_1 & O\\
		O & 0
	\end{bmatrix}$ be a real block matrix of order $n+1$ with block $B_1$ of order $n$ and $O$ be the zero matrix. If $(S_1 -sI + i\beta B_1)$ is non-singular for positive constant $\beta$ then the matrix $S + i\beta B$ is similar to matrix $\bar{S} = \begin{bmatrix}
	S_1 + i\beta B_1 & O\\
	O & s
\end{bmatrix}$.
\end{lemma}
\begin{proof}
		Let $X = \begin{bmatrix}
		X_{11} & X_{12}\\
		X_{21} & X_{22}
	\end{bmatrix}$ be the block matrix of order $n+1$, where $X_{11}$ is of order $n$ block. Let
	\begin{align*}
		\begin{bmatrix}
			S_1 + i\beta B_1 & O\\
			S_2 & s
		\end{bmatrix} \begin{bmatrix}
		X_{11} & X_{12}\\
		X_{21} & X_{22}
	\end{bmatrix} &= \begin{bmatrix}
	X_{11} & X_{12}\\
	X_{21} & X_{22}
	\end{bmatrix} \begin{bmatrix}
	S_1 + i\beta B_1 & O\\
	O & s
	\end{bmatrix}\\
	\begin{bmatrix}
		(S_1 + i\beta B_1)X_{11} & (S_1 + i\beta B_1)X_{12}\\
		S_2X_{11} + sX_{21} & S_2X_{12} + sX_{22}
	\end{bmatrix} &= \begin{bmatrix}
	X_{11}(S_1 + i\beta B_1) & X_{12}l\\
	X_{21}(S_1 + i\beta B_1) &  X_{22}l
	\end{bmatrix}
	\end{align*}
	We choose $X_{11} = I, X_{12} = O, X_{22} = 1$ then we have $X_{21}(S_1 + i\beta B_1) = sX_{21} + S_2$. Therefore the system has solution if the matrix $(S_1 -sI + i\beta B_1)$ is invertiable.
	Hence the matrix $S + i\beta B$ is similar to $\bar{S}$ by the matrix $ \begin{bmatrix}
		I & O\\
		S_2(S_1 -sI + i\beta B_1)^{-1} & 1
	\end{bmatrix}$.
\hfill\end{proof}
\begin{remark} \label{eigenvector_re}
	If $S_1$ be toeplitz and $B_1$ be the exchange matrix $J$ then $|\det(S_1 -sI + i\beta J)| = |\det((S_1 -sI)J + i\beta I)|$. Now  $(S_1 -sI)J$ is symmetric so eigenvalues are all real hence $\det((S_1 -sI)J + i \beta I) \neq 0$ for any $\sigma$.
\end{remark}

%% file: Mesh.tex
Before going to main convergence results for DNWR and NNWR algorithm for the control problem we want to introduce a different kind of mesh grading technique, which is efficient and also helpful for the analysis.
\section{Mesh Grading}\label{Section_mesh}
	We are aware that utilizing a fractional diffusion graded mesh over time yields a higher order of accuracy, denoted as $\mathcal{O}(M^{-(2-\beta)})$, compared to a uniform mesh, which offers an accuracy of $\mathcal{O}(M^{-\beta})$, where $\beta$ represents fractional order and $M$ denotes number of nodes. Thus, it is customary to adopt a one-sided graded mesh, as suggested in \cite{stynes2017error}.
	
	In control problems, solving both the state and adjoint problems necessitates extending the one-sided graded mesh approach to a both-sided grading. Let's consider a total time domain $[0,T]$, divided into $n+1$ non-overlapping sub-intervals, with $n$ being odd, defined by the nodes $0=t_0 < t_1 <\cdots<t_{n+1}=T$. Initially, we employ the usual graded mesh for the time window $[0,T/2]$, utilizing nodes $t_0 < t_1 <\cdots<t_{(n+1)/2}$. Subsequently, for the remaining time window $[T/2,T]$, we apply a similar grading, starting from the right end, i.e., from $T$, using nodes $t_{n+1} > t_{n} >\cdots>t_{(n+1)/2}$. Thus, for $j = 0, 1,\ldots,(n+1)/2$, the time nodes are defined as $t_j = (j/(n+1)/2)^{(2- \alpha)/\alpha} T/2$ and $t_{n+1 - j} = T - t_j$.

%% file: Proof_DNWR.tex
\section{Convergence of Dirichlet-Neumann Waveform Relaxation}\label{sectionDNWR}
Presently, we analyze the convergence assessments of the DNWR algorithm \eqref{DNWR_state}-\eqref{DNWR_adjoint}, selecting the 1D model of the sub-diffusion and diffusion equation for enhanced theoretical simplicity and algebraic straightforwardness in our evaluations. Our analysis is directed towards a heterogeneous scenario where $\kappa(\boldsymbol{x},t)=\kappa_1$ over $\Omega_{1}=(-h_1,0)$ and $\kappa(\boldsymbol{x},t)=\kappa_2$ over $\Omega_{2}=(0,h_2)$.
We introduce the interface error $\omega^{(k)}(t)$ on state equation and $\upsilon^{(k)}(t)$ on adjoint equation and using the lemma \ref{rl_eqi_c} for transforming the Riemann-Liouville derivative to Caputo derivative in the error equation to obtain:
\begin{equation}\label{DNWR_error_state}
	\begin{array}{rcll}
		\begin{cases}
			_{0}D_t^{\alpha} y_{1}^{(k)} = \kappa_1\partial_{xx} y_{1}^{(k)} - \frac{p_1^{(k)}}{\sigma}, & \textrm{in}\; (-h_1, 0),\\
			y_{1}^{(k)}(\boldsymbol{x},0) = 0, & \textrm{in}\; (-h_1, 0),\\
			y_{1}^{(k)}(0,t) = \omega^{(k-1)}(t), \\
			y_{1}^{(k)}(-h_1,t)=0, 
		\end{cases}
	\end{array}\
	\begin{array}{rcll}
		\begin{cases}
			_{0}D_t^{\alpha} y_{2}^{(k)} = \kappa_1\partial_{xx} y_{2}^{(k)} - \frac{p_2^{(k)}}{\sigma}, & \textrm{in}\; (0, h_2),\\
			y_{2}^{(k)}(\boldsymbol{x},0) = 0, & \textrm{in}\; (0, h_2),\\
			\kappa_2\partial_x y_{2}^{(k)}(0,t)  =  \kappa_1\partial_x y_{1}^{(k)}(0,t),\\
			y_{2}^{(k)}(h_2,t) = 0,
		\end{cases}
	\end{array}
\end{equation}

\begin{equation}\label{DNWR_error_adjoint}
	\begin{array}{rcll}
		\begin{cases}
			_{T}D_t^{\alpha}p_{1}^{(k)} = \kappa_1 \partial_{xx} p_{1}^{(k)} + y_1^{(k)}, & \textrm{in}\; (-h_1,0),\\
			p_{1}^{(k)}(\boldsymbol{x},T) = 0, & \textrm{in}\; (-h_1,0),\\
			p_{1}^{(k)}(0,t) = \upsilon^{(k-1)}(t), \\
			p_{1}^{(k)}(-h_1,t) = 0, 
		\end{cases}
	\end{array}\
	\begin{array}{rcll}
		\begin{cases}
			_{T}D_t^{\alpha}p_{2}^{(k)} = \kappa_2 \partial_{xx} p_{2}^{(k)} + y_2^{(k)}, & \textrm{in}\; (0, h_2),\\
			p_{2}^{(k)}(\boldsymbol{x},T) = 0, & \textrm{in}\; (0, h_2),\\
			\kappa_2\partial_x p_{2}^{(k)}(0,t)  =  \kappa_1\partial_x p_{1}^{(k)}(0,t),\\
			p_{2}^{(k)}(h_2, 0) = 0, 
		\end{cases}
	\end{array}
\end{equation}
with the update step:
\begin{align}\label{DNWR_error_update}
	\omega^{(k)}(t) &= \theta y_{2}^{(k)}(0,t) + (1-\theta) \omega^{(k-1)}(t),\\
	\upsilon^{(k)}(t) &= \phi p_{2}^{(k)}(0,t) + (1-\phi) \upsilon^{(k-1)}(t).
\end{align}
Now, we will conduct an analysis of the system using semi-discrete methods. We define time nodes as $0 = t_0 < t_1 < \cdots < t_{n+1} = T$, where $n$ belongs to the set of natural numbers. Here, we have the option to employ either a graded or uniform time mesh.
When employing the L1 scheme for both left and right-sided Caputo fractional derivatives, the discretized matrices for these derivatives will be the same in uniform mesh for any $n \in \mathbb{N}$ and in both sides of a graded mesh when $n$ is odd. However, for the usual one-sided graded mesh, they will differ. This is because both the uniform and two-sided graded meshes adhere to a reflective symmetry, i.e., $t_j = T - t_{n+1-j}$ for $j = 0, 1, \cdots, n+1$. This symmetry condition is not satisfied for the one-sided graded mesh. As the graded mesh inherently accounts for the fractional order, as the fractional order $\alpha$ tends towards $1$, the graded mesh converges towards the uniform mesh.
With these discretization mechanisms in place, we can proceed with the analysis. Let, the discretized time fractional derivative defined as
\begin{equation*} \label{matrix_1}
	L := \begin{bmatrix}
		d_{1,1} &  & & & \\
		d_{2,2} - d_{2,1} & d_{2,1} &  & & \\
		\vdots& \ddots & \ddots & & \\
		d_{n,n} - d_{n,n-1} & \cdots &  & d_{n,1} & \\
		d_{n+1,n+1} - d_{n+1,n} & \cdots &  & d_{n+1,2} - d_{n+1,1} & d_{n+1,1}
	\end{bmatrix}. 
	\textit{ And }
	\bar{J} := \begin{bmatrix}
		& &  & 1 &  \\
		& & 1 &  &  \\
		& \iddots & & \\
		1 &  & & &  \\
		& & & & 
	\end{bmatrix}
\end{equation*} 
be the order $n+1$ square matrix coming from the source term in the right hand side. Here $d_{m,j}$ for $1 \leq j \leq m \leq n+1$ are given in \eqref{diagonal_element}. For the sake of the analysis we choose only the uniform and both sided graded mesh and have the following state and adjoint equation:
\begin{equation}\label{DNWR_discrete_state}
	\begin{array}{rcll}
		\begin{cases}
			L Y_{1}^{(k)}(\boldsymbol{x}) = \kappa_1\partial_{xx} Y_{1}^{(k)}(\boldsymbol{x}) - \frac{1}{\sigma}\bar{J} P_1^{(k)}(\boldsymbol{x}), & \textrm{in}\; (-h_1, 0),\\
			Y_{1}^{(k)}(0) = \Omega^{(k-1)}, \\
			Y_{1}^{(k)}(-h_1)=0, 
		\end{cases}
	\end{array}\
	\begin{array}{rcll}
		\begin{cases}
			L Y_{2}^{(k)}(\boldsymbol{x}) = \kappa_2\partial_{xx} Y_{2}^{(k)}(\boldsymbol{x}) - \frac{1}{\sigma}\bar{J} P_2^{(k)}(\boldsymbol{x}), & \textrm{in}\; (0, h_2),\\
			\kappa_2\partial_x Y_{2}^{(k)}(0)  =  \kappa_1\partial_x Y_{1}^{(k)}(0),\\
			Y_{2}^{(k)}(h_2) = 0,
		\end{cases}
	\end{array}
\end{equation}
and
\begin{equation}\label{DNWR_discrete_adjoint}
	\begin{array}{rcll}
		\begin{cases}
			L P_{1}^{(k)}(\boldsymbol{x}) = \kappa_1 \partial_{xx} P_{1}^{(k)}(\boldsymbol{x}) + \bar{J} Y_1^{(k)} (\boldsymbol{x}), & \textrm{in}\; (-h_1,0),\\
			P_{1}^{(k)}(0) = \Upsilon^{(k-1)}, \\
			P_{1}^{(k)}(-h_1) = 0, 
		\end{cases}
	\end{array}\
	\begin{array}{rcll}
		\begin{cases}
			L P_{2}^{(k)}(\boldsymbol{x}) = \kappa_2 \partial_{xx} P_{2}^{(k)}(\boldsymbol{x}) + \bar{J} Y_2^{(k)} (\boldsymbol{x}), & \textrm{in}\; (0, h_2),\\
			\kappa_2\partial_x P_{2}^{(k)}(0)  =  \kappa_1\partial_x P_{1}^{(k)}(0),\\
			P_{2}^{(k)}(h_2) = 0, 
		\end{cases}
	\end{array}
\end{equation}
with the update in discretize form:
\begin{align}\label{DNWR_discrete_update}
	\Omega^{(k)} &= \theta Y_{2}^{(k)} + (1-\theta) \Omega^{(k-1)},\\
	\Upsilon^{(k)} &= \phi P_{2}^{(k)} + (1-\phi) \Upsilon^{(k-1)}.
\end{align}
Here we have taken the discretization of $Y_i := [y_i^1, y_i^2, ..., y_i^{n+1}]^T$ in forward direction and  $P_i$ into reverse direction i.e. $P_i := [p_i^{n}, p_i^{n-1}, ..., p_i^{0}]^T$. For sake of analysis we combine both state and adjoint system into a single system and choosing the dependent variable $X_i := [Y_i; P_i]$ Using the above two equations we form the system:
\begin{equation}\label{DNWR_discrete_system}
	\begin{array}{rcll}
		\begin{cases}
			\mathbb{L} X_{1}^{(k)}(\boldsymbol{x}) = \kappa_1\partial_{xx} X_{1}^{(k)}(\boldsymbol{x}), & \textrm{in}\; (-h_1, 0),\\
			X_{1}^{(k)}(0) = \Pi^{(k-1)}, \\
			X_{1}^{(k)}(-h_1)=0, 
		\end{cases}
	\end{array}\
	\begin{array}{rcll}
		\begin{cases}
			\mathbb{L} X_{2}^{(k)}(\boldsymbol{x}) = \kappa_2\partial_{xx} X_{2}^{(k)}(\boldsymbol{x}), & \textrm{in}\; (0, h_2),\\
			\kappa_2\partial_x X_{2}^{(k)}(0)  =  \kappa_1\partial_x X_{1}^{(k)}(0),\\
			X_{2}^{(k)}(h_2) = 0,
		\end{cases}
	\end{array}
\end{equation}
where, $\mathbb{L} := \begin{bmatrix}
	L & \frac{1}{\sigma}\bar{J}\\
	-\bar{J} & L
	\end{bmatrix}$. We choose $\Pi^{(k)} := [\Omega^{(k)}; \Upsilon^{(k)}]$ for the combine form of the update condition for the new system and get:
\begin{equation}\label{DNWR_discrete_system_update}
	\Pi^{(k)} = \begin{bmatrix}
		\theta I & 0\\
		0 & \phi I
	\end{bmatrix} X_{2}^{(k)} + \begin{bmatrix}
	(1-\theta) I & 0\\
	0 & (1-\phi) I
\end{bmatrix} \Pi^{(k-1)}.
\end{equation}
Before going to the main convergent result of DNWR algorithm we first need to check that the time discretization is a convergent scheme of the continuous version, i.e. we need to check the consistency and stability of the numerical scheme. Which is achieved by the following theorem.
\begin{theorem} \label{Th_1}
	For the L1 numerical scheme in time fractional derivative the semi-discrete ode system in \eqref{DNWR_discrete_system} is convergent for the pde system \eqref{DNWR_error_state} and \eqref{DNWR_error_adjoint}. 
\end{theorem}
\begin{proof}
	We need to show that the numerical scheme is consistent and stable for the convergence of the solution. In discretization we have used L1-scheme on graded mesh or on uniform mesh, which is consistent. Now to show this scheme is stable we use the similar technique as in lemma \ref{stable}.
	For that we first need the error equation of the semi discrete pde, which can be obtained by the substruction of the semi discrete pde with the continuous pde with the Taylor approximation of the time derivative. which is given as follows:
\begin{equation}\label{DNWR_discrete_convergence}
		\begin{cases}
			\mathbb{L} E_{1}(x) = \partial_{xx} E_{1}(x) + F(x), & \textrm{in}\; (-h_1, 0),\\
			E_{1}(0) = 0, \\
			E_{1}(-h_1)=0, 
		\end{cases}
\end{equation}
Where $F(x)$ is the remaining term of the Taylor series expansion in time. To show the stability we use the Fourier frequency in space. We see that Fourier sine expansion is enough for that purpose. We consider each frequency separately and we have:
\begin{align}\label{DNWR_discrete_convergence1}
		\mathbb{L} \hat{E}_{1} = -\xi^2 \hat{E}_{1} + \hat{F}, \\
		\hat{E}_{1} = (\mathbb{L} + \xi^2 \mathbb{I})^{-1}\hat{F}. \nonumber
\end{align}
Therefor the error is contingent upon the smallest absolute eigenvalue of the matrix $(\mathbb{L} + \xi^2 \mathbb{I})$. Now $\mathbb{L} = \begin{bmatrix}
	L & \frac{1}{\sigma}\bar{J}\\
	-\bar{J} & L
\end{bmatrix}$. From lemma \ref{matrix} we know that eigenvalues of $\mathbb{L}$ are the collection of eigenvalues of $L + \frac{i}{\sqrt{\sigma}}\bar{J}$ and $L - \frac{i}{\sqrt{\sigma}}\bar{J}$. Let, $L = \begin{bmatrix}
L_1 & 0\\
L_2 & d_{n+1,1}
\end{bmatrix}$ be the order $(n+1)$ block matrix with $L_1$ be the lower triangular block of order $n$ and  $\bar{J} = \begin{bmatrix}
J & O\\
O & 0
\end{bmatrix}$ be the order $(n+1)$ block matrix with $J$ be the exchange matrix of order $n$. From corollary \ref{frac_derivative} we know that each diagonal element from L1 scheme is positive, so $d_{n+1,1} > 0$. Again from the same corollary we have the matrix $L_1$ is positive definite, hence from lemma \ref{positive_eigenvalues} we conclude that the eigenvalues of $L + \frac{i}{\sqrt{\sigma}}\bar{J}$ consist positive real part. In a similar manner we can also show the same for $L - \frac{i}{\sqrt{\sigma}}\bar{J}$. 
Hence $(\mathbb{L} + \xi^2 \mathbb{I})$ be a nonsingular matrix and the error is bounded. Similar procedure can be used for the other part of \eqref{DNWR_discrete_system} with Neumann boundary condition, where we have to use cosine function with certain phase difference.
\hfill\end{proof}\\

Let the matrix $\mathbb{L}$ be diagonalizable by the eigenvector matrix $P$ corresponding to the eigenvalue matrix $\Lambda$, which means $\mathbb{L} = P\Lambda P^{-1}$. For the sake of notation simplicity, $f(\Lambda)$, where $f$ is any scalar-valued function, implies that the function $f$ operates solely on the diagonal elements of $\Lambda$ in a pointwise manner, while all other elements are zero. Solutions of \eqref{DNWR_discrete_system} take the form:
\begin{align}
	X_{1}^{(k)}(x) &= P\sinh\left(
		(x+h_1)\sqrt{\Lambda/\kappa_1}\right)\cosech(h_1\sqrt{\Lambda/\kappa_1})P^{-1} \Pi^{(k-1)}, \nonumber\\
	X_{2}^{(k)}(x) &= \sqrt{\frac{\kappa_{1}}{\kappa_{2}}} P
	   \sinh((x-h_2)\sqrt{\Lambda/\kappa_2})\sech(h_2\sqrt{\Lambda/\kappa_2}) \coth(h_1\sqrt{\Lambda/\kappa_1})P^{-1} \Pi^{(k-1)}. \label{DNWRsol1}
\end{align}
Substituting \eqref{DNWRsol1} into \eqref{DNWR_discrete_system_update}, the recurrence relation for $k\in \mathbb{N}$ becomes: 
\begin{align}\label{DNWRsol2}
	\Pi^{(k)} &= -\begin{bmatrix}
		\theta I & 0\\
		0 & \phi I
	\end{bmatrix} \sqrt{\frac{\kappa_{1}}{\kappa_{2}}} P\tanh(h_2\sqrt{\Lambda/\kappa_2})\coth(h_1\sqrt{\Lambda/\kappa_1})
	P^{-1}\Pi^{(k-1)} + \begin{bmatrix}
		(1-\theta) I & 0\\
		0 & (1-\phi) I
	\end{bmatrix}\Pi^{(k-1)}.
\end{align}
For notational simplicity, consider $a := h_1/\sqrt{\kappa_1}, b := h_2/\sqrt{\kappa_2}$. Then for $a \neq b$ \eqref{DNWRsol2} can be written as
\begin{align}\label{DNWRsol3}
	\Pi^{(k)} &= -\begin{bmatrix}
		\theta I & 0\\
		0 & \phi I
	\end{bmatrix} \sqrt{\frac{\kappa_{1}}{\kappa_{2}}} P\left(\tanh(b\sqrt{\Lambda})\coth(a\sqrt{\Lambda})-I\right)
	P^{-1}\Pi^{(k-1)}\\ 
	&\hspace{20mm} + \begin{bmatrix}
		(1-(1+\sqrt{\kappa_1/\kappa_2})\theta) I & 0\\
		0 & (1-(1+\sqrt{\kappa_1/\kappa_2})\phi) I
	\end{bmatrix}\Pi^{(k-1)}. \nonumber
\end{align}
Choose $\theta = \phi = 1/(1+\sqrt{\kappa_1/\kappa_2})$ we have:
\begin{align}\label{DNWRsol4}
	\Pi^{(k)} &= -\frac{\sqrt{\kappa_{1}}}{\sqrt{\kappa_1} + \sqrt{\kappa_2}} P\left(\tanh(b\sqrt{\Lambda})\coth(a\sqrt{\Lambda})-I\right)
	P^{-1}\Pi^{(k-1)}.
\end{align}
\begin{lemma}\label{estimate}
	For any $\alpha, \beta$ nonzero positive real number and $z \in \mathbb{C}$ with $\Re(z))>0$ 
	\begin{equation*}
		\left|\frac{\exp(-\alpha z)}{1 \pm \exp(-\beta z)}\right| \leq \frac{\exp(-\alpha \Re(z))}{1-\exp(-\beta \Re(z))}.
	\end{equation*}	
\hfill\end{lemma}
\begin{proof}
	Let $z = x + iy, x>0$. Then
	\begin{align*}
		\left|\frac{\exp(-\alpha z)}{1 \pm \exp(-\beta z)}\right| &= \frac{\exp(-\alpha x)}{\sqrt{(1 \pm \exp(-\beta x)\cos(\beta y))^2 + (\exp(-\beta x) \sin(\beta y))^2}}\\
		& = \frac{\exp(-\alpha x)}{\sqrt{1\pm 2\exp(-\beta x)\cos(\beta y) + \exp(-2\beta x)}}\\
		& \leq \frac{\exp(-\alpha x)}{1 - \exp(-\beta x)}.
	\end{align*} 
\hfill\end{proof}
\begin{theorem}\label{Th_2}
	Let time discretize combined state and adjoint  matrix $\mathbb{L}$ in \eqref{DNWR_discrete_system} be diagonalizable by the eigenvector matrix $P$ corresponding to the eigenvalues $\{\lambda_i: i = 1,\ldots,2n+2\}$.
	If $\lambda = \min(\Re(\sqrt{\lambda_i}))$, then the interface error term of the DNWR algorithm in \eqref{DNWR_state}-\eqref{DNWR_adjoint} in 1D for relaxation parameters $\theta = \phi = 1/(1+\sqrt{\kappa_1/\kappa_2})$ satisfies the following bound
	\begin{equation*}
		\|\Pi^{(k)}\|_{\infty} < \left(\frac{\sqrt{\kappa_{1}}}{\sqrt{\kappa_1} + \sqrt{\kappa_2}}\right)^k \left(\frac{\cosh((b-a)\lambda)}{\sinh(a\lambda) \sinh(b\lambda)}\right)^k  \|P\|_{\infty} \|P^{-1}\|_{\infty} \|\Pi^{(0)}\|_{\infty}.
	\end{equation*}
\end{theorem}
\begin{proof}
	Let $\mathbb{L}$ be diagonalizable by eigenvector matrix $P$ corresponding to the eigenvalue matrix $\Lambda$, i.e. $\mathbb{L} = P\Lambda P^{-1}$. We know from the proof of Theorem \ref{Th_1} that all the eigenvalues of $\mathbb{L}$ lies in first or fourth quadrant in complex plane therefore the principal square roots of eigenvalues also lie in first or fourth quadrant, hence has positive real part i.e. $\lambda = \min(\Re(\sqrt{\lambda_i}))$ is well defined.
	
	 Now for the proof we take the update part \eqref{DNWRsol4} then left multiplying on both side with $P^{-1}$ to have 
	\begin{align}\label{DNWRsol5}
		P^{-1}\Pi^{(k)} &= -\frac{\sqrt{\kappa_{1}}}{\sqrt{\kappa_1} + \sqrt{\kappa_2}} \left(\tanh(b\sqrt{\Lambda})\coth(a\sqrt{\Lambda})-I\right)P^{-1}\Pi^{(k-1)} \\
		&= -\frac{\sqrt{\kappa_{1}}}{\sqrt{\kappa_1} + \sqrt{\kappa_2}} \frac{\sinh((b-a)\sqrt{\Lambda})}{\sinh(b\sqrt{\Lambda})\cosh(a\sqrt{\Lambda})} P^{-1}\Pi^{(k-1)} =: -\frac{\sqrt{\kappa_{1}}}{\sqrt{\kappa_1} + \sqrt{\kappa_2}} W P^{-1}\Pi^{(k-1)}.	\nonumber
	\end{align}
	Let $(W)_{i i} := \frac{\sinh((b-a)\sqrt{\lambda_i})}{\sinh(b\sqrt{\lambda_i})\cosh(a\sqrt{\lambda_i})}$ be the weight. And we will estimate the weight by simplifying the hyperbolic terms in terms of exponential, which is as follows
	\begin{align}
		 (W)_{ii}
		 &= 2\frac{\left(\exp((b-a)\sqrt{\lambda_i}) -
		 \exp(-(b-a)\sqrt{\lambda_i})\right)}{\left(\exp(a\sqrt{\lambda_i}) - \exp(-a\sqrt{\lambda_i})\right) \left(\exp(b\sqrt{\lambda_i}) + \exp(-b\sqrt{\lambda_i})\right)} \\
 	     &= 2\frac{\left(\exp(-2a\sqrt{\lambda_i}) -
 		\exp(-2b\sqrt{\lambda_i})\right)}{\left(1 - \exp(-2a\sqrt{\lambda_i})\right) \left(1 + \exp(-2b\sqrt{\lambda_i})\right)} \nonumber\\
 		 &= 2\frac{\exp(-2a\sqrt{\lambda_i}) -
 		\exp(-2b\sqrt{\lambda_i})}{\left(1 - \exp(-2a\sqrt{\lambda_i})\right) \left(1 + \exp(-2b\sqrt{\lambda_i})\right)} \nonumber\\
 		&= \frac{2\exp(-2a\sqrt{\lambda_i})}{\left(1 - \exp(-2a\sqrt{\lambda_i})\right) \left(1 + \exp(-2b\sqrt{\lambda_i})\right)} - \frac{2
 		\exp(-2b\sqrt{\lambda_i})}{\left(1 - \exp(-2a\sqrt{\lambda_i})\right) \left(1 + \exp(-2b\sqrt{\lambda_i})\right)} \nonumber\\
 	    & := I_1 - I_2. \nonumber
	\end{align}
	 Now
	 \begin{align}
	 	|I_1| &= 2\left|\frac{\exp(-a\sqrt{\lambda_i})}{\left(1 - \exp(-2a\sqrt{\lambda_i})\right) }\right| \left|\frac{\exp(-a\sqrt{\lambda_i})}{\left(1 + \exp(-2b\sqrt{\lambda_i})\right) }\right|.
	 \end{align}
 	Using the inequality from lemma \ref{estimate}, we have 
 	 \begin{align}
 		|I_1| 
 		& \leq 2 \frac{\exp(-2a \Re(\sqrt{\lambda_i}))}{\left(1 - \exp(-2a\Re(\sqrt{\lambda_i}))\right)\left(1 - \exp(-2b\Re(\sqrt{\lambda_i}))\right)}.  \nonumber
 	\end{align}
 	Further using $\lambda =  \min(\Re(\sqrt{\lambda_i}))$ and the property of monotonic decreasing function, we obtain
 	\begin{align}
 		|I_1| 
 		& \leq 2 \frac{\exp(-2a \lambda))}{\left(1 - \exp(-2a\lambda))\right)\left(1 - \exp(-2b\lambda))\right)}.  \nonumber
 	\end{align}
     By similar argument we have the bound
     \begin{align}
     	|I_2| &= 2\left|\frac{\exp(-b\sqrt{\lambda_i})}{\left(1 - \exp(-2a\sqrt{\lambda_i})\right) }\right| \left|\frac{\exp(-b\sqrt{\lambda_i})}{\left(1 + \exp(-2b\sqrt{\lambda_i})\right) }\right| \\
     	& \leq 2 \frac{\exp(-2b\lambda)}{\left(1 - \exp(-2a\lambda)\right)\left(1 - \exp(-2b\lambda)\right)}.\nonumber
     \end{align}
     Therefore
     \begin{align} \label{DNWRsol6}
		|(W)_{ii}| & \leq 2 \frac{\exp(-2a\lambda) + \exp(-2b\lambda)}{\left(1 - \exp(-2a\lambda)\right)\left(1 - \exp(-2b\lambda)\right)} \\
		& \leq \frac{\cosh((b-a)\lambda)}{\sinh(a\lambda) \sinh(b\lambda)}. \nonumber
     \end{align}
     Hence $\|W\|_{\infty} \leq \frac{\cosh((b-a)\lambda)}{\sinh(a\lambda) \sinh(b\lambda)} $.
	 Now taking the supremum norm on both side in \eqref{DNWRsol5} we have
	 \begin{align}\label{DNWRsol8}
	 	\left\|P^{-1}\Pi^{(k)}\right\|_{\infty} 
	 	&\leq \frac{\sqrt{\kappa_{1}}}{\sqrt{\kappa_1} + \sqrt{\kappa_2}} \frac{\cosh((b-a)\lambda)}{\sinh(a\lambda) \sinh(b\lambda)}\left\|P^{-1}\Pi^{(k-1)}\right\|_{\infty}  \\
	 	&\leq \left(\frac{\sqrt{\kappa_{1}}}{\sqrt{\kappa_1} + \sqrt{\kappa_2}}\right)^k \left(\frac{\cosh((b-a)\lambda)}{\sinh(a\lambda) \sinh(b\lambda)}\right)^k \left\|P^{-1}\Pi^{(0)}\right\|_{\infty}. \nonumber
	 \end{align}
   Again using the norm inequality we have the required result
   \begin{equation*}
   	\|\Pi^{(k)}\|_{\infty} \leq \left(\frac{\sqrt{\kappa_{1}}}{\sqrt{\kappa_1} + \sqrt{\kappa_2}}\right)^k \left(\frac{\cosh((b-a)\lambda)}{\sinh(a\lambda) \sinh(b\lambda)}\right)^k  \|P\|_{\infty} \|P^{-1}\|_{\infty} \|\Pi^{(0)}\|_{\infty}.
   \end{equation*}
\hfill\end{proof}

%% file: Proof_NNWR.tex
\section{Convergence of Neumann-Neumann Waveform Relaxation} \label{NNWR_con}
To analyze the convergence for the system described by equations \eqref{NNWR_state}-\eqref{NNWR_adjoint} in 1D setting, we initially partition the domain $(L_1, L_2)$ into $N$ non-overlapping subdomains $(x_{i-1}, x_i)$, where $i = 1, 2, ..., N$. Each subdomain is characterized by a length $a_i$, and we consider distinct constant generalized diffusion coefficients $\kappa_i$ within each subdomain $\Omega_i$. Additionally, we introduce $\omega_i^{(k)}(t)$ and $\upsilon_i^{(k)}(t)$ to represents the error terms at ${x_i}_{i=1}^{N-1}$ for the state and adjoint equations, respectively. Than set boundary conditions as $\omega_0^{(k)}(t)=0, \upsilon_0^{(k)}(t)=0$ at $x=x_0$, and $\omega_N^{(k)}(t)=0, \upsilon_N^{(k)}(t)=0$ at $x=x_N$. Now using the Lemma \ref{rl_eqi_c} to convert the Riemann-Liouville derivative into the Caputo derivative within the error equation yields:
\begin{equation}\label{NNWR_error_state}
	\begin{array}{rcll}
		\begin{cases}
			_{0}D_t^{\alpha}y_{i}^{(k)} = \kappa_i \frac{d^2}{dx^2} y_{i}^{(k)} - \frac{p_i^{(k)}}{\sigma}, & \textrm{in}\; (x_{i-1},x_i),\\
			y_{i}^{(k)}(\boldsymbol{x},0) = 0, & \textrm{in}\; (x_{i-1},x_i),\\
			y_{i}^{(k)} = \omega_{i-1}^{(k-1)}, & \textrm{at }x_{i-1}, i \neq 1,\\
			y_{i}^{(k)} = \omega_{i}^{(k-1)}, & \textrm{at }x_{i}, i \neq N,
		\end{cases}
	\end{array}
	\begin{array}{rcll}
		\begin{cases}
			_{0}D_t^{\alpha}\psi_{i}^{(k)} = \kappa_i \frac{d^2}{dx^2} \psi_{i}^{(k)} - \frac{\chi_i^{(k)}}{\sigma}, & \textrm{in}\; (x_{i-1},x_i),\\
		    \psi_{i}^{(k)}(\boldsymbol{x},0) = 0, & \textrm{in}\; (x_{i-1},x_i),\\
			-\kappa_i\frac{d}{dx}\psi_{i}^{(k)} = \kappa_{i-1}\frac{d}{dx}y_{i-1}^{(k)} - \kappa_{i}\frac{d}{dx}y_{i}^{(k)} & \textrm{at }x_{i-1}, i \neq 1,\\
			\kappa_i\frac{d}{dx}\psi_{i}^{(k)} = \kappa_i\frac{d}{dx}y_{i}^{(k)} - \kappa_{i+1}\frac{d}{dx}y_{i+1}^{(k)}, & \textrm{at }x_{i}, i \neq N-1.
		\end{cases}
	\end{array}
\end{equation}
\begin{equation}\label{NNWR_error_adjoint}
	\begin{array}{rcll}
		\begin{cases}
			_{T}D_t^{\alpha}p_{i}^{(k)} = \kappa_i \frac{d^2}{dx^2} p_{i}^{(k)} + y_i^{(k)}, & \textrm{in}\; (x_{i-1},x_i),\\
			p_{i}^{(k)}(\boldsymbol{x},T) = 0, & \textrm{in}\; (x_{i-1},x_i),\\
			p_{i}^{(k)} = \upsilon_{i-1}^{(k-1)}, & \textrm{at }x_{i-1}, i \neq 1,\\
			p_{i}^{(k)} = \upsilon_{i}^{(k-1)}, & \textrm{at }x_{i}, i \neq N,
		\end{cases}
	\end{array}
	\begin{array}{rcll}
		\begin{cases}
			_{T}D_t^{\alpha}\chi_{i}^{(k)} = \kappa_i \frac{d^2}{dx^2} \chi_{i}^{(k)} + \psi_{i}^{(k)}, & \textrm{in}\; (x_{i-1},x_i),\\
			\chi_{i}^{(k)}(\boldsymbol{x},T) = 0, & \textrm{in}\; (x_{i-1},x_i),\\
			-\kappa_i\frac{d}{dx}\chi_{i}^{(k)} = \kappa_{i-1}\frac{d}{dx}p_{i-1}^{(k)}+\kappa_i\frac{d}{dx}p_{i}^{(k)}, & \textrm{at }x_{i-1}, i \neq 1,\\
			\kappa_i\frac{d}{dx}\chi_{i}^{(k)} = \kappa_i\frac{d}{dx}p_{i}^{(k)}+\kappa_{i+1}\frac{d}{dx}p_{i+1}^{(k)}, & \textrm{at }x_{i}, i \neq N,
		\end{cases}
	\end{array}
\end{equation}
with the update conditions at interfaces for $i=1,\ldots,N-1$:
\begin{align}\label{NNWR_error_update}
	\omega_{i}^{(k)}(x_i,t)&=\omega_{i}^{(k-1)}(x_i,t)-\theta_{i}
	\left( \psi_{i}^{(k)}(x_i,t)+\psi_{i+1}^{(k)}(x_i,t)\right),\\
	\upsilon_{i}^{(k)}(x_i,t)&=\upsilon_{i}^{(k-1)}(x_i,t)-\phi_{i}
	\left( \chi_{i}^{(k)}(x_i,t)+\chi_{i+1}^{(k)}(x_i,t)\right).
\end{align}
After employing the same discretization technique for the time fractional derivative, which encompasses both uniform and two-sided graded meshes as outlined in section \ref{sectionDNWR}, the NNWR system can be readily expressed as:
\begin{equation}\label{NNWR_discrete_system}
	\begin{array}{rcll}
		\begin{cases}
			\mathbb{L} X_{i}^{(k)} = \kappa_i \frac{d^2}{dx^2} X_{i}^{(k)}, & \textrm{in}\; (x_{i-1},x_i),\\
			X_{i}^{(k)} = \Pi_{i-1}^{(k-1)}, & \textrm{at }x_{i-1}, i \neq 1,\\
			X_{i}^{(k)} = \Pi_{i}^{(k-1)}, & \textrm{at }x_{i}, i \neq N,
		\end{cases}
	\end{array}
	\begin{array}{rcll}
		\begin{cases}
			\mathbb{L}\Xi_{i}^{(k)} = \kappa_i \frac{d^2}{dx^2} \Xi_{i}^{(k)}, & \textrm{in}\; (x_{i-1},x_i),\\
			-\kappa_i\frac{d}{dx}\Xi_{i}^{(k)} = \kappa_{i-1}\frac{d}{dx}X_{i-1}^{(k)}+\kappa_i\frac{d}{dx}X_{i}^{(k)}, & \textrm{at }x_{i-1}, i \neq 1,\\
			\kappa_i\frac{d}{dx}\Xi_{i}^{(k)} = \kappa_i\frac{d}{dx}X_{i}^{(k)}+\kappa_{i+1}\frac{d}{dx}X_{i+1}^{(k)}, & \textrm{at }x_{i}, i \neq N.
		\end{cases}
	\end{array}
\end{equation}
The formula, subsequently employed to update the interface values for $i = 1,2, \ldots,N-1$, be:
\begin{align}\label{NNWR_discrete_update}
	\Pi_{i}^{(k)} &= \Pi_{i}^{(k-1)}- \begin{bmatrix}
		\theta_{i} I & 0\\
		0 & \phi_{i} I
	\end{bmatrix}
	\left( \Xi_{i}^{(k)}(x_i)+\Xi_{i+1}^{(k)}(x_i)\right)
\end{align}
After evaluating the solution of Dirichlet subproblems outlined in equation \eqref{NNWR_discrete_system} for $i=1,\ldots,N$, we obtain:
\begin{align}\label{NNWR_sol1}
		X^{(k)}_{i} = 
		P\sinh\left((x-x_{i-1})\sqrt{\Lambda/\kappa_{i}}\right)\cosech\left(h_i\sqrt{\Lambda/\kappa_{i}}\right)P^{-1}\Pi_{i}^{(k-1)}\\
		+ 	P\sinh\left((x_{i}-x)\sqrt{\Lambda/\kappa_{i}}\right)\cosech\left(h_i\sqrt{\Lambda/\kappa_{i}}\right)P^{-1}\Pi_{i-1}^{(k-1)}. \nonumber
\end{align}
And using this Dirichlet solutions \eqref{NNWR_sol1}, Neumann subproblems in \eqref{NNWR_discrete_system} give:
\begin{align*}
		\hat{\Xi}^{(k)}_{1}(x) &=  P\sinh\left((x-x_{0})\sqrt{\Lambda/\kappa_{1}}\right)C^{(k-1)}_{1},\\ 
		\hat{\Xi}^{(k)}_{i}(x) &= P\cosh\left((x-x_{i-1})\sqrt{\Lambda/\kappa_{i}}\right)C^{(k-1)}_{i} + P\cosh\left((x_{i}-x)\sqrt{\Lambda/\kappa_{i}}\right)D^{(k-1)}_{i}, \mbox{for $i = 2,\ldots,N-1,$} \nonumber\\
		\hat{\Xi}^{(k)}_{N}(x) &= P\sinh\left((x_{N}-x)\sqrt{\Lambda/\kappa_{N}}
		\right)D^{(k-1)}_{N}. \nonumber
\end{align*}
For the sake of simplifying notations, we consider: $\sigma_{i}:=\sinh\left(a_{i}\sqrt{\Lambda/\kappa_{i}}\right)$ and $\gamma_{i}:=\cosh\left(a_{i}\sqrt{\Lambda/\kappa_{i}}\right)$, and using boundary conditions in Neumann subproblems to evaluate the ode constants, which for $i = 2,\ldots,N-1$ are:
\begin{align*}
		C^{(k-1)}_{1} & = \gamma_{1}^{-1}\left(\left(
		\gamma_{1}\sigma_{1}^{-1} + \sqrt{\frac{\kappa_{2}}{\kappa_{1}}}\gamma_{2}\sigma_{2}^{-1}\right)P^{-1}\Pi_{1}^{(k-1)}
		- \sqrt{\frac{\kappa_{2}}{\kappa_{1}}}\sigma_{2}^{-1}P^{-1}\Pi_{2}^{(k-1)}\right),\\
		C^{(k-1)}_{i} & = \sigma_{i}^{-1}\left(
		\left(\gamma_{i}\sigma_{i}^{-1}+\sqrt{\frac{\kappa_{i+1}}{\kappa_{i}}}\gamma_{i+1}\sigma_{i+1}^{-1}\right)P^{-1}\Pi_{i}^{(k-1)} 
		-\sigma_{i}^{-1}P^{-1}\Pi_{i-1}^{(k-1)}-\sqrt{\frac{\kappa_{i+1}}{\kappa_{i}}}\sigma_{i+1}^{-1}P^{-1}\Pi_{i+1}^{(k-1)}\right), \nonumber\\		
		D^{(k-1)}_{i} & = \sigma_{i}^{-1}\left(
		\left(\gamma_{i}\sigma_{i}^{-1}+\sqrt{\frac{\kappa_{i-1}}{\kappa_{i}}}\gamma_{i-1}\sigma_{i-1}^{-1}\right)P^{-1}\Pi_{i-1}^{(k-1)}
		-\sqrt{\frac{\kappa_{i-1}}{\kappa_{i}}}\sigma_{i-1}^{-1}P^{-1}\Pi_{i-2}^{(k-1)} - \sigma_{i}^{-1}P^{-1}\Pi_{i}^{(k-1)}\right), \nonumber\\		
		D^{(k-1)}_{N} & = \gamma_{N}^{-1}\left(\left(
		\sqrt{\frac{\kappa_{N-1}}{\kappa_{N}}}\gamma_{N-1}\sigma_{N-1}^{-1} + \gamma_{N}\sigma_{N}^{-1}\right)P^{-1}\Pi_{N-1}^{(k-1)}
		-\sqrt{\frac{\kappa_{N-1}}{\kappa_{N}}}\sigma_{N-1}^{-1}P^{-1}\Pi_{N-2}^{(k-1)}\right).\nonumber
\end{align*}
By substituting the values of $\Xi_{i}$ into the updating part \eqref{NNWR_discrete_update} and employing  $\gamma_{i}^{2}-\sigma_{i}^{2} = I$, where $I$ imply identity, we obtain:
\begin{align}
		\Pi_{1}^{(k)} &=\Pi_{1}^{(k-1)}-\begin{bmatrix}
			\theta_{1} I & 0\\
			0 & \phi_{1} I
		\end{bmatrix}
		P\left(
		\left(2I + \sqrt{\frac{\kappa_{1}}{\kappa_{2}}}\gamma_{1}\sigma_{1}^{-1}\gamma_{2}\sigma_{2}^{-1}
		+\sqrt{\frac{\kappa_{2}}{\kappa_{1}}}\sigma_{1}\gamma_{1}^{-1}\gamma_{2}\sigma_{2}^{-1}\right)P^{-1}\Pi_{1}^{(k-1)}\right. \label{NN_1}\\
		& \hspace{10pt}\left.\kern-\nulldelimiterspace + \sigma_{2}^{-1}\left(
		\sqrt{\frac{\kappa_{3}}{\kappa_{2}}} \gamma_{3}\sigma_{3}^{-1} -\sqrt{\frac{\kappa_{2}}{\kappa_{1}}}\sigma_{1}\gamma_{1}^{-1}\right)P^{-1}\Pi_{2}^{(k-1)}
		-\sqrt{\frac{\kappa_{3}}{\kappa_{2}}}\sigma_{2}^{-1}\sigma_{3}^{-1}P^{-1}\Pi_{3}^{(k-1)}\right), \nonumber
\end{align}
\begin{align}
		\Pi_{i}^{(k)} &= \Pi_{i}^{(k-1)} - \begin{bmatrix}
			\theta_{i} I & 0\\
			0 & \phi_{i} I
		\end{bmatrix} 
	    P\left(
		\left(2I + \left(\sqrt{\frac{\kappa_i}{\kappa_{i+1}}}+\sqrt{\frac{\kappa_{i+1}}{\kappa_{i}}}\right)\gamma_{i}\sigma_{i}^{-1}\gamma_{i+1}\sigma_{i+1}^{-1}\right)P^{-1}\Pi_{i}^{(k-1)}\right.\\
		&\hspace{2pt}\left.\kern-\nulldelimiterspace + \sigma_{i+1}^{-1}
		\left(\sqrt{\frac{\kappa_{i+2}}{\kappa_{i+1}}}\gamma_{i+2}\sigma_{i+2}^{-1} - \sqrt{\frac{\kappa_{i+1}}{\kappa_{i}}}\gamma_{i}\sigma_{i}^{-1}\right)P^{-1}\Pi_{i+1}^{(k-1)}
		+ \sigma_{i}^{-1}
		\left(\sqrt{\frac{\kappa_{i-1}}{\kappa_{i}}} \gamma_{i-1}\sigma_{i-1}^{-1}-\sqrt{\frac{\kappa_{i}}{\kappa_{i+1}}}\gamma_{i+1}\sigma_{i+1}^{-1}\right)P^{-1}\Pi_{i-1}^{(k-1)} \right.\nonumber\\
		& \hspace{10pt}\left.\kern-\nulldelimiterspace - \sqrt{\frac{\kappa_{i+2}}{\kappa_{i+1}}}\sigma_{i+1}^{-1}\sigma_{i+2}^{-1}P^{-1}\Pi_{i+2}^{(k-1)}
		-\sqrt{\frac{\kappa_{i-1}}{\kappa_{i}}}\sigma_{i}^{-1}\sigma_{i-1}^{-1}P^{-1}\Pi_{i-2}^{(k-1)}\right), \quad i = 2,\ldots,N-2, \nonumber
\end{align}
\begin{align}
		\Pi_{N-1}^{(k)} &= \Pi_{N-1}^{(k-1)} - \begin{bmatrix}
			\theta_{N-1} I & 0\\
			0 & \phi_{N-1} I
		\end{bmatrix} 
		P\left(
		\left(2I + \sqrt{\frac{\kappa_{N}}{\kappa_{N-1}}}\gamma_{N-1}\sigma_{N-1}^{-1}\gamma_{N}\sigma_{N}^{-1}
		+\sqrt{\frac{\kappa_{N-1}}{\kappa_{N}}}\sigma_{N}\gamma_{N}^{-1}\gamma_{N-1}\sigma_{N-1}^{-1}\right)P^{-1}\Pi_{N-1}^{(k-1)}\right. \label{NN_3}\\
		&\hspace{5pt}\left.\kern-\nulldelimiterspace + \sigma_{N-1}^{-1}\left(
		\sqrt{\frac{\kappa_{N-2}}{\kappa_{N-1}}}\gamma_{N-2}\sigma_{N-2}^{-1} - \sqrt{\frac{\kappa_{N-1}}{\kappa_{N}}}\sigma_{N}\gamma_{N}^{-1}\right)P^{-1}\Pi_{N-2}^{(k-1)} - \sqrt{\frac{\kappa_{N-2}}{\kappa_{N-1}}}\sigma_{N-1}^{-1}\sigma_{N-2}^{-1}P^{-1}\Pi_{N-3}^{(k-1)}\right). \nonumber	
\end{align}
Using $\theta_i = \phi_i = \bar{\theta}_i := 1/(2+\sqrt{\kappa_i/\kappa_{i+1}}+\sqrt{\kappa_{i+1}/\kappa_{i}})$ to eliminate the linear terms in \eqref{NN_1}-\eqref{NN_3}, we have:
\begin{align} 
	\Pi_{1}^{(k)} &= -\bar{\theta}_{1}
	P\left(
	\left(\sqrt{\frac{\kappa_{1}}{\kappa_{2}}}(\gamma_{1}\sigma_{1}^{-1}\gamma_{2}\sigma_{2}^{-1} - I)
	+\sqrt{\frac{\kappa_{2}}{\kappa_{1}}}(\sigma_{1}\gamma_{1}^{-1}\gamma_{2}\sigma_{2}^{-1} - I)\right)P^{-1}\Pi_{1}^{(k-1)}\right. \label{NN_4}\\
	& \hspace{10pt}\left.\kern-\nulldelimiterspace + \sigma_{2}^{-1}\left(
	\sqrt{\frac{\kappa_{3}}{\kappa_{2}}} \gamma_{3}\sigma_{3}^{-1} -\sqrt{\frac{\kappa_{2}}{\kappa_{1}}}\sigma_{1}\gamma_{1}^{-1}\right)P^{-1}\Pi_{2}^{(k-1)}
	-\sqrt{\frac{\kappa_{3}}{\kappa_{2}}}\sigma_{2}^{-1}\sigma_{3}^{-1}P^{-1}\Pi_{3}^{(k-1)}\right), \nonumber
\end{align}
\begin{align}
	\Pi_{i}^{(k)} &= -\bar{\theta}_{i}
	P\left(
	 \left(\sqrt{\frac{\kappa_i}{\kappa_{i+1}}}+\sqrt{\frac{\kappa_{i+1}}{\kappa_{i}}}\right)\left(\gamma_{i}\sigma_{i}^{-1}\gamma_{i+1}\sigma_{i+1}^{-1} - I\right)P^{-1}\Pi_{i}^{(k-1)}\right.\\
	&\hspace{2pt}\left.\kern-\nulldelimiterspace + \sigma_{i+1}^{-1}
	\left(\sqrt{\frac{\kappa_{i+2}}{\kappa_{i+1}}}\gamma_{i+2}\sigma_{i+2}^{-1} - \sqrt{\frac{\kappa_{i+1}}{\kappa_{i}}}\gamma_{i}\sigma_{i}^{-1}\right)P^{-1}\Pi_{i+1}^{(k-1)}
	+ \sigma_{i}^{-1}
	\left(\sqrt{\frac{\kappa_{i-1}}{\kappa_{i}}} \gamma_{i-1}\sigma_{i-1}^{-1}-\sqrt{\frac{\kappa_{i}}{\kappa_{i+1}}}\gamma_{i+1}\sigma_{i+1}^{-1}\right)P^{-1}\Pi_{i-1}^{(k-1)} \right.\nonumber\\
	& \hspace{10pt}\left.\kern-\nulldelimiterspace - \sqrt{\frac{\kappa_{i+2}}{\kappa_{i+1}}}\sigma_{i+1}^{-1}\sigma_{i+2}^{-1}P^{-1}\Pi_{i+2}^{(k-1)}
	-\sqrt{\frac{\kappa_{i-1}}{\kappa_{i}}}\sigma_{i}^{-1}\sigma_{i-1}^{-1}P^{-1}\Pi_{i-2}^{(k-1)}\right), \quad i = 2,\ldots,N-2, \nonumber
\end{align}
\begin{align}
	\Pi_{N-1}^{(k)} &= - \bar{\theta}_{N-1}
	P\left(
	\left(\sqrt{\frac{\kappa_{N}}{\kappa_{N-1}}}(\gamma_{N-1}\sigma_{N-1}^{-1}\gamma_{N}\sigma_{N}^{-1} - I)
	+\sqrt{\frac{\kappa_{N-1}}{\kappa_{N}}}(\sigma_{N}\gamma_{N}^{-1}\gamma_{N-1}\sigma_{N-1}^{-1} - I)\right)P^{-1}\Pi_{N-1}^{(k-1)}\right. \label{NN_6}\\
	&\hspace{5pt}\left.\kern-\nulldelimiterspace + \sigma_{N-1}^{-1}\left(
	\sqrt{\frac{\kappa_{N-2}}{\kappa_{N-1}}}\gamma_{N-2}\sigma_{N-2}^{-1} - \sqrt{\frac{\kappa_{N-1}}{\kappa_{N}}}\sigma_{N}\gamma_{N}^{-1}\right)P^{-1}\Pi_{N-2}^{(k-1)} - \sqrt{\frac{\kappa_{N-2}}{\kappa_{N-1}}}\sigma_{N-1}^{-1}\sigma_{N-2}^{-1}P^{-1}\Pi_{N-3}^{(k-1)}\right). \nonumber	
\end{align}	
Before going to the convergence result we need the subsequent lemma for estimation.
\begin{lemma} \label{NNWR_lam_1}
	For $l_0,l_1,l_2 >0$, and $z \in \mathbb{C}$ such that $Re(z) \geq x_0 > 0$   we have the following estimates:
	\begin{enumerate}
		\item [(i)] 
		\begin{align*}
				\left|\frac{\cosh(l_1z)\cosh(l_2z)}{\sinh(l_1z)\sinh(l_2z)}-1\right| \leq
				\coth(l_1x_0)\coth(l_2x_0) - 1
		\end{align*}
		\item [(ii)] 
		\begin{align*}
				\left|\frac{\sinh(l_1z)\cosh(l_2z)}{\sinh(l_2z)\cosh(l_1z)}-1\right|
				\leq \coth(l_1x_0)\coth(l_2x_0) - 1 
		\end{align*}
		\item [(iii)] 
		\begin{align*}
				\left|\frac{1}{\sinh(l_1z)\sinh(l_2z)}\right| 
				\leq \cosech(l_1x_0)\cosech(l_2x_0)
		\end{align*}
		\item [(iv)]
		\begin{align*}
				\left|\frac{1}{\sinh(l_0z)}\left(\rho_1\frac{\cosh(l_1z)}{\sinh(l_1z)}-\rho_2\frac{\cosh(l_2z)}{\sinh(l_2z)}\right)\right| 
			   \leq 
			   \cosech(l_0x_0)\left(\rho_1\coth(l_1x_0) + \rho_2\coth(l_2x_0)\right)
		\end{align*}
		\item [(v)] 
		\begin{align*}
				\left|\frac{1}{\sinh(l_0z)}\left(\rho_1\frac{\cosh(l_1z)}{\sinh(l_1z)}-\rho_2\frac{\sinh(l_2z)}{\cosh(l_2z)}\right)\right| 
				\leq \cosech(l_0x_0)\left(\rho_1\coth(l_1x_0) + \rho_2\coth(l_2x_0)\right)
		\end{align*}
	\end{enumerate}	
\end{lemma}	
\begin{proof}
	The procedure for proving all five parts is similar, with the only distinction being the combination of exponentials. 
\begin{enumerate}
	\item [(i)] Rewriting the hyperbolic functions in terms of exponential than using the inequality yield
	\begin{align} 
			&\left|\frac{\cosh(l_1z)\cosh(l_2z)}{\sinh(l_1z)\sinh(l_2z)}-1\right| \label{mid_11}\\
			&= \left|\frac{\cosh((l_1-l_2)z)}{\sinh(l_1z)\sinh(l_2z)}\right| \nonumber\\
			&\leq 2\left|\frac{\exp(-2l_1z)}{(1-\exp(-2l_1z))(1-\exp(-2l_2z))}\right| +  2\left|\frac{\exp(-2l_2z)}{(1-\exp(-2l_1z))(1-\exp(-2l_2z))}\right|. \nonumber
	\end{align}
	Now for the first part, using the Lemma \ref{estimate} we have
	\begin{align} 
		\left|\frac{\exp(-2l_1z)}{(1-\exp(-2l_1z))(1-\exp(-2l_2z))}\right|  
		& \leq  \left|\frac{\exp(-l_1z)}{(1-\exp(-2l_1z))}\right| \left|\frac{\exp(-l_1z)}{(1-\exp(-2l_2z))}\right| \label{mid_12}\\
		& \leq  \frac{\exp(-l_1x_0)}{(1-\exp(-2l_1x_0))}\frac{\exp(-l_1x_0)}{(1-\exp(-2l_2x_0))} \nonumber\\
		 &= \frac{\exp(-2l_1x_0)}{(1-\exp(-2l_1x_0))(1-\exp(-2l_2x_0))}. \nonumber
	\end{align}
	In similar manner from the second part we obtain the estimate
	\begin{align}
		\left|\frac{\exp(-2l_2z)}{(1-\exp(-2l_1z))(1-\exp(-2l_2z))}\right| 
		\leq \frac{\exp(-2l_2x_0)}{(1-\exp(-2l_1x_0))(1-\exp(-2l_2x_0))}. \label{mid_13}
	\end{align}
	Therefore combining \eqref{mid_12} and \eqref{mid_13} in \eqref{mid_11} provides
	\begin{align*}
		\left|\frac{\cosh(l_1z)\cosh(l_2z)}{\sinh(l_1z)\sinh(l_2z)}-1\right| 
		\leq \frac{\cosh(l_1x_0)\cosh(l_2x_0)}{\sinh(l_1x_0)\sinh(l_2x_0)}-1.
	\end{align*}
	Proofs for parts (ii) and (iii) follow a similar approach as part (i), and are therefore omitted.
	\item [(iv)] 
	By employing the triangle inequality gives
	\begin{align*}
			\left|\frac{1}{\sinh(l_0z)}\left(\rho_1\frac{\cosh(l_1z)}{\sinh(l_1z)}-\rho_2\frac{\cosh(l_2z)}{\sinh(l_2z)}\right)\right|
			 \leq \rho_1\left|\frac{\cosh(l_1z)}{\sinh(l_0z) \sinh(l_1z)}\right| + \rho_2\left|\frac{\cosh(l_2z)}{\sinh(l_0z)\sinh(l_2z)}\right|.
	\end{align*}
	Rewrite the hyperbolic expressions in exponential form, then applying Lemma \ref{estimate} to the initial part yields
	\begin{align*}
			\rho_1\left|\frac{\cosh(l_1z)}{\sinh(l_0z) \sinh(l_1z)}\right| &= 2\rho_1\left|\frac{\exp(-l_0z)+\exp(-(l_0+2l_1)z)}{(1-\exp(-2l_0z))(1-\exp(-2l_1z))}\right| \\
			&\leq 2\rho_1\frac{\exp(-l_0x_0)+\exp(-(l_0+2l_1)x_0)}{(1-\exp(-2l_0x_0))(1-\exp(-2l_1x_0))} \\
			& \leq \rho_1\frac{\cosh(l_1x_0)}{\sinh(l_0x_0) \sinh(l_1x_0)}.
	\end{align*}
	In similar manner we obtain the second part's estimate, which is
	\begin{align*}
		\rho_2\left|\frac{\cosh(l_2z)}{\sinh(l_0z) \sinh(l_2z)}\right| 
		\leq \rho_2\frac{\cosh(l_2x_0)}{\sinh(l_0x_0) \sinh(l_2x_0)}.
	\end{align*}
	Therefore combining both parts we obtain
	\begin{align*}
		\left|\frac{1}{\sinh(l_0z)}\left(\rho_1\frac{\cosh(l_1z)}{\sinh(l_1z)}-\rho_2\frac{\cosh(l_2z)}{\sinh(l_2z)}\right)\right|
		\leq \frac{1}{\sinh(l_0x_0)}\left(\rho_1\frac{\cosh(l_1x_0)}{\sinh(l_1x_0)}+\rho_2\frac{\cosh(l_2x_0)}{\sinh(l_2x_0)}\right).
	\end{align*}
	\item [(v)] The proof follows a similar approach as part (iv) and is consequently omitted.
\end{enumerate}
\hfill\end{proof}

Based on the aforementioned estimates, we are now able to introduce the principal convergent outcome of the NNWR algorithm.
\begin{theorem}[Convergence of NNWR] \label{NNWR_th}
		Let time discretize combine state and adjoint  matrix $\mathbb{L}$ in \eqref{NNWR_discrete_system} be diagonalizable by the eigenvector matrix $P$ corresponding to the eigenvalues $\{\lambda_i: i = 1,\ldots,2n+2\}$.
	If $\lambda = \min(\Re(\sqrt{\lambda_i}))$, then the interface error term of the NNWR algorithm in \eqref{NNWR_state}-\eqref{NNWR_adjoint} in 1D for relaxation parameters $\theta_i = \phi_{i} =\bar{\theta}_i= 1/(2+\sqrt{\kappa_i/\kappa_{i+1}}+\sqrt{\kappa_{i+1}/\kappa_{i}})$ satisfies the following bound
		\begin{equation*}
		\max_{1 \leq i \leq N-1}||\Pi_{i}^{(k)}||_{\infty} \leq \bar{d}^k \cosech^{2k}\left(\frac{h}{2} \lambda\right) \|P\|_{\infty} \|P^{-1}\|_{\infty} \max_{1 \leq j \leq N-1}||\Pi_{i}^{(0)}||_{\infty},
	\end{equation*}
	where $h = \min_{1 \leq i\leq N} a_i/\sqrt{\kappa_i}$ and $\bar{d} = \max_{1\leq i\leq N-1}\bar{\theta}_i d_i$, the values of $d_i$ are given in \eqref{NNWR_5_5} and \eqref{NNWR_5_6}.
\end{theorem}
\begin{proof}
	First we use the same argument as in the proof of Theorem \ref{Th_2} to show that $\lambda = \min(\Re(\sqrt{\lambda_i}))$ is well defined. Then left multiplying on both side with $P^{-1}$ of the update terms in \eqref{NN_4}-\eqref{NN_6} to obtain 
	\begin{align} \label{NNWR_1}
			P^{-1}\Pi_{1}^{(k)} &= -\bar{\theta}_1\left(T_{1,1}P^{-1}\Pi_{1}^{(k-1)} + T_{1,2}P^{-1}\Pi_{2}^{(k-1)} - T_{1,3}P^{-1}\Pi_{3}^{(k-1)}\right),  \\
			P^{-1}\Pi_{i}^{(k)} &= -\bar{\theta}_i\left(T_{i,i}P^{-1}\Pi_{i}^{(k-1)} + T_{i,i+1}P^{-1}\Pi_{i+1}^{(k-1)} + \hat{t}_{i,i-1}P^{-1}\Pi_{i-1}^{(k-1)} - T_{i,i+2}P^{-1}\Pi_{i+2}^{(k-1)} - T_{i,i-2}P^{-1}\Pi_{i-2}^{(k-1)}\right), \nonumber\\
			P^{-1}\Pi_{N-1}^{(k)} &= -\bar{\theta}_{N-1}\left(T_{N-1,N-1}P^{-1}\Pi_{N-1}^{(k-1)} + T_{N-1,N-2}P^{-1}\Pi_{N-2}^{(k-1)} - T_{N-1,N-3}P^{-1}\Pi_{N-3}^{(k-1)}\right), \nonumber
	\end{align}
	where the weights are given by
	\begin{gather} \label{NNWR_2}
			T_{1,1} = \left(\sqrt{\frac{\kappa_{1}}{\kappa_{2}}}(\gamma_{1}\sigma_{1}^{-1}\gamma_{2}\sigma_{2}^{-1} - I)
			+\sqrt{\frac{\kappa_{2}}{\kappa_{1}}}(\sigma_{1}\gamma_{1}^{-1}\gamma_{2}\sigma_{2}^{-1} - I)\right), \\
			T_{1,2} = \sigma_{2}^{-1}\left(
			\sqrt{\frac{\kappa_{3}}{\kappa_{2}}} \gamma_{3}\sigma_{3}^{-1} -\sqrt{\frac{\kappa_{2}}{\kappa_{1}}}\sigma_{1}\gamma_{1}^{-1}\right), \,
			T_{1,3} = \sqrt{\frac{\kappa_{3}}{\kappa_{2}}}\sigma_{2}^{-1}\sigma_{3}^{-1}, \nonumber 
	\end{gather}
	for $i = 2,\cdots,N-2$,
	\begin{gather} 
			T_{i,i} = \left(\sqrt{\frac{\kappa_i}{\kappa_{i+1}}}+\sqrt{\frac{\kappa_{i+1}}{\kappa_{i}}}\right)\left(\gamma_{i}\sigma_{i}^{-1}\gamma_{i+1}\sigma_{i+1}^{-1} - I\right), \,
			T_{i,i+1} = \sigma_{i+1}^{-1}
			\left(\sqrt{\frac{\kappa_{i+2}}{\kappa_{i+1}}}\gamma_{i+2}\sigma_{i+2}^{-1} - \sqrt{\frac{\kappa_{i+1}}{\kappa_{i}}}\gamma_{i}\sigma_{i}^{-1}\right), \\
			T_{i,i-1} = \sigma_{i}^{-1}
			\left(\sqrt{\frac{\kappa_{i-1}}{\kappa_{i}}} \gamma_{i-1}\sigma_{i-1}^{-1}-\sqrt{\frac{\kappa_{i}}{\kappa_{i+1}}}\gamma_{i+1}\sigma_{i+1}^{-1}\right), \,
			T_{i,i+2} = \sqrt{\frac{\kappa_{i+2}}{\kappa_{i+1}}}\sigma_{i+1}^{-1}\sigma_{i+2}^{-1}, \,
			T_{i,i-2} = \sqrt{\frac{\kappa_{i-1}}{\kappa_{i}}}\sigma_{i}^{-1}\sigma_{i-1}^{-1}, \nonumber
	\end{gather}
	and
	\begin{gather} \label{NNWR_2_2}
			T_{N-1,N-1} = \left(\sqrt{\frac{\kappa_{N}}{\kappa_{N-1}}}(\gamma_{N-1}\sigma_{N-1}^{-1}\gamma_{N}\sigma_{N}^{-1} - I)
			+ \sqrt{\frac{\kappa_{N-1}}{\kappa_{N}}}(\sigma_{N}\gamma_{N}^{-1}\gamma_{N-1}\sigma_{N-1}^{-1} - I)\right), \\
			T_{N-1,N-2} = \sigma_{N-1}^{-1}\left(
			\sqrt{\frac{\kappa_{N-2}}{\kappa_{N-1}}}\gamma_{N-2}\sigma_{N-2}^{-1} - \sqrt{\frac{\kappa_{N-1}}{\kappa_{N}}}\sigma_{N}\gamma_{N}^{-1}\right),\,
			T_{N-1,N-3} = \sqrt{\frac{\kappa_{N-2}}{\kappa_{N-1}}}\sigma_{N-1}^{-1}\sigma_{N-2}^{-1}. \nonumber
	\end{gather}
	These weights $T_{i,j}$ in \eqref{NNWR_2}-\eqref{NNWR_2_2} are in the same form as in Lemma \ref{NNWR_lam_1}. Therefore, using the notation $h_j = \frac{a_j}{\sqrt{\kappa_j}}: j = 1,\ldots,N$, we obtain the subsequent estimations for $i = 1,2,\cdots,N-1$:
	\begin{align*}
			||T_{i,i}(.)||_{\infty} &\leq \left(\sqrt{\frac{\kappa_i}{\kappa_{i+1}}}+\sqrt{\frac{\kappa_{i+1}}{\kappa_i}}\right) 
			\left(\coth(h_i \lambda) \coth(h_{i+1} \lambda) - 1\right)=:w_{i,i},\\
	\end{align*}
	for $i = 1,\cdots,N-2,$
	\begin{align*}
			||T_{i,i+1}(.)||_{\infty} &\leq \cosech(h_{i+1}\lambda)\left(\sqrt{\frac{\kappa_{i+2}}{\kappa_{i+1}}}\coth(h_{i+2}\lambda) + \sqrt{\frac{\kappa_{i+1}}{\kappa_i}}\coth(h_i \lambda)\right)=:w_{i,i+1}, \\
			||T_{i,i+2}(.)||_{\infty} &\leq \sqrt{\frac{\kappa_{i+2}}{\kappa_{i+1}}} \cosech(h_{i+1}\lambda) \cosech(h_{i+2}\lambda)=:w_{i,i+2},
	\end{align*}
	for $i = 2,\cdots,N-1,$
	\begin{align*}
			||T_{i,i-1}(.)||_{\infty} &\leq \cosech(h_{i}\lambda)\left(\sqrt{\frac{\kappa_{i-1}}{\kappa_{i}}}\coth(h_{i-1}\lambda) + \sqrt{\frac{\kappa_{i}}{\kappa_{i+1}}}\coth(h_{i+1} \lambda)\right)=:w_{i,i-1}, \\
			||T_{i,i-2}(.)||_{\infty} &\leq \sqrt{\frac{\kappa_{i-1}}{\kappa_{i}}} \cosech(h_{i-1}\lambda) \cosech(h_{i}\lambda)=:w_{i,i-2}.
	\end{align*}
	 Now taking the supremum norm on both side in \eqref{NNWR_1} we have
	\begin{align} \label{NNWR_5}
		||P^{-1}\Pi_{1}^{(k)}||_{\infty} &\leq c_1 \bar{\theta}_1 \max_{1 \leq j \leq N-1}||P^{-1}\Pi_{j}^{(k-1)}||_{\infty}, \nonumber\\
		||P^{-1}\Pi_{i}^{(k)}||_{\infty} &\leq c_i \bar{\theta}_i \max_{1 \leq j \leq N-1}||P^{-1}\Pi_{j}^{(k-1)}||_{\infty},  \\
		||P^{-1}\Pi_{N-1}^{(k)}||_{\infty} &\leq c_{N-1} \bar{\theta}_{N-1} \max_{1 \leq j \leq N-1}||P^{-1}\Pi_{j}^{(k-1)}||_{\infty}. \nonumber
	\end{align}
	where
	\begin{align} \label{NNWR_5_1}
		c_1 &=  w_{1,1} + w_{1,1+1} + w_{1,1+2}, \nonumber\\
		c_i &= w_{i,i} + w_{i,i-1} + w_{i,i+1} + w_{i,i-2} + w_{i,i+2}, \quad i = 2,\cdots,N-2,\\
		c_{N-1} &= w_{N-1,N-1} + w_{N-1,N-2} + w_{N-1,N-3}. \nonumber
	\end{align}
	Now we can find the weight bounds in terms of minimum subdomain length defined as $h := \min_{1 \leq i\leq N} h_i$. Utilizing the monotonically decreasing nature of the $\coth$ and $\cosech$ functions with positive arguments. The estimates are as follows:
	\begin{align} \label{NNWR_5_3}
		w_{i,i} 
		&= \left(\sqrt{\frac{\kappa_i}{\kappa_{i+1}}}+\sqrt{\frac{\kappa_{i+1}}{\kappa_i}}\right) 
		\left(\coth(h_i \lambda) \coth(h_{i+1} \lambda) - 1\right) \\
		&\leq \left(\sqrt{\frac{\kappa_i}{\kappa_{i+1}}}+\sqrt{\frac{\kappa_{i+1}}{\kappa_i}}\right)
		\left(\coth^2(h \lambda) - 1\right) \nonumber\\
		&= \left(\sqrt{\frac{\kappa_i}{\kappa_{i+1}}}+\sqrt{\frac{\kappa_{i+1}}{\kappa_i}}\right) \cosech^2(h\lambda), \nonumber
	\end{align}
	similarly,
	\begin{align} \label{NNWR_5_4}
	w_{i,i+1} 
	&\leq \left(\sqrt{\frac{\kappa_{i+2}}{\kappa_{i+1}}} + \sqrt{\frac{\kappa_{i+1}}{\kappa_i}}\right)\cosech(h\lambda) \coth(h\lambda), \,
	w_{i,i+2}
	\leq \sqrt{\frac{\kappa_{i+2}}{\kappa_{i+1}}} \cosech^2(h\lambda),\\
	w_{i,i-1}
	&\leq \left(\sqrt{\frac{\kappa_{i-1}}{\kappa_{i}}} + \sqrt{\frac{\kappa_{i}}{\kappa_{i+1}}}\right)\cosech(h\lambda) \cot(h\lambda), \,
	w_{i,i-2}
	\leq \sqrt{\frac{\kappa_{i-1}}{\kappa_{i}}} \cosech^2(h\lambda). \nonumber
	\end{align}
    Therefore using the estimates from \eqref{NNWR_5_3}-\eqref{NNWR_5_4} in \eqref{NNWR_5_1} we have
    \begin{align} \label{NNWR_5_5}
    	c_1 
    	&\leq \left(\sqrt{\frac{\kappa_1}{\kappa_{2}}}+\sqrt{\frac{\kappa_{2}}{\kappa_1}} + \sqrt{\frac{\kappa_{3}}{\kappa_{2}}}\right)\cosech(h\lambda) (\cosech(h\lambda) + \coth(h\lambda)) \\
    	&= \left(\sqrt{\frac{\kappa_1}{\kappa_{2}}}+\sqrt{\frac{\kappa_{2}}{\kappa_1}} + \sqrt{\frac{\kappa_{3}}{\kappa_{2}}}\right)\frac{2\cosh^2(\frac{h}{2} \lambda)}{\sinh^2(h\lambda)} \nonumber\\
    	&= \left(\sqrt{\frac{\kappa_1}{\kappa_{2}}}+\sqrt{\frac{\kappa_{2}}{\kappa_1}} + \sqrt{\frac{\kappa_{3}}{\kappa_{2}}}\right)\frac{\cosech^2(\frac{h}{2} \lambda)}{2} 
    	 =: d_1 \cosech^2\left(\frac{h}{2} \lambda\right), \nonumber
    \end{align}
    and similarly 
    \begin{align} \label{NNWR_5_6}
    	c_i &= \left(\sqrt{\frac{\kappa_i}{\kappa_{i+1}}}+\sqrt{\frac{\kappa_{i+1}}{\kappa_i}} + \sqrt{\frac{\kappa_{i+2}}{\kappa_{i+1}}} + \sqrt{\frac{\kappa_{i-1}}{\kappa_{i}}}\right) \frac{\cosech^2(\frac{h}{2}\lambda)}{2}
    	=: d_i \cosech^2\left(\frac{h}{2} \lambda\right),\\
    	c_{N-1} &= \left(\sqrt{\frac{\kappa_{N-1}}{\kappa_{N}}}+\sqrt{\frac{\kappa_{N}}{\kappa_{N-1}}} + \sqrt{\frac{\kappa_{N-2}}{\kappa_{N-1}}}\right)\frac{\cosech^2(\frac{h}{2} \lambda)}{2}
    	=: d_{N-1} \cosech^2\left(\frac{h}{2} \lambda\right). \nonumber
    \end{align}
	Therefore using the bounds of $c_i$ in \eqref{NNWR_5} we have the final result
	\begin{equation*}
		\max_{1 \leq i \leq N-1}||\Pi_{i}^{(k)}||_{\infty} \leq \bar{d}^k \cosech^{2k}\left(\frac{h}{2} \lambda\right) \|P\|_{\infty} \|P^{-1}\|_{\infty} \max_{1 \leq j \leq N-1}||\Pi_{i}^{(0)}||_{\infty},
	\end{equation*}
	where $\bar{d} := \max_{1\leq i\leq N-1}\bar{\theta}_i d_i$.	
	\hfill\end{proof}

%% file: Numerical.tex
\section{Numerical Experiments}\label{Numerical_ex}
We are examining the model problem described by equations \eqref{model_problem_3}-\eqref{model_problem_4}, with zero initial condition and zero forcing term for the state equation. In this setup, we select 1D domain $\Omega = [-1,1]$, time span $T = 1$, generalized diffusion coefficient $\kappa = 1$, and regularization parameter $\sigma = 10^{-6}$. Homogeneous Dirichlet boundary conditions are applied to the boundary $\partial\Omega$. Our target solution follows a similar form to that presented in the article \cite{ciaramella2024convergence}, specifically:

\begin{equation}\label{num_eqn1}
	y_Q(x,t) = (1+t)\sin(\pi x) \left(\exp(-8(x-1)^2) + \exp(-8(x+1)^2) - \exp(-8) - \exp(-72)\right).
\end{equation}
However, we make some adjustments as the head-tail coupling condition in time is not required.

\begin{figure}
	\centering
	\includegraphics[width=0.30\textwidth]{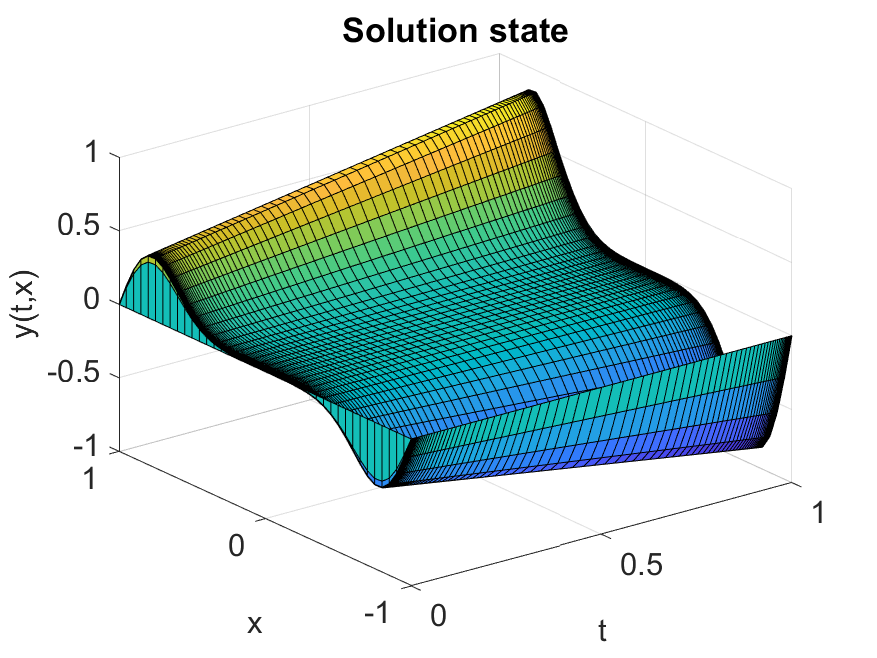}
	\includegraphics[width=0.30\textwidth]{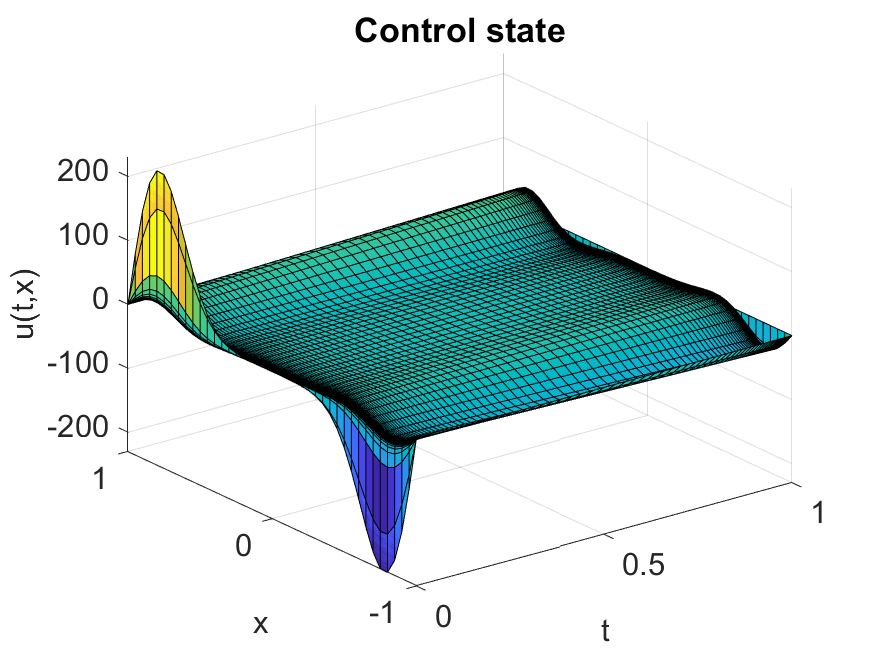}		
	\includegraphics[width=0.30\textwidth]{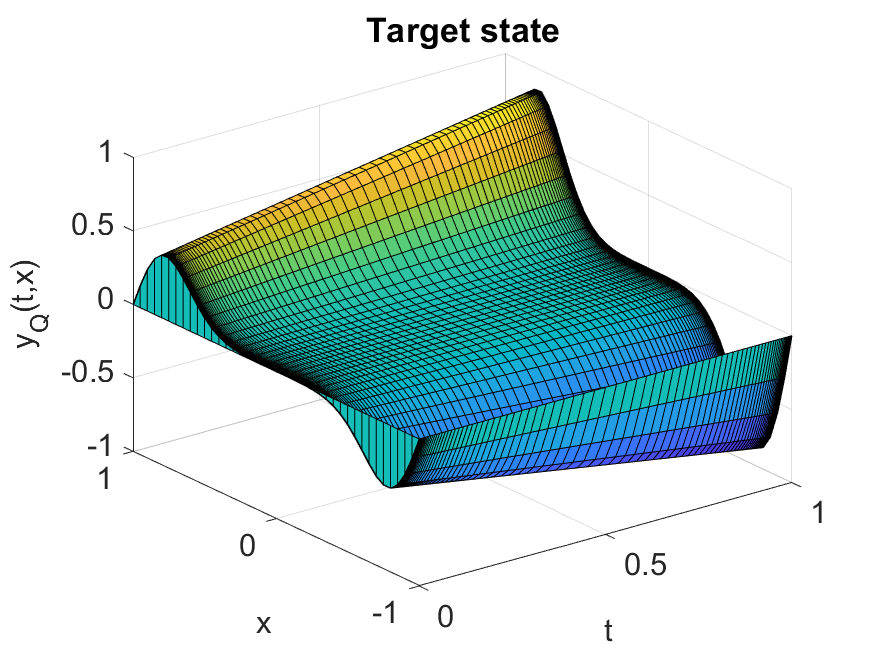}
	\caption{Monodomain solutions profile of solution state, control state, and target state.}
	\label{NumFig00}
\end{figure}

For sub-diffusion scenarios, we utilize a central finite difference approach for spatial calculations and implement an L1 scheme on a both-sided graded mesh to model the fractional time derivative. The test is executed until the relative error of the Dirichlet trace becomes below the threshold  of $10^{-10}$. We set the spatial grid size to $\Delta x = 0.05$ and the temporal grid size to $\Delta t = 0.01$. With the fractional order $\alpha = 0.5$, we initially derive the solution, control, and target profiles as depicted in Figure~\ref{NumFig00}. The figure illustrates that the control component exerts more influence near the sharp change at time $t=0$, due to the zero initial condition.

\subsection{DNWR Algorithm}
\begin{figure}
	\centering
	\includegraphics[width=0.30\textwidth]{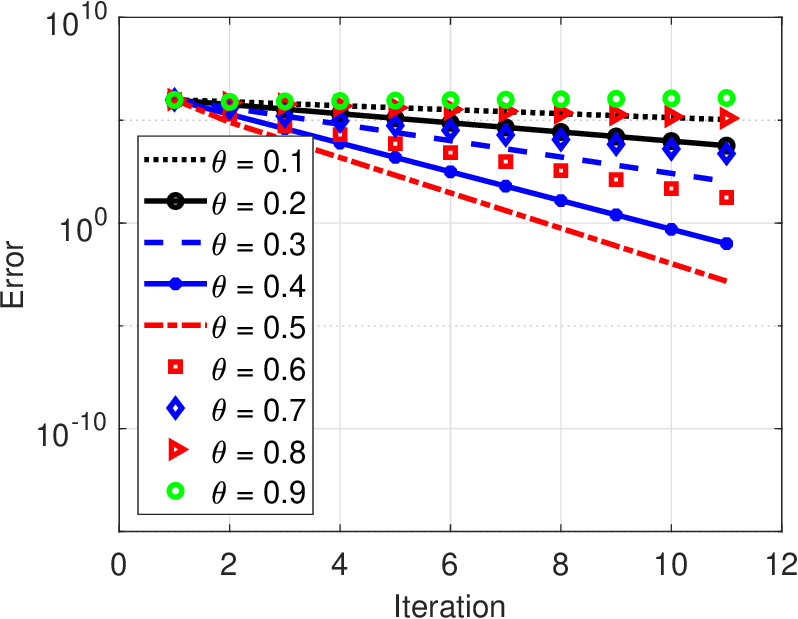}
	\includegraphics[width=0.30\textwidth]{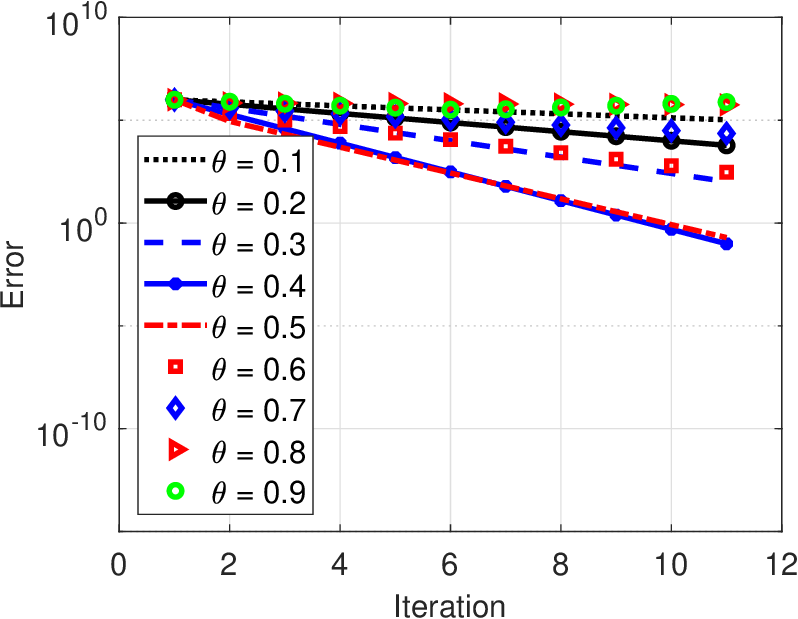}		
	\includegraphics[width=0.30\textwidth]{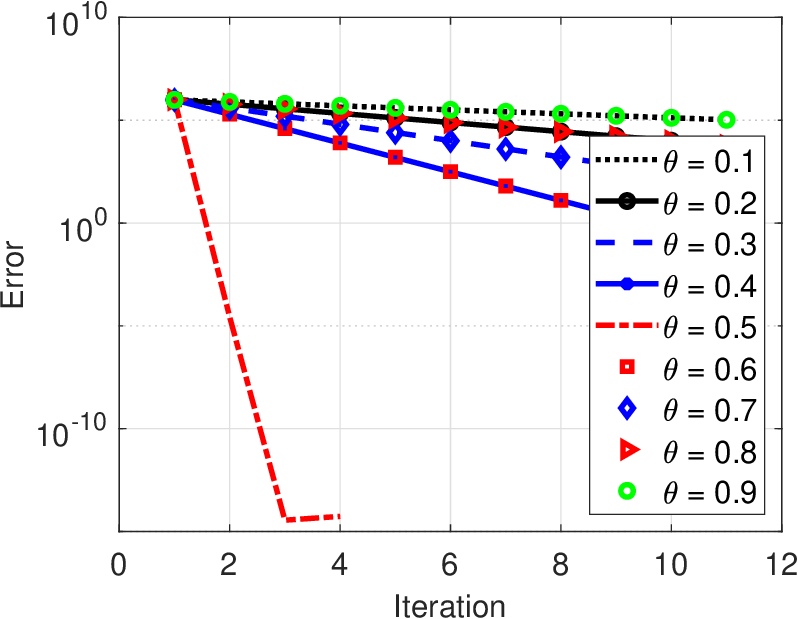}
	\caption{DNWR: The numerical convergence rate varies across different values of $\theta$, while $h_1<h_2$ and $T=1$ remain fixed, considering a fractional order $\alpha = 0.3$. On the left, we employ a uniform mesh in time; in the middle, a one-sided graded mesh in time is utilized; and on the right, both-sided graded mesh in time is applied.}
	\label{NumFig01}
\end{figure}

\begin{figure}
	\centering
	\includegraphics[width=0.30\textwidth]{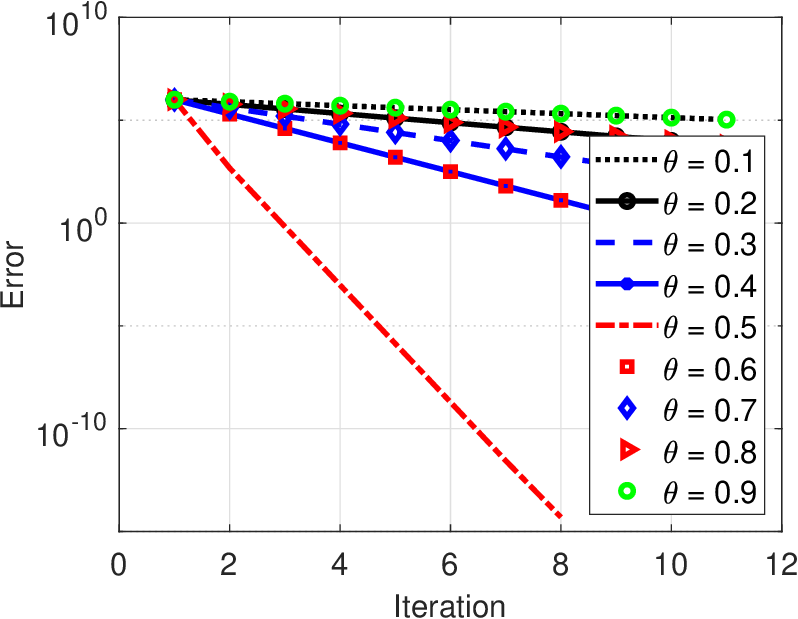}
	\includegraphics[width=0.30\textwidth]{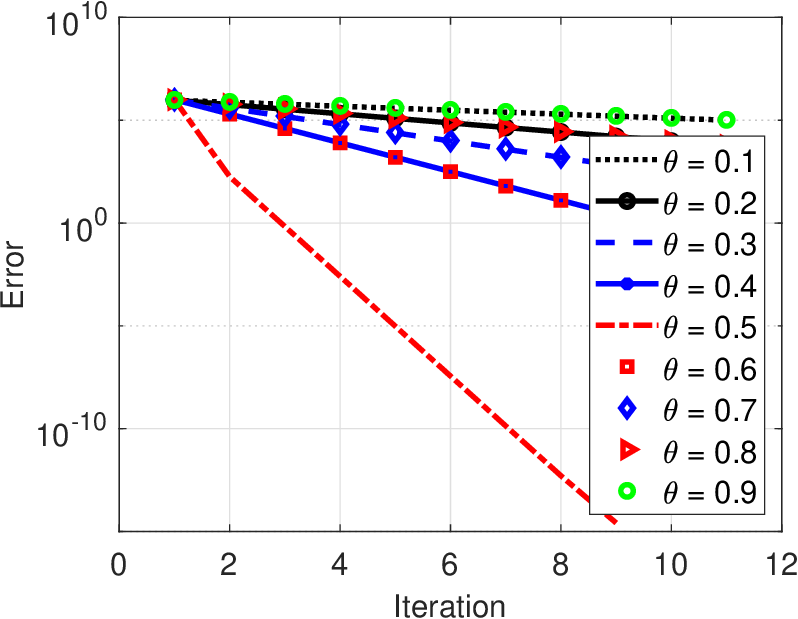}		
	\includegraphics[width=0.30\textwidth]{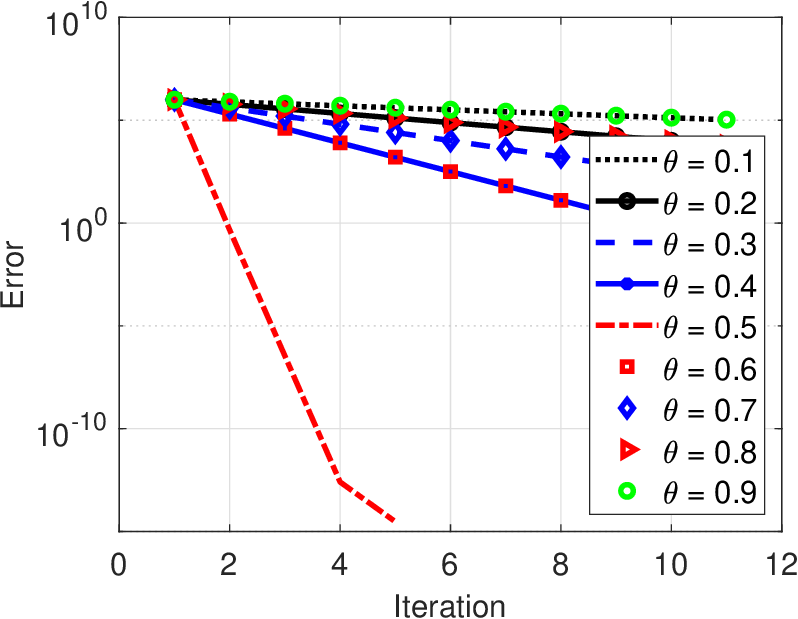}
	\caption{DNWR: The numerical convergence rate varies across different values of $\theta$, while $h_1<h_2$ and $T=1$ remain fixed, considering a fractional order $\alpha = 0.8$. On the left, we employ a uniform mesh in time; in the middle, a one-sided graded mesh in time is utilized; and on the right, both-sided graded mesh in time is applied.}
	\label{NumFig02}
\end{figure}

In Theorem \ref{Th_2}, for the DNWR case, we have selected the relaxation parameter $\theta = \phi = 1/(1+\sqrt{\kappa_1/\kappa_2})$ to eliminate the linear terms from the error component. Now, we substantiate our choice numerically using Figures \ref{NumFig01} and \ref{NumFig02}. In Figure \ref{NumFig01}, we vary $\theta = \phi$ within the range $(0,1)$ and observe the convergence rate as iterations progress. Since the convergence rate relies on the subdomain length and the diffusion coefficient within each subdomain, we set $h_1 = 0.5$, $h_2 = 1.5$, and $\kappa_1 = \kappa_2 = 1$. From Figure \ref{NumFig01}, it's evident that at $\theta = \phi = 0.5$, for a fractional order of $\alpha = 0.3$, the convergence rate surpasses others for all three mesh discretization categories: uniform, one-sided, and both-sided graded mesh. Here one-sided graded means the usual grading technique used for time-fractional PDEs. Similarly, we conducted the same analysis for Figure \ref{NumFig02} but with a fractional order of $\alpha = 0.8$, yielding identical results; namely, $\theta = 0.5$ yields the superior convergence rate.

Let's start by defining the convergence factor $\rho$ as follows:
\begin{equation*}
	\rho := \frac{\sqrt{\kappa_{1}}}{\sqrt{\kappa_1} + \sqrt{\kappa_2}} \left\|\frac{\sinh((b-a)\sqrt{\Lambda})}{\sinh(a\sqrt{\Lambda}) \cosh(b\sqrt{\Lambda})}\right\|_{\infty}.
\end{equation*}
As per \eqref{DNWRsol4}, we understand that the convergence rate of the DNWR algorithm is directly influenced by $\rho$, a scalar quantity dependent on the regularization parameter $\sigma$ and the length of subdomains. Therefore, by analyzing how $\rho$ changes with various factors, we can gain valuable insights into the behavior of the algorithm. In the upcoming set of experiments, we'll delve into this analysis.

\begin{figure}
	\centering
	\includegraphics[width=0.40\textwidth]{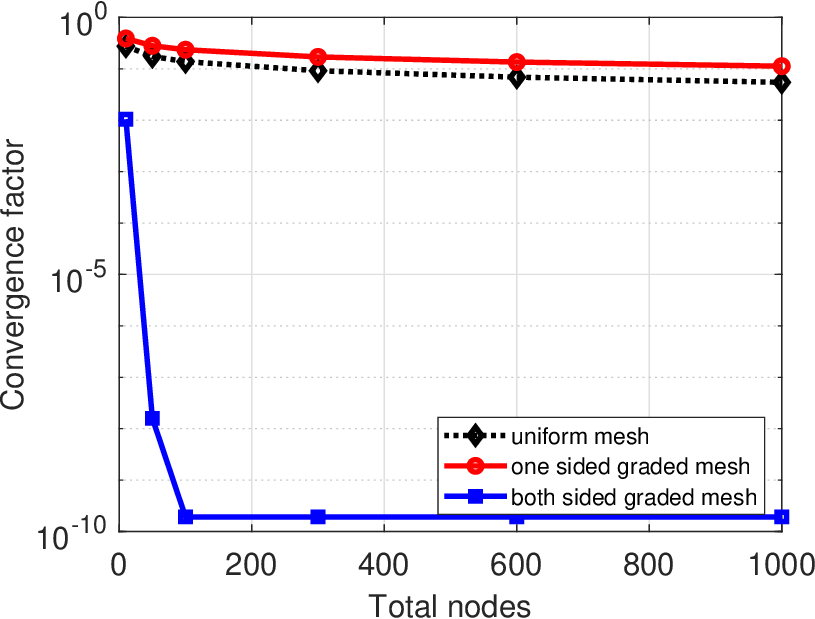}
	\includegraphics[width=0.40\textwidth]{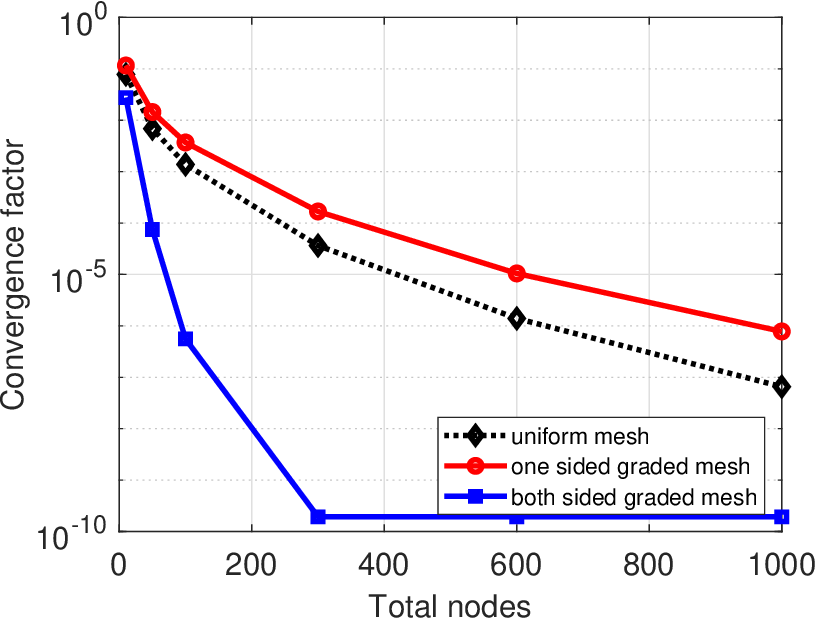}
	\caption{DNWR: As the number of time nodes increases within a fixed time window of $T = 1$ and a constant regularization parameter of $\sigma = 10^{-6}$, the convergence factor varies. On the left side, with a fractional order of $\alpha = 0.3$, while on the right side, with a fractional order of $\alpha = 0.8$.}
	\label{NumFig04}
\end{figure}

In Figure~\ref{NumFig04}, we compare how the convergence factor $\rho$ changes with increasing numbers of nodes across various discretization schemes. 

\begin{figure}
	\centering
	\includegraphics[width=0.40\textwidth]{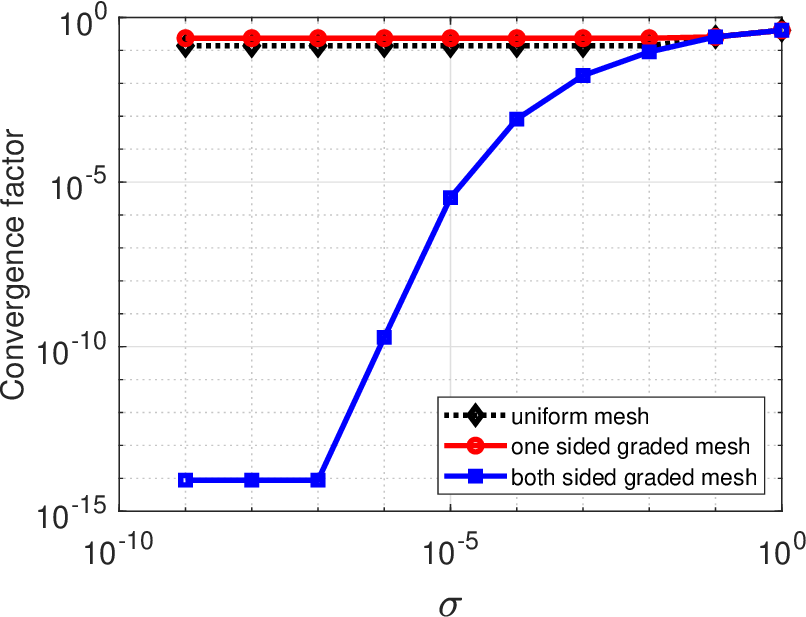}
	\includegraphics[width=0.40\textwidth]{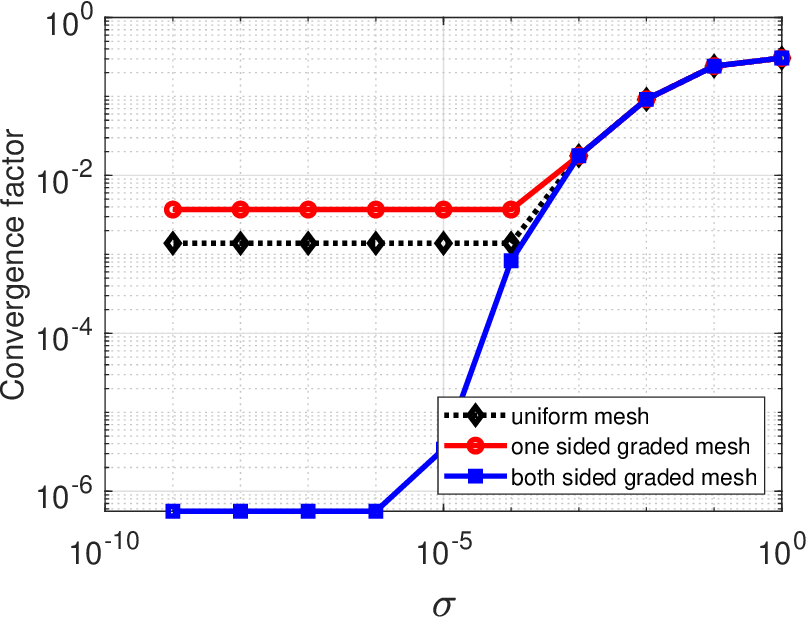}
	\caption{DNWR: As the regularization parameter $\sigma$ changes, the convergence factor varies for a constant number of nodes $N_t = 100$ and for a constant time window $T = 1$. On the left, we examine this phenomenon with a fractional order of $\alpha = 0.3$, while on the right, we consider a fractional order of $\alpha = 0.8$}
	\label{NumFig05}
\end{figure}

In Figure~\ref{NumFig05}, we analyze the relationship between the convergence factor $\rho$ and the reduction of the regularization parameter $\sigma$. This comparison reveals that the convergence rate of the DNWR algorithm increases as $\sigma$ decreases for both-sided graded mesh configurations. In contrast, for the other two cases, the convergence rate remains relatively constant for small fractional orders and saturates after reaching a certain value of $\sigma$ for large fractional orders. 

\begin{figure}
	\centering
	\includegraphics[width=0.30\textwidth]{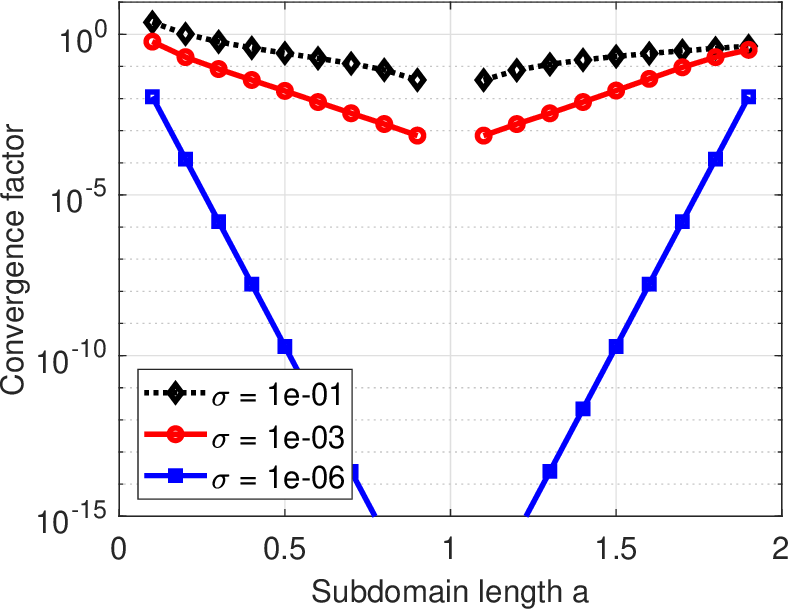}
	\includegraphics[width=0.30\textwidth]{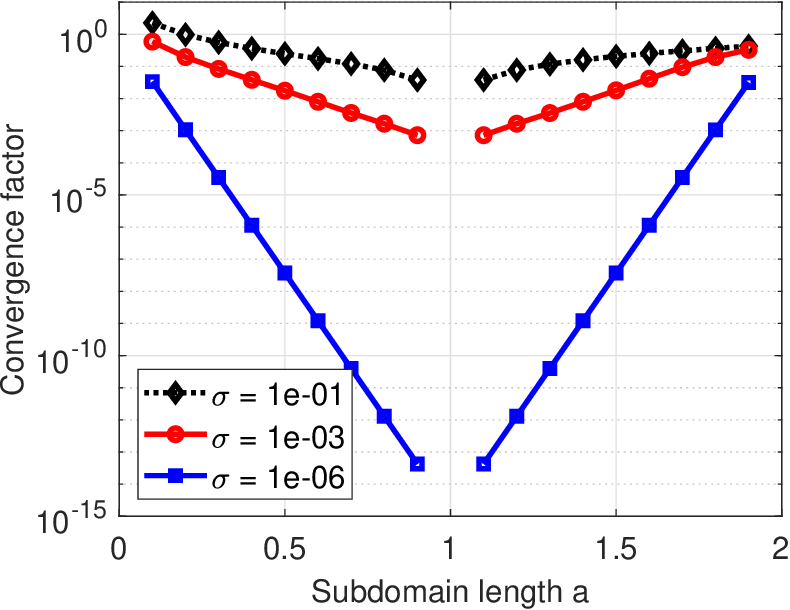}
	\includegraphics[width=0.30\textwidth]{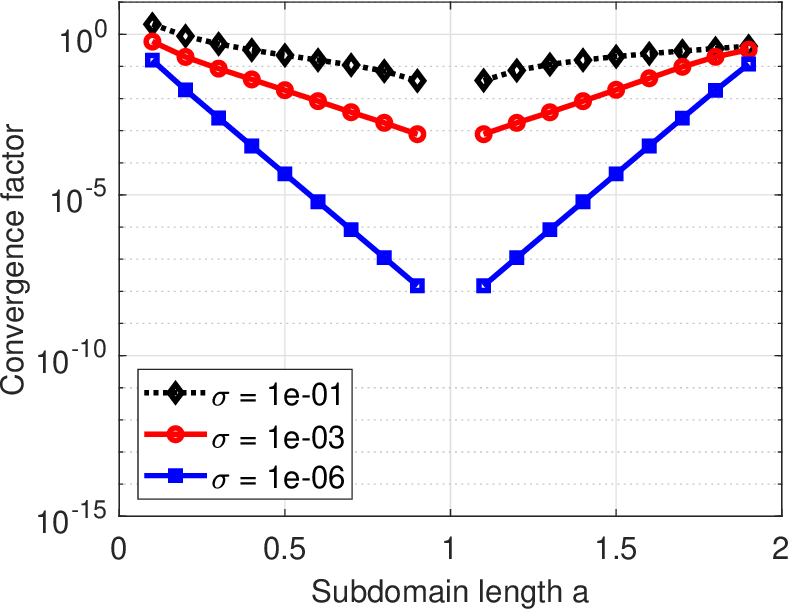}
	\caption{DNWR: As the length of the first subdomain $a$ changes and the length of the second subdomain $b$ is $b = 2-a$, with a constant number of nodes over time $N_t = 100$ and a fixed time interval $T = 1$, the convergence factor is examined. This experiment is conducted for three scenarios: on the left, with a fractional order of $\alpha = 0.3$; in the middle, with a fractional order of $\alpha = 0.7$; and on the right, with a fractional order of $\alpha = 1$.}
	\label{NumFig06}
\end{figure}

In Figure~\ref{NumFig06}, we compare the convergence factor with the length of the first subdomain. Here, we set $\kappa_{1} = \kappa_{2} = 1$, and the length of the first subdomain, $h_1$, varies from 0 to 2, with the total domain length being two. The plot indicates that for very small subdomain lengths, the convergence factor is greater than one, implying divergence of the algorithm. Additionally, at $h_1 = 1$, which corresponds to equal subdomain lengths, the algorithm exhibits an exact solution after two iterations, hence $\rho$ is not defined at this point. In this analysis, we have employed only both-sided graded mesh for time discretization.

\begin{figure}
	\centering
	\includegraphics[width=0.30\textwidth]{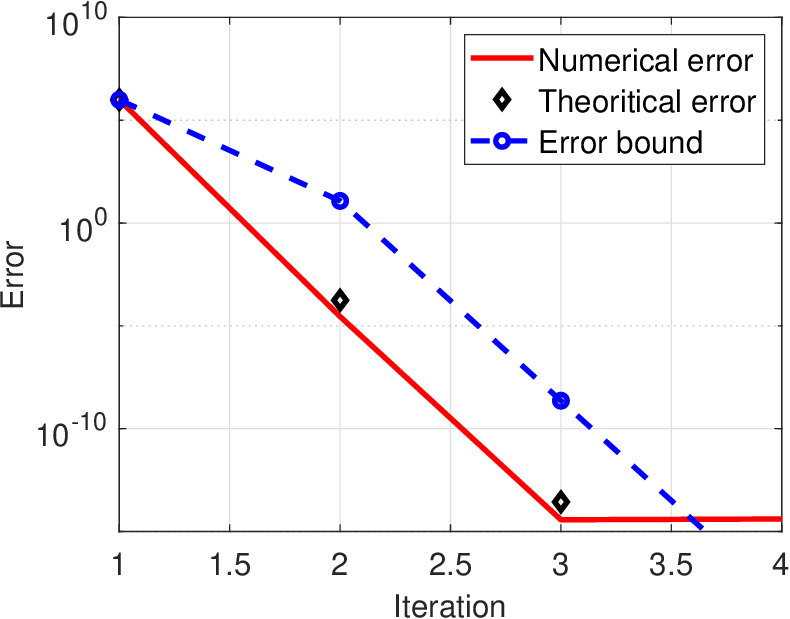}
	\includegraphics[width=0.30\textwidth]{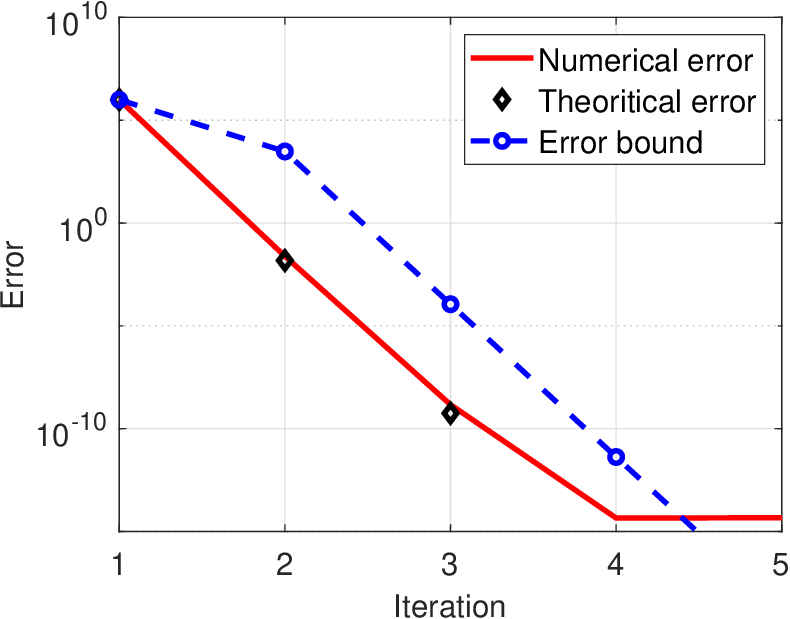}
	\includegraphics[width=0.30\textwidth]{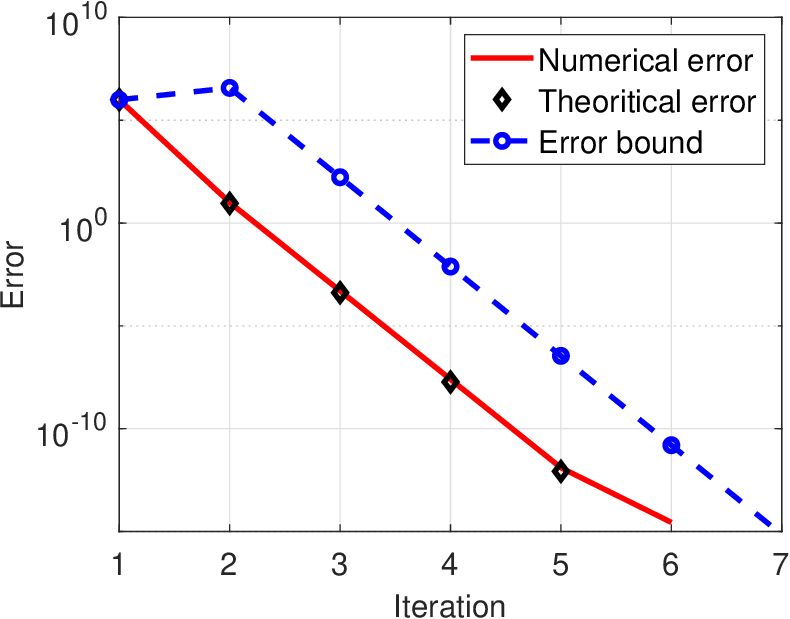}
	\caption{DNWR: Examining the convergence rates derived from numerical measurements, theoretical analysis, and estimated error bounds for a fixed number of nodes $N_t = 100$, with a constant time window $T = 1$ and regularization parameter $\sigma = 10^{-6}$. The left panel represents a fractional order of $\alpha = 0.3$, the middle panel corresponds to $\alpha = 0.7$, and the right panel depicts $\alpha = 1$.}
	\label{NumFig07}
\end{figure}

In Figure~\ref{NumFig07}, we compare the numerically measured convergence rate, theoretical convergence rate, and estimated error bounds using a both-sided graded mesh in time discretization. The experiments show that the error bound is quite sharp.  For the estimated error bounds, we calculated $\min\Re\left(\sqrt{\lambda\left(\mathbb{L}\right)}\right)$ using MATLAB's eigenvalue finder command.

\subsection{NNWR Algorithm in 1D}
For the NNWR experiments, we employ the model problem described by equations \eqref{model_problem_3}-\eqref{model_problem_4},  with zero initial condition and zero forcing term for the state equation. Depending on the experiments, we partition the spatial domain $\Omega = (-4,4)$ into multiple non-overlapping subdomains. We maintain the same target profile as specified in \eqref{num_eqn1}.

\begin{figure}
	\centering
	\includegraphics[width=0.30\textwidth]{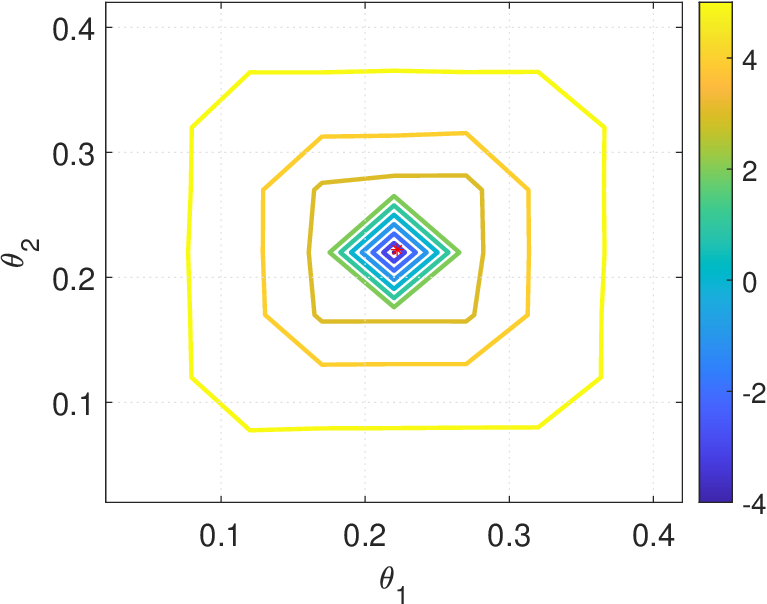}
	\includegraphics[width=0.30\textwidth]{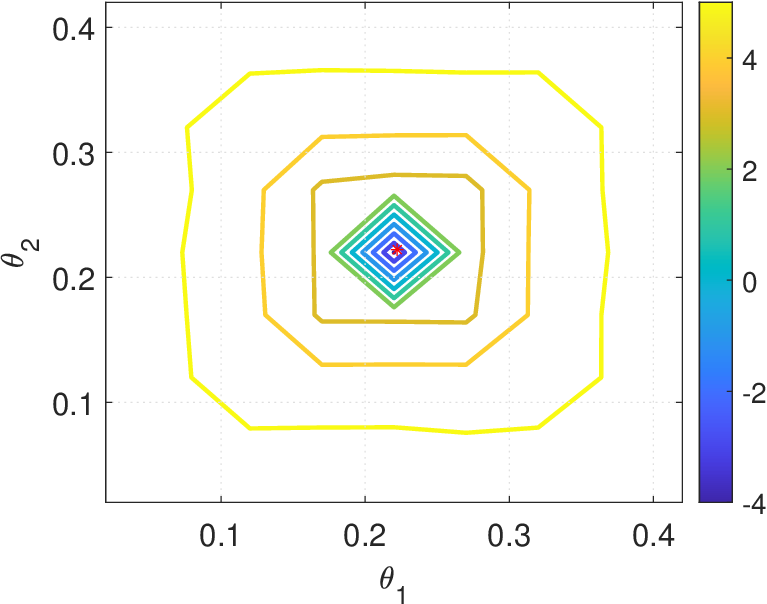}
	\includegraphics[width=0.30\textwidth]{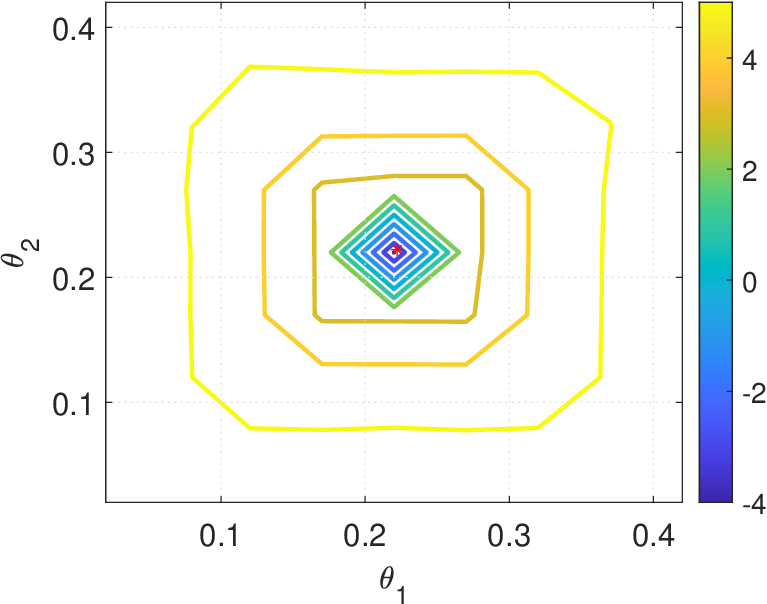}	
	\caption{NNWR: Numerical convergence rate in log scale as shown in color bar with different values of $\theta_1,\theta_2$ after five iterations at $T=1$ in a three-subdomain scenario with subdomain lengths $a_1 = 1$, $a_2 = 4$, and $a_3 = 3$, along with diffusion coefficients $\kappa_1 = 0.25$, $\kappa_2 = 1$, and $\kappa_3 = 0.25$. Specifically, the convergence rates are analyzed for fractional order on the left $\alpha=0.3$, middle $\alpha=0.7$, and right $\alpha=1$.} 
	\label{NumFig3}
\end{figure}

Prior to commencing the estimation bound experiments, we seek to validate our selected relaxation parameters for the NNWR algorithm. To achieve this, we conduct experiments using three subdomains of varying lengths: $a_1 = 1$, $a_2 = 4$, and $a_3 = 3$, each with different diffusion coefficients $\kappa_1 = 0.25$, $\kappa_2 = 1$, and $\kappa_3 = 0.25$. In our analysis, we compute $\theta_i^* = \phi_i^* = 1/(2+\sqrt{\kappa_i/\kappa_{i+1}}+\sqrt{\kappa_{i+1}/\kappa_i})$. Figure \ref{NumFig3} indicates that the values $\theta_1^* = 0.22$ and $\theta_2^* = 0.22$, marked as red stars, yield the most favorable convergence after $5$ iterations compared to all choices of $\theta_1$ and $\theta_2$ within the interval $(0,0.5)$. Additionally, these values appear to be independent of the fractional order $\alpha$. Our analysis justifies the selection of these relaxation parameters as they effectively eliminate the linear error terms. Consequently, for subsequent experiments, we consistently adopt $\theta_i^* = \phi_i^* = 1/(2+\sqrt{\kappa_i/\kappa_{i+1}}+\sqrt{\kappa_{i+1}/\kappa_i})$.

\begin{figure}
	\centering
	\includegraphics[width=0.40\textwidth]{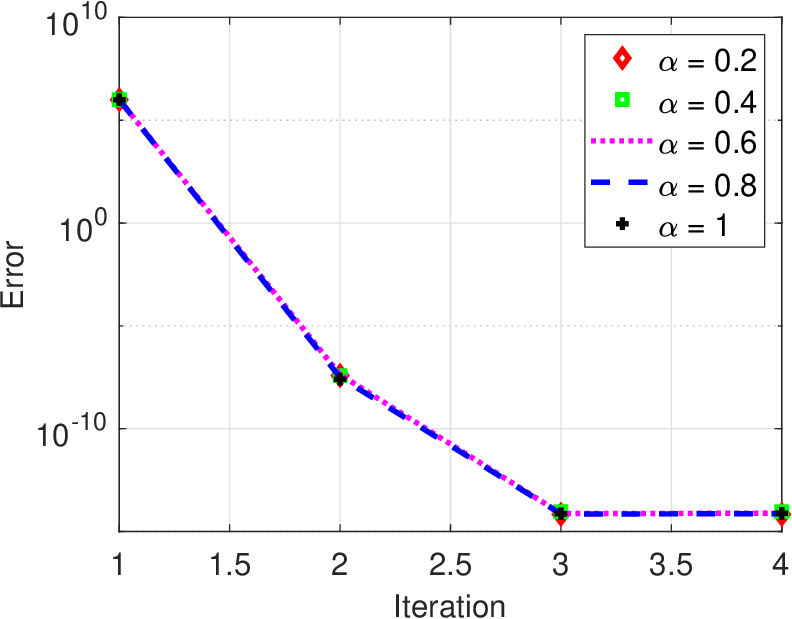}
	\includegraphics[width=0.40\textwidth]{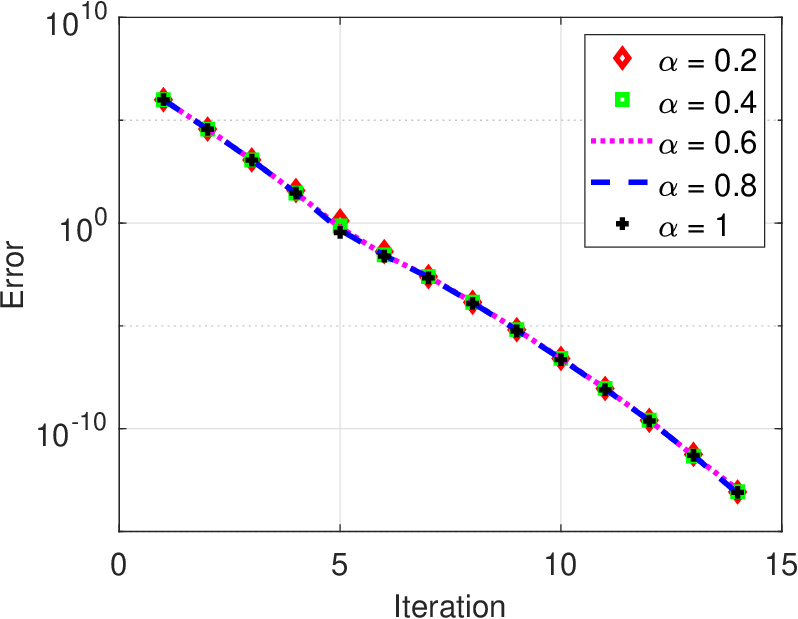}
	\caption{NNWR: The numerical convergence rates for various fractional orders $\alpha$ are investigated under two scenarios: on the left, with equal kappa and subdomain sizes; on the right, with unequal kappa and subdomain sizes. Here, $\Omega_1 = (-4,-3)$, $\Omega_2 = (-3,-1.5)$, $\Omega_3 = (-1.5,1)$, $\Omega_4 = (1,2.5)$, $\Omega_5 = (2.5,4)$, and $\kappa_1 = 0.25$, $\kappa_2 = 1$, $\kappa_3 = 0.25$, $\kappa_4 = 4$, $\kappa_5 = 1$. The fixed time is set at $T=1$.}
	\label{NumFig4}
\end{figure}

\begin{figure}
	\centering
	\includegraphics[width=0.40\textwidth]{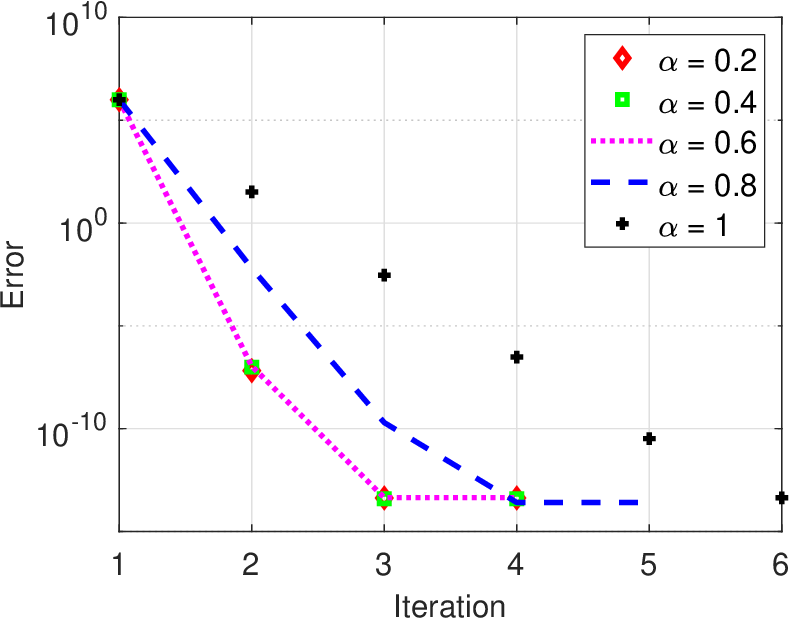}
	\includegraphics[width=0.40\textwidth]{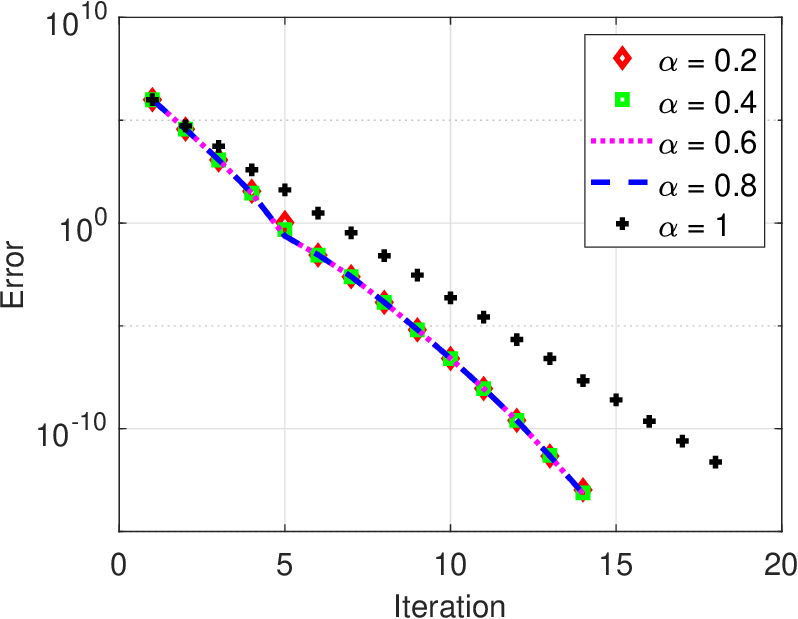}
	\caption{NNWR: The numerical convergence rates for various fractional orders $\alpha$ are investigated under two scenarios: on the left, with equal kappa and subdomain sizes; on the right, with unequal kappa and subdomain sizes. Here, $\Omega_1 = (-4,-3)$, $\Omega_2 = (-3,-1.5)$, $\Omega_3 = (-1.5,1)$, $\Omega_4 = (1,2.5)$, $\Omega_5 = (2.5,4)$, and $\kappa_1 = 0.25$, $\kappa_2 = 1$, $\kappa_3 = 0.25$, $\kappa_4 = 4$, $\kappa_5 = 1$. The fixed time is set at $T=10$.}
	\label{NumFig5}
\end{figure}

In Figure \ref{NumFig4}, we compare the convergence rates of the NNWR algorithm for the state solution across different values of the fractional order $\alpha$, considering a small time window of $T=1$. In both the cases of equal and unequal subdomains, the domain $\Omega = (-4,4)$ is partitioned into five subdomains. In the case of equal subdomains, we set constant diffusion coefficients $\kappa_i = 1$ for each subdomain. However, for the unequal subdomains, we define $\Omega_1 = (-4,-3)$, $\Omega_2 = (-3,-1.5)$, $\Omega_3 = (-1.5,1)$, $\Omega_4 = (1,2.5)$, and $\Omega_5 = (2.5,4)$, with corresponding diffusion coefficients $\kappa_1 = 0.25$, $\kappa_2 = 1$, $\kappa_3 = 0.25$, $\kappa_4 = 4$, and $\kappa_5 = 1$.
The numerical experiments suggest that the convergence rate does not depend on the fractional order, implying that the change in $\min\Re\left(\sqrt{\lambda\left(\mathbb{L}\right)}\right)$ does not affect the convergence factor significantly. In Figure \ref{NumFig5}, we repeat the same experiment but with a larger time window, namely $T = 10$.

\begin{figure}
	\centering
	\includegraphics[width=0.30\textwidth]{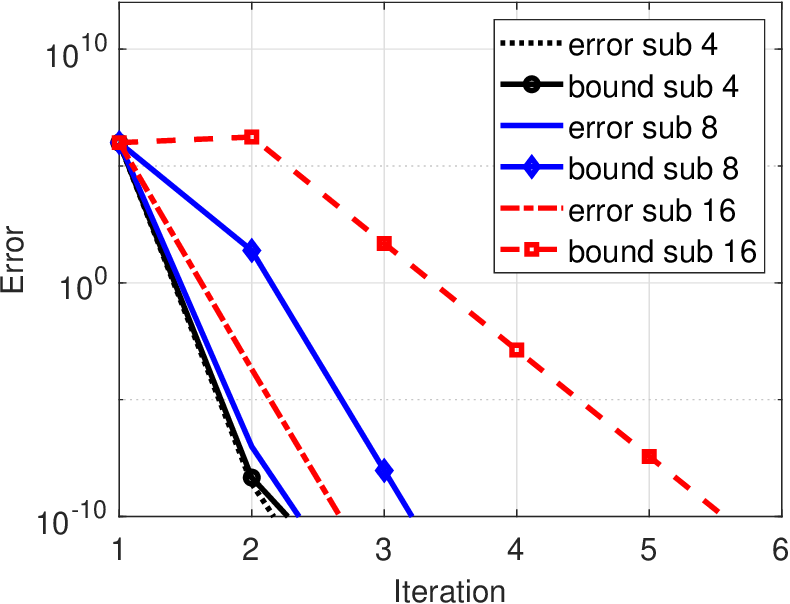}
	\includegraphics[width=0.30\textwidth]{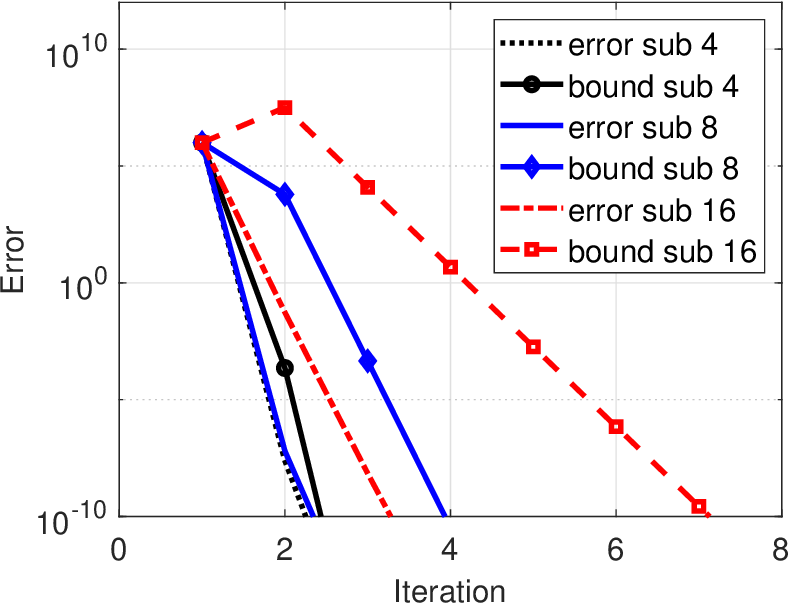}
	\includegraphics[width=0.30\textwidth]{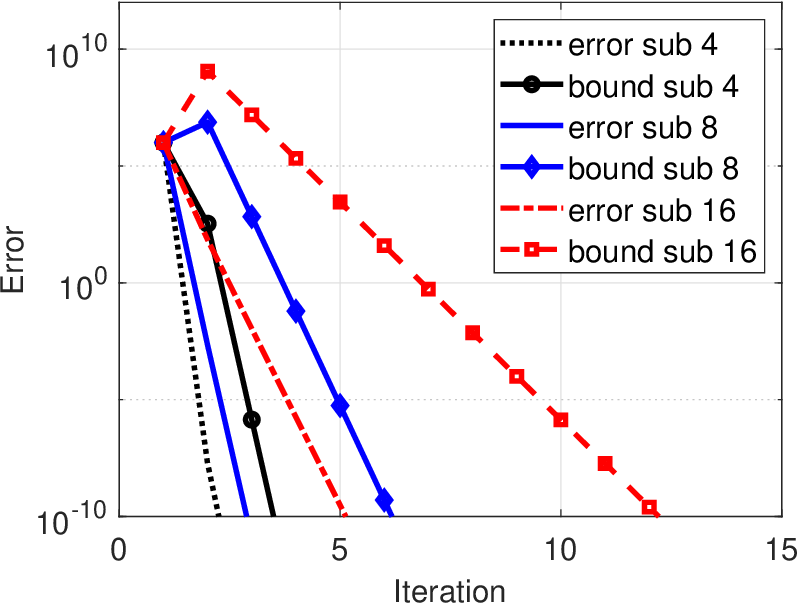}
	\caption{NNWR: Convergence rate measured numerically compared to the theoretical error bound at $T=1$, with equal subdomain lengths and a diffusion coefficient of $\kappa = 1$. On the left: $\alpha = 0.3$, in the middle: $\alpha = 0.7$, and on the right: $\alpha = 1$.}
	\label{NumFig6}
\end{figure}

In Figure \ref{NumFig6}, we compare the numerical convergence rate with the estimated bound from Theorem \ref{NNWR_th} for a small time window $T = 1$. To set up our experiments, we divide the spatial domain $\Omega = (-4,4)$ into $4$, $8$, and $16$ equal subdomains respectively with diffusion coefficients $\kappa = 1$. We observe that for a small number of subdomains, the estimate is quite sharp. However, as the number of subdomains grows, the estimated bound needs to satisfy more inequalities, as seen in \eqref{NNWR_5}, which may not be sharp in all cases, resulting in decreased sharpness.

\begin{figure}
	\centering
	\includegraphics[width=0.30\textwidth]{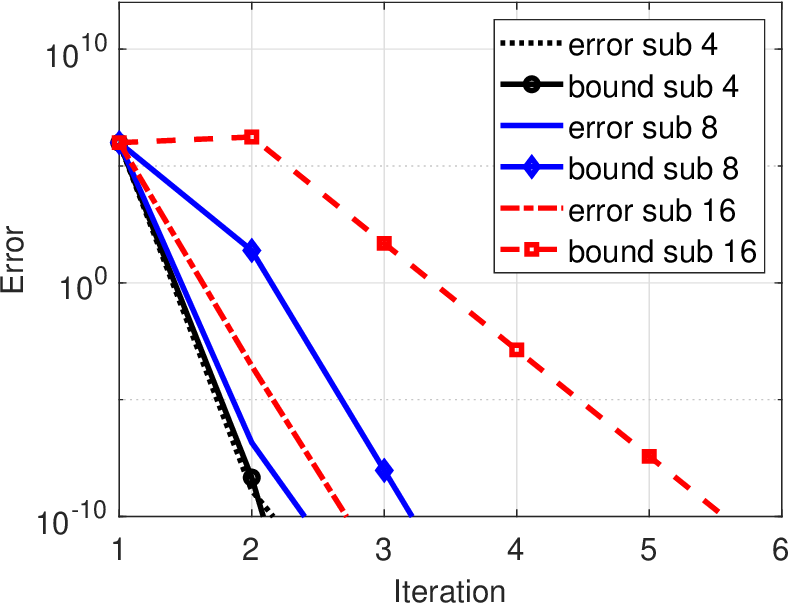}
	\includegraphics[width=0.30\textwidth]{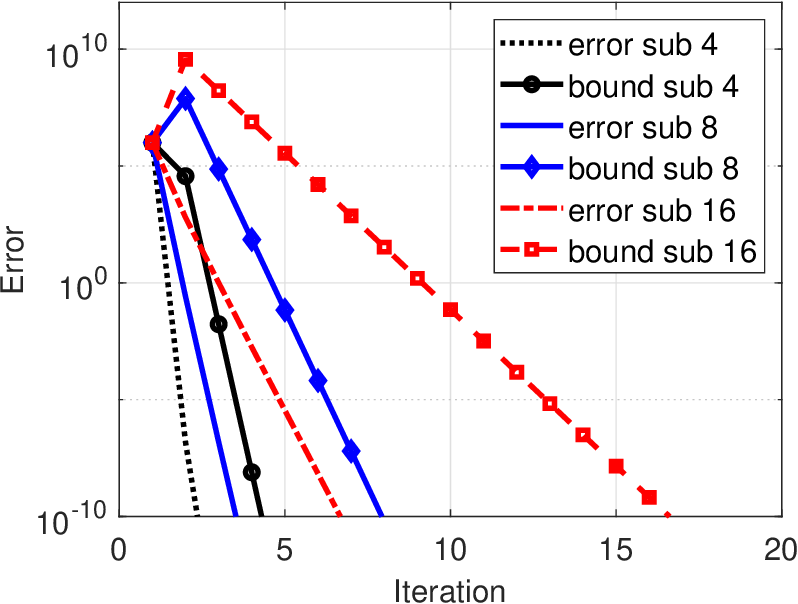}
	\includegraphics[width=0.30\textwidth]{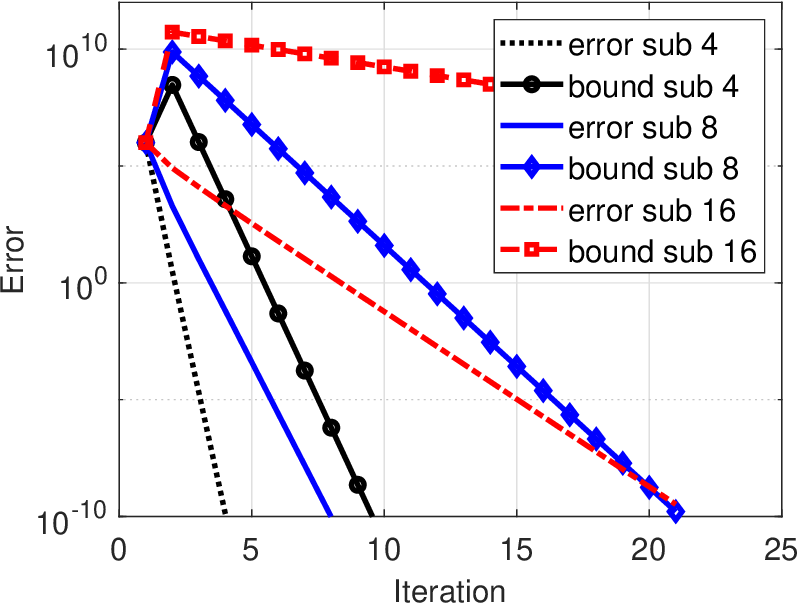}
	\caption{NNWR: Convergence rate measured numerically compared to the theoretical error bound at $T=10$, with equal subdomain lengths and a diffusion coefficient of $\kappa = 1$. On the left: $\alpha = 0.3$, in the middle: $\alpha = 0.7$, and on the right: $\alpha = 1$.}
	\label{NumFig7}
\end{figure}

In Figure \ref{NumFig7}, we conduct the same experiments but for a time window $T = 10$. Here, we also observe the same effect of subdomains on the estimated bound. The combined effect of these two factors is highly visible in the estimated bound plot for $16$ subdomains with a fractional order of $\alpha = 1$.

\begin{figure}
	\centering
	\includegraphics[width=0.30\textwidth]{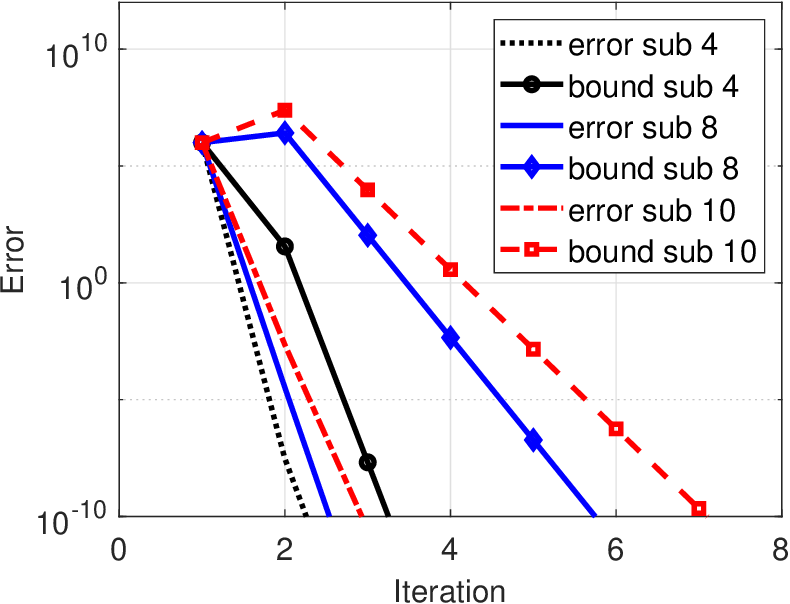}
	\includegraphics[width=0.30\textwidth]{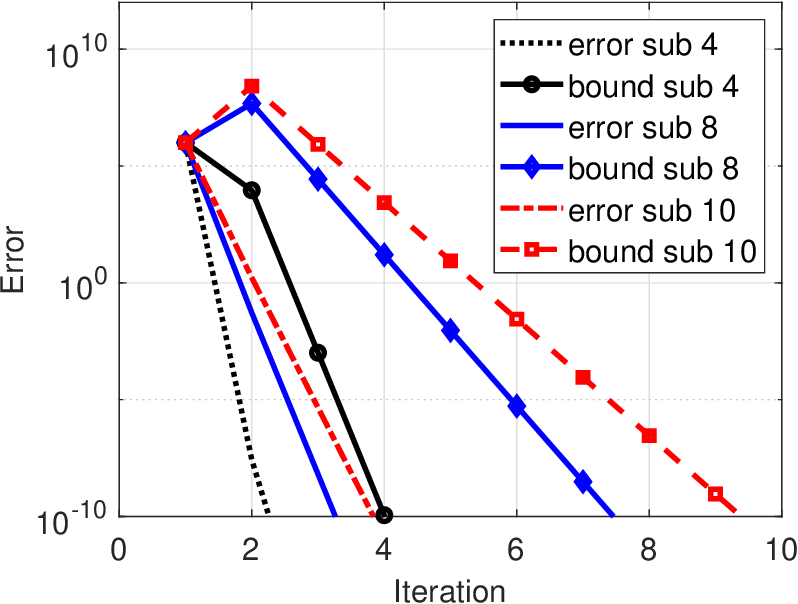}
	\includegraphics[width=0.30\textwidth]{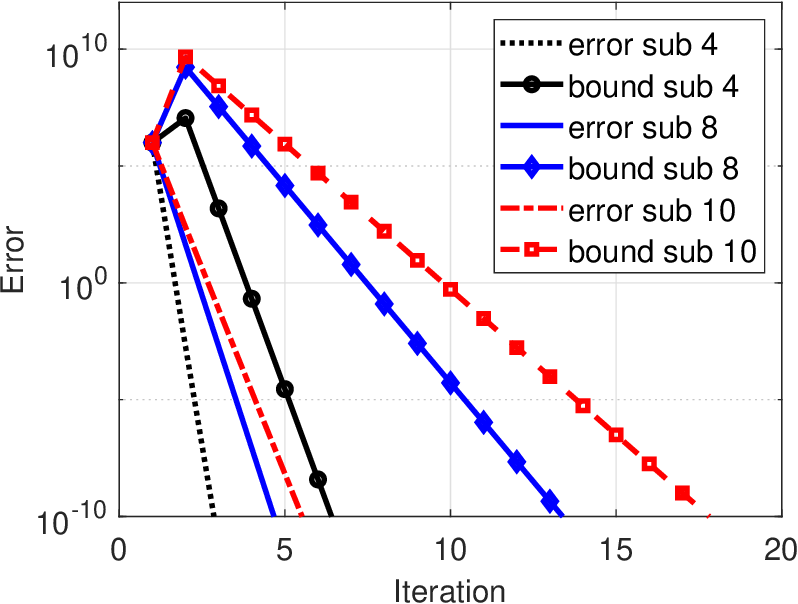}
	\caption{NNWR: Convergence rate measured numerically compared to the theoretical error bound at $T=1$, with equal subdomain lengths but different diffusion coefficients. On the left: $\alpha = 0.3$, in the middle: $\alpha = 0.7$, and on the right: $\alpha = 1$.}
	\label{NumFig8}
\end{figure}

In Figure \ref{NumFig8}, we present a comparison between the numerical convergence rate and the estimated bound derived from Theorem \ref{NNWR_th} for a small time window $T = 1$. To set up our experiments, we partition the spatial domain $\Omega = (-4,4)$ into $4, 8$, and $10$ equal subdomains, each with different diffusion coefficients $\kappa_i = \kappa_{N+1-i} = 4^{2-i}$, where $i = 1, 2, \ldots, N/2$. Here, we employ relaxation parameters $\theta_i = \phi_i = 1/(2+\sqrt{\kappa_i/\kappa_{i+1}}+\sqrt{\kappa_{i+1}/\kappa_i})$, with $i = 1, 2, \ldots, N-1$. Upon analyzing the plot, we observe a consistent resemblance in characteristics, indicating agreement with our estimates.

%% file: Conclusion.tex
\section{Conclusions}
We aim to investigate the convergence behavior of the NNWR and DNWR algorithms concerning optimal control problems constrained by sub-diffusion and diffusion partial differential equations. Our analysis is carried out in semi-discrete settings utilizing the L1 scheme in time. We propose a simple yet innovative concept (as far as we know), the both-sided graded mesh, which not only simplifies the analysis but also improves the convergence rate.


%% file: paper.bib
@book{podlubny1998fractional,
	title={Fractional differential equations: an introduction to fractional derivatives, fractional differential equations, to methods of their solution and some of their applications},
	author={Podlubny, Igor},
	year={1998},
	publisher={Elsevier}
}

@article{mophou2011optimal,
	title={Optimal control of a fractional diffusion equation with state constraints},
	author={Mophou, Gis{\`e}le M and N’Gu{\'e}r{\'e}kata, Gaston M},
	journal={Computers \& Mathematics with Applications},
	volume={62},
	number={3},
	pages={1413--1426},
	year={2011},
	publisher={Elsevier}
}

@article{mophou2011optimal_a,
	title={Optimal control of fractional diffusion equation},
	author={Mophou, Gis{\`e}le M},
	journal={Computers \& Mathematics with Applications},
	volume={61},
	number={1},
	pages={68--78},
	year={2011},
	publisher={Elsevier}
}

@article{stynes2017error,
	title={Error analysis of a finite difference method on graded meshes for a time-fractional diffusion equation},
	author={Stynes, Martin and O'Riordan, Eugene and Gracia, Jos{\'e} Luis},
	journal={SIAM Journal on Numerical Analysis},
	volume={55},
	number={2},
	pages={1057--1079},
	year={2017},
	publisher={SIAM}
}

@article{ostrowski1962some,
	title={Some theorems on the inertia of general matrices},
	author={Ostrowski, Alexander and Schneider, Hans},
	journal={Journal of Mathematical analysis and applications},
	volume={4},
	number={1},
	pages={72--84},
	year={1962},
	publisher={Academic Press}
}

@article{stynes2021survey,
	title={A survey of the L1 scheme in the discretisation of time-fractional problems},
	author={Stynes, Martin},
	journal={Submitted for publication},
	year={2021}
}

@article{ciaramella2024convergence,
	title={Convergence analysis and optimization of a Robin Schwarz waveform relaxation method for time-periodic parabolic optimal control problems},
	author={Ciaramella, Gabriele and Halpern, Laurence and Mechelli, Luca},
	journal={Journal of Computational Physics},
	volume={496},
	pages={112572},
	year={2024},
	publisher={Elsevier}
}

@article{rossikhin2010application,
	title={Application of fractional calculus for dynamic problems of solid mechanics: novel trends and recent results},
	author={Rossikhin, Yuriy A and Shitikova, Marina V},
	year={2010}
}

@article{povstenko2020fractional,
	title={Fractional thermoelasticity problem for an infinite solid with a penny-shaped crack under prescribed heat flux across its surfaces},
	author={Povstenko, Y and Kyrylych, T},
	journal={Philosophical Transactions of the Royal Society A},
	volume={378},
	number={2172},
	pages={20190289},
	year={2020},
	publisher={The Royal Society Publishing}
}

@article{qin2017multi,
	title={Multi-term time-fractional Bloch equations and application in magnetic resonance imaging},
	author={Qin, Shanlin and Liu, Fawang and Turner, Ian and Vegh, Viktor and Yu, Qiang and Yang, Qianqian},
	journal={Journal of Computational and Applied Mathematics},
	volume={319},
	pages={308--319},
	year={2017},
	publisher={Elsevier}
}

@book{schwarz1870ueber,
	title={Ueber einen Grenz{\"u}bergang durch alternirendes Verfahren},
	author={Schwarz, Hermann Amandus},
	year={1870},
	publisher={Z{\"u}rcher u. Furrer}
}

@inproceedings{lions1990schwarz,
	title={On the Schwarz alternating method. III: a variant for nonoverlapping subdomains},
	author={Lions, Pierre-Louis},
	booktitle={Third international symposium on domain decomposition methods for partial differential equations},
	volume={6},
	pages={202--223},
	year={1990},
	organization={SIAM Philadelphia}
}

@article{cai1999restricted,
	title={A restricted additive Schwarz preconditioner for general sparse linear systems},
	author={Cai, Xiao-Chuan and Sarkis, Marcus},
	journal={Siam journal on scientific computing},
	volume={21},
	number={2},
	pages={792--797},
	year={1999},
	publisher={SIAM}
}

@incollection{dryja1990some,
	title={Some domain decomposition algorithms for elliptic problems},
	author={Dryja, Maksymilian and Widlund, Olof B},
	booktitle={Iterative methods for large linear systems},
	pages={273--291},
	year={1990},
	publisher={Elsevier}
}

@article{gander2007optimized,
	title={Optimized Schwarz waveform relaxation methods for advection reaction diffusion problems},
	author={Gander, Martin J and Halpern, Laurence},
	journal={SIAM Journal on Numerical Analysis},
	volume={45},
	number={2},
	pages={666--697},
	year={2007},
	publisher={SIAM}
}

@article{gander1998space,
	title={Space-time continuous analysis of waveform relaxation for the heat equation},
	author={Gander, Martin J and Stuart, Andrew M},
	journal={SIAM Journal on Scientific Computing},
	volume={19},
	number={6},
	pages={2014--2031},
	year={1998},
	publisher={SIAM}
}

@article{giladi2002space,
	title={Space-time domain decomposition for parabolic problems},
	author={Giladi, Eldar and Keller, Herbert B},
	journal={Numerische Mathematik},
	volume={93},
	number={2},
	pages={279--313},
	year={2002},
	publisher={Springer}
}

@article{gander2021dirichlet,
	title={Dirichlet--Neumann waveform relaxation methods for parabolic and hyperbolic problems in multiple subdomains},
	author={Gander, Martin J and Kwok, Felix and Mandal, Bankim C},
	journal={BIT Numerical Mathematics},
	volume={61},
	number={1},
	pages={173--207},
	year={2021},
	publisher={Springer}
}

@article{mandal2017neumann,
	title={Neumann--Neumann waveform relaxation algorithm in multiple subdomains for hyperbolic problems in 1D and 2D},
	author={Mandal, Bankim C},
	journal={Numerical Methods for Partial Differential Equations},
	volume={33},
	number={2},
	pages={514--530},
	year={2017},
	publisher={Wiley Online Library}
}

@article{etna_vol45_pp424-456,
	author  = {Martin J. Gander and Felix Kwok and Bankim C. Mandal},
	title   = {Dirichlet-Neumann and Neumann-Neumann waveform relaxation algorithms for parabolic problems},
	journal = {Electron. Trans. Numer. Anal.},
	volume  = {45},
	year    = {2016},
	pages   = {424--456},
}

@article{boltyanskiy1962mathematical,
	title={Mathematical theory of optimal processes},
	author={Boltyanskiy, VG and Gamkrelidze, Revaz V and Mishchenko, YEF and Pontryagin, LS},
	year={1962}
}

@book{lions1971optimal,
	title={Optimal control of systems governed by partial differential equations},
	author={Lions, Jacques Louis},
	volume={170},
	year={1971},
	publisher={Springer}
}

@article{sussmann1997300,
	title={300 years of optimal control: from the brachystochrone to the maximum principle},
	author={Sussmann, Hector J and Willems, Jan C},
	journal={IEEE Control Systems Magazine},
	volume={17},
	number={3},
	pages={32--44},
	year={1997},
	publisher={IEEE}
}

@article{fernandez2003control,
	title={Control theory: History, mathematical achievements and perspectives},
	author={Fern{\'a}ndez Cara, Enrique and Zuazua Iriondo, Enrique},
	journal={Bolet{\'\i}n de la Sociedad Espa{\~n}ola de Matem{\'a}tica Aplicada, 26, 79-140.},
	year={2003},
	publisher={Sociedad Espa{\~n}ola de Matem{\'a}tica Aplicada}
}

@article{sana2023dirichlet,
	title={Dirichlet-Neumann and Neumann-Neumann waveform relaxation algorithms for heterogeneous sub-diffusion and diffusion-wave equations},
	author={Sana, Soura and Mandal, Bankim C},
	journal={Computers \& Mathematics with Applications},
	volume={150},
	pages={102--124},
	year={2023},
	publisher={Elsevier}
}

@article{sana2023dirichlet_a,
	title={Dirichlet-Neumann Waveform Relaxation Algorithm for Time Fractional Diffusion Equation in Heterogeneous Media},
	author={Sana, Soura and Mandal, Bankim C},
	journal={arXiv preprint arXiv:2301.12909},
	year={2023}
}

@article{qian2017certified,
	title={A certified trust region reduced basis approach to PDE-constrained optimization},
	author={Qian, Elizabeth and Grepl, Martin and Veroy, Karen and Willcox, Karen},
	journal={SIAM Journal on Scientific Computing},
	volume={39},
	number={5},
	pages={S434--S460},
	year={2017},
	publisher={SIAM}
}

@article{quiroga2015adjoint,
	title={Adjoint method for a tumor invasion PDE-constrained optimization problem in 2D using adaptive finite element method},
	author={Quiroga, Andr{\'e}s Agust{\'\i}n Ignacio and Fern{\'a}ndez, Dami{\'a}n and Torres, Germ{\'a}n Ariel and Turner, Cristina Vilma},
	journal={Applied Mathematics and Computation},
	volume={270},
	pages={358--368},
	year={2015},
	publisher={Elsevier}
}

@inproceedings{bensoussan1973methode,
	title={Decomposition method applied to the optimal control of distributed systems},
	author={Bensoussan, Alain and Glowinski, Roland and Lions, Jacques-Louis},
	journal={IFIP Technical Conference on Optimization Techniques},
	pages={141--151},
	year={1973},
	publisher={Springer}
}

@article{benamou1996domain,
	title={A domain decomposition method with coupled transmission conditions for the optimal control of systems governed by elliptic partial differential equations},
	author={Benamou, Jean-David},
	journal={SIAM journal on numerical analysis},
	volume={33},
	number={6},
	pages={2401--2416},
	year={1996},
	publisher={SIAM}
}

@article{benamou1997domain,
	title={A domain decomposition method for the Helmholtz equation and related optimal control problems},
	author={Benamou, Jean-David and Despr{\`e}s, Bruno},
	journal={Journal of Computational Physics},
	volume={136},
	number={1},
	pages={68--82},
	year={1997},
	publisher={Elsevier}
}
